\documentstyle[11pt,twoside,amsfonts,amssymb]{article}
  
\begin{titlepage}

\title{Towards the trace formula for convex-cocompact groups}

\author{Ulrich Bunke\thanks{Mathematisches Institut, Universit\"at G\"ottingen, Bunsenstr. 3-5, 37073 G\"ottingen, GERMANY, E-mail:bunke@uni-math.gwdg.de} and
Martin
Olbrich\thanks{Mathematisches Institut, Universit\"at G\"ottingen, Bunsenstr. 3-5, 37073 G\"ottingen, GERMANY, E-mail:olbrich@uni-math.gwdg.de  }
}

\end{titlepage}

\newcommand{\proof}{{\it Proof.$\:\:\:\:$}}
\newcommand{\dist}{{\mathrm{dist}}}
\newcommand{\kaaa}{{\frak k}}

\newcommand{\R}{{\Bbb R}}
\newcommand{\Z}{{\Bbb Z}}
\newcommand{\C}{{\Bbb C}}

\newcommand{\gaaa}{{\frak g}}

\newcommand{\aaaa}{{\frak a}}
\newcommand{\naaa}{{\frak n}}

\newcommand{\res}{{\mathrm{ res}}}

\newcommand{\cD}{{\cal D}}

\newcommand{\cO}{{\cal O}}

\newcommand{\cU}{{\cal U}}
\newcommand{\Hom}{{ \mathrm{ Hom}}}
\newcommand{\vol}{{\mathrm{ vol}}}

\newcommand{\End}{{ \mathrm{ End}}}

\newcommand{\Ree}{{\mathrm{ Re }}}

\newcommand{\ee}{{\mathrm{ e}}}

\newcommand{\tr}{{ \mathrm{ tr}}}

\newcommand{\Ad}{{\mathrm{ Ad}}}
\newcommand{\ad}{{ \mathrm{ ad}}}

\newcommand{\id}{{ \mathrm{ id}}}
\newcommand{\ord}{{ \mathrm{ ord}}}
\newcommand{\nat}{{\Bbb  N}}
\newcommand{\supp}{{ \mathrm{ supp}}}

\newcommand{\aca}{{\aaaa_\C^\ast}}

\newcommand{\cR}{{\cal R}}
\def\hB{\hspace*{\fill}$\Box$\newline\noindent}

\newcommand{\cS}{{\cal S}}

\newcommand{\Tr}{{\mathrm{Tr}}}

\newcommand{\cP}{{\cal P}}

\def\imath{i}
\newtheorem{prop}{Proposition}[section]
\newtheorem{lem}[prop]{Lemma}
\newtheorem{ddd}[prop]{Definition}
\newtheorem{theorem}[prop]{Theorem}
\newtheorem{kor}[prop]{Corollary}

\newtheorem{con}[prop]{Conjecture}

\begin{document}
\maketitle
\tableofcontents

\section{Introduction}

In this paper we develop a part of the harmonic analysis associated with  
a convex cocompact subgroup $\Gamma$ of a semisimple Lie group $G$ of real rank
one that could play the same role as the trace formula in the case of cocompact
groups or groups of finite covolume.  In these classical situations a smooth,
compactly supported, and $K$-finite  function $f$ on $G$ acts by right
convolution $R_\Gamma(f)$ on the Hilbert space $L^2(\Gamma\backslash G)$. The
trace formula is an expression of the trace of the restriction of this
operator to the discrete subspace in terms of $f$ and its Fourier transform
$\hat{f}$. The part involving $f$ is called the geometric side and usually
written as a sum of orbital integrals. The Fourier transform enters the trace
formula in the case of non-cocompact subgroups, where it is combined with the
scattering matrix. We will call this part the contribution of the scattering
matrix.

In the present paper we assume that $\Gamma$ is a convex cocompact subgroup of
$G$. Let $X$ be the symmetric space of $G$ and $\partial X$ be its
geodesic boundary. If $\Gamma\subset G$ is a discrete torsion-free
subgroup, then there is a $\Gamma$-equivariant decomposition $\partial
X=\Omega\cup \Lambda$, where $\Lambda$ is the limit set of $\Gamma$. We call
$\Gamma$ convex-cocompact if $\Gamma\backslash X\cup\Omega$ is a compact
manifold with boundary.

Since the discrete spectrum of $L^2(\Gamma\backslash G)$ is rather sparse
we take the point of view that the the contribution of the scattering matrix
is essentially (up to the contribution of the discrete spectrum) the Fourier
transform of the geometric side of the trace formula.

Thus our starting point is the geometric side. It is
a distribution $\Psi$ on $G$ given as a sum of suitably normalized orbital
integrals associated to the hyperbolic conjugacy classes of $\Gamma$ (see
Subsection \ref{psi}). The objective of the trace formula in the case
of convex cocompact $\Gamma$ is an
explicit expression for the Fourier transform of $\Psi$. We are looking for a
"measure" $\Phi$ on the unitary dual $\hat{G}$ such that
\begin{equation}\label{formula}
\Psi(f)=\int_{\hat{G}} \theta_\pi(f) \Phi(d\pi)\ ,
\end{equation}
where $\theta_\pi(f):=\Tr\:\hat{f}(\pi)$ is the character of $\pi$.
In the present paper we will formulate a precise
Conjecture \ref{zzz} about $\Phi$, but we are not able to prove
the formula (\ref{formula}) in the general case.

The unitary dual $\hat{G}$ has a natural topology. 
Now observe that the intersection of the support of
$\hat{f}$ and the support of the Plancherel measure of
$L^2(\Gamma\backslash G)$ is the spectrum of $R_\Gamma(f)$.
The  Fourier transform of a compactly supported function $f$
is never compactly supported. 
In order to do our computations
we have to approximate $R_\Gamma(f)$  
by operators which have compact spectrum.
The missing piece is some estimate which eventually  allows for dropping the
cut-off. Conjecture \ref{zzz} can easily be verified in the case of a negative
critical exponent.

In the present paper we will prove a formula which is similar to
(\ref{formula}),
but where $\Psi$ has a different interpretation. Let $R(f)$
denote the right-convolution operator on $L^2(G)$ induced by $f$. Then both,
$R(f)$ and $R_\Gamma(f)$, have smooth integral kernels $K_{R(f)}$,
$K_{R_\Gamma(f)}$, and, by Lemma \ref{asl1}, the value $\Psi(f)$ is nothing
else than the integral  $$\Psi^\prime(f):=\int_{\Gamma\backslash G}
(K_{R_\Gamma(f)}(  g,  g)-K_{R(f)}(g,g) )\mu_G(dg)\ .$$  
We will show that
$\Psi^\prime$ can be applied to functions   with compactly-supported Fourier
transform, and our main Theorem \ref{main} is a formula
\begin{equation}\label{formula1} \Psi^\prime(f)=\int_{\hat{G}} \theta_\pi(f)
\Phi(d\pi)\ . \end{equation} together with an explicit expression for $\Phi$.

We will describe $\Phi$ in terms of multiplicities $N_\Gamma(\pi)$ and a density $L_\Gamma(\pi)$. 
If $\pi$ is a representation of the complementary series (a non-tempered
unitary representation of $G$),  then the integer
$N_\Gamma(\pi):=\Phi(\{\pi\})$ is the  multiplicity of $\pi$ in
$L^2(\Gamma\backslash G)$.
If $\pi$ is a discrete series representation, then by Proposition \ref{zpro}
$N_\Gamma(\pi):=\Phi(\{\pi\})$ is an integer. It is an interesting problem to study this number in detail.

There are interesting operators with non-compact spectrum to which
$\Psi^\prime$ can be applied. Let $K\subset G$ be a maximal compact
subgroup and $\Omega$ be the Casimir operator of $G$. We fix a $K$-type
$\gamma$ and consider $(z-\Omega)^{-1}$ on $L^2(\Gamma\backslash G)(\gamma)$
and $L^2(G)(\gamma)$ if  $z$ is not in the spectrum. Let $K_\Gamma(z)$ and
$K(z)$ denote the corresponding integral kernels. The difference
$K_\Gamma(z)-K(z)$ is smooth on the diagonal and goes into $\Psi^\prime$ if
$\Ree(z)\ll 0$. The consideration of these operators provides the link between
the continuous part of $\Phi$ and the Selberg zeta functions.
In fact, using the analysis of the present paper we can show 
the meromorphic continuation of the logarithmic derivative 
and the functional equation (Theorem \ref{zetfun}) of the Selberg zeta functions
(By other methods we apriori know that the Selberg zeta functions are meromorphic, see below.). Our work extends
previous results \cite{guillope92} in the two-dimensional case and
\cite{pattersonperry95} in the spherical case of $G=SO(1,n)$, $n\ge 2$.

Using symbolic dynamics of the geodesic flow and the
thermodynamic formalism one can show that the Selberg zeta functions themselves  are meromorphic functions of finite order \cite{pattersonperry95}. This
gives information about the growth of $\Phi$ and on the counting
function of resonances. It in particular shows that $\Phi$ can be applied
to Schwartz functions like the Fourier transform $\hat{f}$ of a $K$-finite
smooth function $f$ of compact support on $G$.

\section{The distribution $\Psi$}
 
\subsection{Invariant distributions}
Let $G$ be a semisimple Lie group. We fix once and for all
a Haar measure $\mu_G$ on $G$. 
In this subsection we describe
two sorts of conjugation invariant distributions on $G$, namely
orbital integrals and characters of irreducible representations.

Let $\gamma\in G$ be a semisimple element. The orbit $\cO_\gamma:=\{g\gamma g^{-1}\:|\:g\in G\}$ 
of $\gamma$ under conjugation by $G$ is a submanifold 
of $G$ which can be identified with $G_\gamma\backslash G$, where
$G_\gamma$ denotes the centralizer of $\gamma$. The inclusion
$i_\gamma:G_\gamma\backslash G\cong \cO_\gamma\hookrightarrow G$ is a proper map. Therefore
the pull-back by $i_\gamma$ is a continuous map
$$i_\gamma^*:C^\infty_c(G)\rightarrow C^\infty_c(G_\gamma\backslash G)\ .$$
If we choose a Haar measure $\mu_{G_\gamma}$ on $G_\gamma$, then we obtain an induced measure
$\mu_{G_\gamma\backslash G}$ on $G_\gamma\backslash G$ such that
$$\int_G f(g) \mu_G(dg)=\int_{G_\gamma\backslash G}\int_{G_\gamma} f(hg) \mu_{G_\gamma}(dh) \mu_{G_\gamma\backslash G}(dg)\ .$$
The orbital integral $\theta_\gamma$
associated to $\gamma$ and the choice of the Haar measure $\mu_{G_\gamma}$
is, by definition, the composition of $i_\gamma^*$ and the measure $\mu_{G_\gamma\backslash G}$,
i.e.
$$\theta_\gamma(f):= \mu_{G_\gamma\backslash G}\circ i_\gamma^*(f)=  \int_{G_\gamma\backslash G} f(g\gamma g^{-1}) \mu_{G_\gamma\backslash G}(dg)\ .$$

We now introduce the character $\theta_\pi$ 
associated to an irreducible admissible representation $\pi$ of $G$ on a Hilbert space
$V_\pi$. If $f\in C^\infty_c(G)$, then 
$$\pi(f):=\int_G f(g) \pi(g) \mu_G(dg)$$
is a trace class operator on $V_\pi$. The character $\theta_\pi$
is the distribution on $G$ given by
$$\theta_\pi(f):= \Tr\: \pi(f)\ .$$

 \subsection{An invariant distribution asociated to $\Gamma$}\label{psi}

Let $G$ be a semisimple linear connected Lie group of real rank one.
We consider a torsion-free discrete convex-cocompact, non-cocompact
subgroup $\Gamma\subset G$ (see \cite{bunkeolbrich982}, Sec. 2).
Let $\tilde{\cO}_\Gamma$ denote the disjoint
union of manifolds $G_\gamma\backslash G$, where $\gamma$ runs over a set
$\widetilde{C\Gamma}$ of  representatives of the set $C\Gamma\setminus\{[1]\}$ of non-trivial conjugacy classes of $\Gamma$ :
$$\tilde{\cO}_\Gamma=\bigcup_{\gamma\in \widetilde{C\Gamma}}
 G_\gamma\backslash G\ .$$ 
The natural map
$i_\Gamma: \tilde{\cO}_\Gamma\rightarrow G$ is proper, and we obtain a continuous map
$$i^*_\Gamma:C^\infty_c(G)\rightarrow C^\infty_c(\tilde{\cO}_\Gamma)\ .$$
For each $\gamma\in \widetilde{C\Gamma}$ we fix a Haar measure
$\mu_{G_\gamma}$. Then we define a measure
$\mu_{\Gamma}$ on $\tilde{\cO}_\Gamma$ such that its restriction
to $G_\gamma\backslash G$ is $\vol(\Gamma_\gamma\backslash G_\gamma) \mu_{G_\gamma\backslash G}$. Note that this measure only depends
on the Haar measure $\mu_G$ and not on the choices of $\mu_{G_\gamma}$. 
\begin{ddd}
The geometric side of the trace formula is the distribution $\Psi$ on $G$ given by
$$\Psi := \mu_{\Gamma}\circ i^*_\Gamma\ .$$
\end{ddd}
 In terms of orbital integrals we can write
$$\Psi(f)=\sum_{\gamma\in \widetilde{C\Gamma}} \vol(\Gamma_\gamma \backslash G_\gamma)\theta_\gamma(f)\ .$$
Note that this distribution is in fact a measure, invariant under conjugation,
and it only depends on $\Gamma$ and the Haar measure
$\mu_G$.

\subsection{The Fourier inversion formula}

Let $\hat{G}$ denote the unitary dual of $G$. This is the set of
equivalence classes of irreducible unitary representations of $G$ equipped
with a natural structure of a measurable space. For $\pi\in\hat{G}$ 
the operator $\pi(f)$ is, by definition, the value of the Fourier tranform of $f$ at $\pi$, which we will also denote by $\hat{f}(\pi)$.

It is a consequence of the Plancherel theorem for $G$
that there is a measure $p$ on $\hat{G}$ such that for any $f\in C^\infty_c(G)$ and $g\in G$ we have
$$f(g)=\int_{\hat{G}} \Tr\: \pi(g)^{-1}\hat{f}(\pi)\: p(d\pi)\ .$$
Note that $p(d\pi)$ depends on the choice of the Haar measure $\mu_G$.
Later in the present paper we will state a more precise version of the
Plancherel theorem.  

If $h$ is a function on $\hat{G}$ such that $h(\pi)$ is a trace class operator
on $V_\pi$ for almost all $\pi$ (mod. $p$),  then we form
the function 
$$\check{h}(g):=\int_{\hat{G}} \Tr\:  \pi(g)^{-1} h(\pi)\: p(d\pi)$$
if the integral exists for all $g\in G$.

\subsection{The Fouriertransform of $\Psi$}\label{bm11}

The contents of a trace formula for convex-cocompact groups $\Gamma$
would be an expression of 
$\Psi(f)$ in terms of the  Fourier transform $\hat{f}$.
In other words, we are interested in the Fourier transform
of the distribution $\Psi$.
Since $\Psi$ is invariant this expression should only involve
the characters $\theta_\pi(f)=\Tr \hat f(\pi)$. Thus 
there should exist a certain measure
$\Phi$ on $\hat{G}$ such that the following equality holds true 
for all $f\in C^\infty_c(G)$:
$$\Psi(f)=\int_{\hat{G}} \theta_\pi(f) \Phi(d\pi) \                                                                                                                                                                                                                                                                           \ .$$
Note that there is a Paley-Wiener theorem for $G$ which characterizes
the range of the Fourier transform as a certain Paley-Wiener space. 
Apriori, $\Phi$ is a functional
on this Paley-Wiener space, and it would be  a non-trivial statement
that this functional is in fact induced by a measure on $\hat{G}$. 

\subsection{The distribution $\Psi$ as a regularized trace}

In the present paper we will not compute the Fourier transform $\Phi$ of $\Psi$
in the sense of Subsection \ref{bm11}. Rather we will compute the candidate for
$\Phi$ using a different interpretation of $\Psi$.

Let $R$ denote the right-regular representation of $G$
on $L^2(G)$. It extends 
to the convolution algebra $L^1(G)$ by the formula
$$R(f)=\int_G f(g) R(g) \mu_G(dg)\ .$$
If $f\in C^\infty_c(G)$, then $R(f)$ is an integral operator
with smooth integral kernel $K_{R(f)}(g,h)=f(g^{-1}h)$.
In a similar manner we have an unitary right-regular representation
$R_\Gamma$ of $G$ on the Hilbert space $L^2(\Gamma\backslash G)$ which can
be extended to $L^1(G)$ using the formula
$$R_\Gamma(f)=\int_G f(g) R_\Gamma(g) \mu_G(dg)\ .$$ If $f\in C^\infty_c(G)$,
then $R_\Gamma(f)$ is an integral operator with smooth kernel
$$K_{R_\Gamma(f)}(  g, h)=\sum_{\gamma\in G} f(g^{-1}\gamma h)\ .$$
Indeed, for $\phi\in L^2(\Gamma\backslash G)$ we have
\begin{eqnarray*}
R(f)\phi(  g)&=&\int_G \phi(  gh) f(h) \mu_G(dh)\\
&=&\int_G \phi( h) f(g^{-1}h) \mu_G(dh)\\
&=&\int_{\Gamma\backslash G} \sum_{\gamma\in\Gamma}\phi(  h)
f(g^{-1}\gamma h) \mu_G(dh)\ .
\end{eqnarray*}
\begin{lem}\label{asl1}
For $f\in C^\infty_c(G)$ we have
$$\Psi(f)=\int_{\Gamma\backslash G} [K_{R_\Gamma(f)}(\Gamma g,\Gamma g) - 
K_{R(f)}(g,g)] \mu_{G}(dg)\ .$$
\end{lem}
\proof
We compute
\begin{eqnarray*}
\Psi(f)&=&\sum_{ \gamma\in \widetilde{C\Gamma}} \vol(\Gamma_\gamma\backslash G_\gamma) \theta_\gamma(f)\\
&=&\sum_{\gamma\in \widetilde{C\Gamma}} \vol(\Gamma_\gamma\backslash G_\gamma)  \int_{G_\gamma\backslash G}
f(g^{-1}\gamma g)\mu_{G_\gamma\backslash G}(dg)\\
&=&
\sum_{\gamma\in \widetilde{C\Gamma}}    \int_{\Gamma_\gamma\backslash G}
f(g^{-1}\gamma g)\mu_G(dg)\\
&=&  
\int_{\Gamma \backslash G}\sum_{\gamma\in \widetilde{C\Gamma}}  \sum_{h\in \Gamma_\gamma \backslash \Gamma}
f(g^{-1} h^{-1} \gamma h g)\mu_G(dg)\\
&=&   
\int_{\Gamma \backslash G}\sum_{1\not=\gamma\in \Gamma}   
f(g^{-1}  \gamma  g)\mu_G(dg)\\
&=&\int_{\Gamma\backslash G} 
 [K_{R_\Gamma(f)}(  g,  g) - 
K_{R(f)}(g,g)]  \mu_G(dg)
\ .
\end{eqnarray*}
\hB
The expression of $\Psi(f)$ in terms of the integral kernels
of $R(f)$ and $R_\Gamma(f)$ can be used to define $\Psi$
on other classes of functions or even on certain distributions. 

Using the Plancherel theorems
for $L^2(G)$ and $L^2(\Gamma\backslash G)$ the right-regular
representations $R$ and $R_\Gamma$ can be extended.
If $f$ is $K$-finite and $\hat{f}$ is smooth and has compact support,
then we will see that $g\mapsto [K_{R_\Gamma(f)}(  g,  g) - 
K_{R(f)}(g,g)]$ belongs to $L^1(\Gamma\backslash G)$,
and thus
$$\Psi^\prime(f):=\int_{\Gamma\backslash G} 
 [K_{R_\Gamma(f)}(  g,  g) - 
K_{R(f)}(g,g)]  \mu_G(dg)$$
is well-defined. The main result of the present paper
is an expression of $\Psi^\prime(f)$ in terms of $\hat{f}$
for those functions.

As mentioned in the introduction we are going to apply
$\Psi^\prime$ to the difference of distribution kernels of the resolvents
$(z-\Omega)^{-1}$
of the Casimir operator restricted to  a $K$-type of
$L^2(\Gamma\backslash G)$ and $L^2(G)$, respectively.
In this example the single kernels are not smooth, but their difference is so
on the diagonal of $\Gamma\backslash G$. 
Strictly speaking, the integral defining $\Psi^\prime$ exists for $\Ree(z)\ll 0$.
For other values of $z$ we introduce a truncated version $\Psi^\prime_R$,
$R>0$, and we define the value of $\Psi^\prime$ as the constant term of the
asymptotic expansion of $\Psi^\prime_R$  as $R\to\infty$.
It seems to be an interesting problem to characterize the functions of
$\Omega$ (restricted to a $K$-type) with the property that
$\Psi^\prime_R$ (applied to the corresponding distribution kernels)
admits such an asymptotic expansion. 

Given a discrete series representation $\pi$ of $G$ we can consider the
corresponding isotypic components of  $L^2(G)$ and
$L^2(\Gamma\backslash G)$. If we further consider a $K$-type
of $\pi$,  then the  projections onto these
components have smooth integral kernels. As a byproduct of the investigation
of the resolvents we can show that $\Psi^\prime$ can be applied to these
integral kernels and that its values are integers.

\section{The Plancherel theorem and integral kernels}

\subsection{The Plancherel theorems for $L^2(G)$ and $L^2(\Gamma\backslash G)$.
Support of Plancherel measures}

We start with describing the rough structure of the unitary dual $\hat{G}$.
First there is a countable family of square integrable unitary representations,
the discrete series $\hat{G}_d$. 
The discretely decomposable subspace
$L^2(G)_d\subset L^2(G)$ is composed out of 
these representations each occuring with infinite multiplicity.

The orthogonal complement $L^2(G)_{ac}$ of $L^2(G)_d$ is given by a countable
direct sum of direct integrals of unitary principal series representations.
We are going to describe their parametrization.
Let $G=KAN$ be an Iwasawa decomposition of $G$. The abelian group
$A$ is isomorphic to the multiplicative group $\R^+$. Let $\aaaa$ and $\naaa$ denote the Lie algebras of $A$ and $N$.  Then $\dim_\R(\aaaa)=1$  and the roots 
of $(\aaaa,\naaa)$, fix  an order on $\aaaa$. 
Let $M=Z_K(A)$ denote the centralizer of $A$ in $K$. 
The unitary
principal series representations $\pi^{\sigma,\lambda}$ of $G$ are 
parametrized by the set $(\sigma,\lambda)\in \hat{M}\times \imath\aaaa^*$.
Let $W$ denote the Weyl group $N_K(A)/M$, where $N_K(A)$ denotes
the normalizer of $A$ in $K$. It is isomorphic to $\Z_2$,
and we can choose a representative  of the non-trivial element 
$w\in N_K(A)$ such that $w^{-1}=w$. 
One knows that $\pi^{\sigma,\lambda}$ is equivalent to $\pi^{\sigma^w,-\lambda}$, where $\sigma^w$ denotes
the Weyl conjugate representation of $\sigma$ given by
$\sigma^w(m):=\sigma(m^w)$.
For $\lambda\not=0$
the representation $\pi^{\sigma,\lambda}$ is irreducible.
If $\sigma$ is equivalent to $\sigma^w$, i.e. $\sigma$ is Weyl invariant,
then it may happen that $\pi^{\sigma,0}$ is reducible. In this case
it decomposes into a sum $\pi^{\sigma,+}\oplus \pi^{\sigma,-}$ of limits of discrete series representations.

The set of  equivalence classes of unitary representations of $G$ which we have
listed above is the set of tempered representations. We refer to Sec. 8 of
\cite{bunkeolbrich982} for a discussion of the notion of temperedness
for $L^2(G)$ and $L^2(\Gamma \backslash G)$.
 
The Plancherel theorem for $L^2(G)$ is, of course, explicitly known for a long
time \cite{harishchandra76}. The Plancherel measure $p$ is supported on  the set of tempered representations (compare \cite{bernstein88}). In particular, it is  absolutely continuous with respect to the Lebesgue measure on $\imath\aaaa$. Thus we can neglect the point $\lambda=0$. Then $L^2(G)_{ac}$ decomposes as a direct integral of
unitary principal series representations over $\hat{M}\times\imath\aaaa^*_+$
with infinite multiplicity, and the Plancherel measure has full support. Note
that the multiplicity space of the representation $\pi$ can be realized as
$V_\pi^*$.

By $\hat{G}_{ac}$ we denote the set of irreducible unitary principal series
representations $\pi^{\sigma,\lambda}$, $\lambda\not=0$.
The remaining unitary representations 
$\hat{G}_c=\hat{G}\setminus (\hat{G}_d\cup\hat{G}_{ac})$ can be realized
as subspaces of principal series representations $\pi^{\sigma,\lambda}$
with $\lambda\in\aaaa^*_+\cup\{0\}$. 
The case of limits of discrete series $\hat{G}_{ld}$ (in this case
$\lambda=0$) was mentioned above. The representations with parameter
$\lambda > 0$ are not tempered and belong to the complementary series
$\hat{G}_{cs}$.

In \cite{bunkeolbrich982} we studied the Plancherel theorem for $L^2(\Gamma\backslash G)$. Let us recall its rough structure.
The support of the corresponding Plancherel measure $p_\Gamma$ is the union
of $\hat{G}_d$, $\hat{G}_{ac}$, and a countable subset of $\hat{G}_c$.
$L^2(\Gamma\backslash G)$ decomposes into  sum of subspaces
$L^2(\Gamma\backslash G)_{cusp}$, $L^2(\Gamma\backslash G)_{ac}$,
and $L^2(\Gamma\backslash G)_c$. Here $L^2(\Gamma\backslash G)_{cusp}$
is discretely decomposable into representations of the discrete series,
each occuring with infinite multiplicity, $L^2(\Gamma\backslash G)_c$
is discretely decomposable into representations belonging to $\hat{G}_c$,
each occuring with finite multiplicity, and $L^2(\Gamma\backslash G)_{ac}$
is a direct integral of unitary principal series representations
with infinite multiplicity over the parameter set $\hat{M}\times\imath\aaaa^*_+$. On this set the Plancherel measure $p_\Gamma$
is absolutely continuous to the Lebesgue measure and has full support.
The multiplicity space $M_\pi$ can be realized as a subspace of the $\Gamma$-invariant distribution vectors of $V_{\tilde\pi}$, i.e.,
$M_\pi\subset {}^\Gamma V_{\tilde\pi,-\infty}$, where $\tilde{\pi}$ denotes the
dual representation of $\pi$. For $\pi\in \hat{G}_{ac}$ we are going to describe
$M_\pi$ explicitly in Subsection \ref{zc1}.
 
\subsection{Extension of $R$ and $R_\Gamma$}

The Plancherel theorem for $G$ provides a $G$-equivariant unitary equivalence
\begin{equation}\label{vb1}
U:L^2(G)\stackrel{\sim}{\rightarrow} \int_{\hat{G}} V_{\pi}^*\hat{\otimes} V_\pi \:p(d\pi)\ ,\end{equation}
where $G$ acts on $L^2(G)$ by the right-regular representation $R$, and
the action on the direct integral is given by $g\mapsto \{\pi\mapsto \id_{V_{\pi}^*} \otimes \pi(g)\}$.
We can identify $V_{\pi}^*\hat{\otimes} V_\pi$ with the space
of Hilbert-Schmidt operators on $V_{\pi}$.
For $\phi\in C^\infty_c(G)$ we set $$U(\phi):=\{\pi\mapsto \hat{\tilde{\phi}}(\pi)\}\ , $$
where $\tilde{\phi}(g):=\phi(g^{-1})$.
This fixes the normalization of the Plancherel measure $p$.

The inverse transformation
maps the family $\pi\mapsto h(\pi)$ to the function 
$$U^{-1}(h)(g)=\int_{\hat{G}} \Tr\:
\pi(g) h(\pi) p(d\pi)\ .$$
If $f\in L^1(G)$, than $R(f)$ is given by 
$g\mapsto \{\pi\mapsto \id_{V_{\pi}^*} \otimes \hat{f}(\pi)\}$.
 
The Plancherel theorem for $\Gamma\backslash G$ provides a
$G$-equivariant unitary equivalence 
\begin{equation}\label{vb2}
U_\Gamma:L^2(\Gamma\backslash G)\stackrel{\sim}{\rightarrow}
\int_{\hat{G}} M_{\pi}\hat{\otimes} V_\pi \:p_\Gamma(d\pi)\ ,\end{equation}
where $G$ acts on $L^2(\Gamma\backslash G)$
by the right-regular representation $R_\Gamma$,
and the representation of $G$ on the direct integral
is given by $g\mapsto \{\pi\mapsto \id_{M_{\pi}} \otimes \pi(g)\}$. Again, if  $f\in L^1(G)$, than $R_\Gamma(f)$ is given by 
$g\mapsto \{\pi\mapsto \id_{M_\pi} \otimes \hat{f}(\pi)\}$. 
In order to write down an explicit formula for $U_\Gamma$
we first identify $M_\pi^*$ with $M_{\tilde\pi}$
and embed $M_{\pi}\hat{\otimes} V_\pi$ into $\Hom(M_{\tilde\pi},V_\pi)$.
For $\phi\in C^\infty_c(\Gamma\backslash G)$ we define
$$U_\Gamma(\phi)(\pi):=\{M_{\tilde\pi}\ni v\mapsto \int_{\Gamma\backslash G}
 \phi(g) \pi(g^{-1}) v \:\mu_G(dg) \in V_\pi\}\ .$$
This fixes  the normalization of $p_\Gamma$.

Let now $h$ be a function on $\supp(p)$ such that $h(\pi)$ is a bounded
operator on $V_\pi$. If $h$ is essentially bounded (and measurable
in the appropriate sense), then it acts on the direct integral (\ref{vb1})
by $\pi\mapsto \{\id_{V_{\pi}^*}\otimes h(\pi)\}$ 
and thus defines a bounded operator $\check{R}(h)$ on $L^2(G)$
commuting with the left-regular action of $G$.

In a similar manner, if $h$ is a function on $\supp(p_\Gamma)$ such that $h(\pi)$ is a bounded operator on $V_\pi$, and $h$ is essentially bounded, then it acts on the direct integral (\ref{vb2}) 
by $\pi\mapsto \{\id_{M_{\pi}}\otimes h(\pi)\}$ 
and thus defines a bounded operator $\check{R}_\Gamma(h)$ on $L^2(\Gamma\backslash G)$.

Let us now assume that the $h(\pi)$ are of trace-class, and that
$\int_{\hat{G}} \|h(\pi)\|_{1}\: p(d\pi)$ is finite, where
$\|.\|_1$ denotes the trace norm $\|A\|_1= \Tr \:|A|$ for
a trace class operator $A$ on $V_\pi$.
Then $\check{R}(h)$ is an integral operator with integral kernel
$K_{\check{R}(h)}(g,k)= \check{h}(g^{-1}k)$.
Indeed, for $\phi\in C^\infty_c(G)$ we have
\begin{eqnarray*}
(\check{R}(h)\phi)(g)&=&\int_{\hat{G}} \Tr \:\pi(g) h(\pi)\hat{\tilde{\phi}}(\pi) p(d\pi)\\
&=&\int_{\hat{G}}\Tr \int_G \pi(g) h(\pi) \pi(k^{-1})\phi(k) \mu_G(dk) p(d\pi)\\
&=&\int_G \int_{\hat{G}}\Tr \:\pi(k^{-1} g ) h(\pi) p(d\pi)\phi(k) \mu_G(dk) \\
&=&\int_G \check{h}(g^{-1}k) \phi(k) \mu_G(dk)\ . 
\end{eqnarray*}
We are looking for a similar formula for the integral kernel of
$\check{R}_\Gamma(h)$ in Subsection \ref{zc1}.

\subsection{The absolute continuous part of $L^2(\Gamma\backslash G)$. Integral kernels for $\check{R}_\Gamma(h)$.}\label{zc1}

In this subsection we describe in detail the Plancherel decomposition
of $L^2(\Gamma\backslash G)_{ac}$. The goal is to exhibit a class
of functions $\pi\mapsto h(\pi)$ with the property that
$\check{R}_\Gamma(h)$ is an integral operator.

Let $P=MAN$ be a fixed parabolic subgroup. If $(\sigma,\lambda)\in\hat{M}\times \aca$, then we define the representation $\sigma_\lambda$
of $P$ by $\sigma_\lambda(man):=\sigma(m) a^{\rho-\lambda}$,
where $\rho\in\aaaa^*$ is given by $2\rho(H)=\tr\:\ad(H)_{|\naaa}$, $H\in\aaaa$,
and for $\lambda\in\aca$ and $a=\exp(H)\in A$ we put $a^\lambda=\ee^{\lambda(H)}$.
We realize the principal series representation
$H^{\sigma,\lambda}:=V_{\pi^{\sigma,\lambda}}$ as the  subspace of
$C^{-\infty}(G\times_P V_{\sigma_\lambda})$ of functions 
such that $f_{|K}\in L^2(K\times_MV_\sigma)$. 
Then $H^{\sigma,\lambda}_{\pm\infty}=C^{\pm\infty}(G\times_P V_{\sigma_\lambda})$
are the spaces of smooth (resp. distribution) vectors of $\pi^{\sigma,\lambda}$.
By restriction to $K$ we obtain canonical isomorphisms
$H^{\sigma,\lambda}\cong L^2(K\times_MV_\sigma)$.
It therefore makes sense to speak of smooth functions $f$ on 
$\imath\aaaa^*$ such that $f(\lambda)\in H^{\sigma,\lambda}$.

Note that $\partial X$ can be identified with $G/P$.
Let $G/P=\Omega\cup\Lambda$ be the $\Gamma$-equivariant decomposition
of the space $G/P$ into the (open) domain of discontinuity $\Omega$ and the (closed) limit set
$\Lambda$. As a convex-cocompct subgroup
$\Gamma$ acts freely and cocompactly on $\Omega$. We put
$B:=\Gamma\backslash\Omega$. Furthermore, we define the bundle
$V_B(\sigma_\lambda)\rightarrow B$ by $V_B(\sigma_\lambda):=\Gamma\backslash (G\times_P V_{\sigma_\lambda})_{|\Omega}$. If $\lambda\in \imath\aaaa^*$, then we have
a natural Hilbert space $L^2(B,V_B(\sigma_\lambda))$. Again, fixing a volume form on $B$ we obtain identifications of the spaces 
$L^2(B,V_B(\sigma_\lambda))$ with the fixed space  $L^2(B,V_B(\sigma_0))$ so that it makes sense to speak of smooth functions $f$ on $\imath\aaaa^*$ such that $f(\lambda)\in  L^2(B,V_B(\sigma_\lambda))$. We refer to \cite{bunkeolbrich982}, Sec.3, for more details.

In \cite{bunkeolbrich982} we defined a family of extension maps
$ext: L^2(B,V_B(\sigma_\lambda))\rightarrow {}^\Gamma H^{\sigma,\lambda}_{-\infty}$. For $\lambda\in\imath\aaaa^*$, $\lambda\not=0$, the extension map provides an explicit identification of the
space of multiplicities $M_{\pi^{\sigma,\lambda}}\subset {}^\Gamma H^{\tilde\sigma,-\lambda}_{-\infty}$ with $L^2(B,V_B(\tilde{\sigma}_{-\lambda}))$.
 
The Plancherel measures $p$ and $p_\Gamma$ on $\{\sigma\}\times\imath\aaaa^*_+$ are given
by $\frac{\dim(V_\sigma)}{2\pi\omega_X}p_\sigma(\lambda)d\lambda$, where $p_\sigma$ is a smooth symmetric function on $\imath\aaaa^*$ of polynomial growth (see \cite{bunkeolbrich982}, Lemma 5.5. (3)), and  $\omega_X:=\lim_{a\to\infty} a^{-2\rho} \vol_{G/K}(KaK)$
(see \cite{bunkeolbrich982}, Sec. 11). Note that $d\lambda$
is the {\em real} Lebesgue measure on $\imath\aaaa$.

We now describe the embedding 
$$U_\Gamma^{-1}:\frac{\dim(V_\sigma)}{2\pi\omega_X}\int_{\{\sigma\}\times\imath\aaaa^*_+} L^2(B,V_B(\tilde{\sigma}_{-\lambda}))\otimes H^{\sigma,\lambda} p_\sigma(\lambda) d\lambda \rightarrow  L^2(\Gamma\backslash G)_{ac}\ .$$

If $v\otimes w\in L^2(B,V_B(\tilde{\sigma}_{-\lambda}))\otimes H^{\sigma,\lambda}_{\infty}$, then we define
$<v\otimes w>:=\langle ext(v),w\rangle$.
Let $\phi$ be a smooth function of compact support on $\imath\aaaa^*_+\cup\{0\}$
such that $\phi(\lambda)\in L^2(B,V_B(\tilde{\sigma}_{-\lambda}))\otimes 
H^{\sigma,\lambda}_{\infty}$, then we have
$$U_\Gamma^{-1}(\phi)(g)=\frac{\dim(V_\sigma)}{2\pi\omega_X}\int_{\imath\aaaa^*_+}
<(1\otimes \pi^{\sigma,\lambda}(g))\phi(\lambda)> p_\sigma(\lambda) d\lambda\ .$$
Note that $ext$ may be singular at $\lambda=0$. In this case it has
a first-order pole and $p_\sigma(0)=0$
(see \cite{bunkeolbrich982}, Prop. 7.4) such that the integral is still
well-defined.

We now fix a $K$-type $\gamma\in\hat{K}$.
Let $H^{\sigma,\lambda}(\gamma)$ denote the $\gamma$-isotypic component
of $H^{\sigma,\lambda}$. 
By Frobenius reciprocity we have a canonical identification
\begin{equation}\label{fr1}
H^{\sigma,\lambda}(\gamma)\stackrel{\alpha}{=} V_\gamma \otimes \Hom_K(V_\gamma,H^{\sigma,\lambda})\stackrel{1\otimes \beta}{=} V_\gamma\otimes \Hom_M(V_\gamma,V_\sigma)\ .\end{equation}
Here $\alpha^{-1}(v\otimes U):=U(v)$ and $\beta(U)(v):=U(v)(1)$.
Any operator $A\in \End(V_\gamma\otimes \Hom_M(V_\gamma,V_\sigma))$ gives rise
to a finite-dimensional operator $F(A)\in \End(H^{\sigma,\lambda})$
which is trivial on the orthogonal complement of $H^{\sigma,\lambda}(\gamma)$.

Let $q$ be a smooth function of compact support on $\hat{M}\times \imath\aaaa^*$ such that
$q(\sigma,\lambda)\in  \End(V_\gamma\otimes \Hom_M(V_\gamma,V_\sigma))$. 
 We call $q$ symmetric if it is compatible with the equivalences
$J^w_{\sigma,\lambda}:H^{\sigma,\lambda}\rightarrow H^{\sigma^w,-\lambda}$,
i.e. if $F(q(\sigma^w,-\lambda))=J^w_{\sigma,\lambda}\circ F(q(\sigma ,\lambda))\circ (J^w_{\sigma,\lambda})^{-1}$.
If $q$ is symmetric, then we can define the function $h_q$ on $\hat{G}$
such that $h_q(\pi)\in\End(V_\pi)$ by
$h_q(\pi^{\sigma,\lambda}):=F(q(\sigma,\lambda))$
for $(\sigma,\lambda)\in\hat{M}\times \imath\aaaa^*$, $\lambda\not=0$, and by $h_q(\pi)=0$
for all other representations. 

Let $\pi_*:H^{\sigma,\lambda}_\infty\rightarrow L^2(B,V_B(\sigma_\lambda))$
denote the push-down map which can be considered here as the adjoint of the extension $ext:L^2(B,V_B(\tilde{\sigma}_{-\lambda}))\rightarrow H^{\tilde\sigma,-\lambda}_{-\infty}$. The composition $$\pi^{\sigma,\lambda}(g) h_q(\pi^{\sigma,\lambda}) \pi^{\sigma,\lambda}(k^{-1}) ext\: \pi_*$$ is a finite-dimensional map
from $H^{\sigma,\lambda}_\infty$ to $H^{\sigma,\lambda}_\infty$.
It is therefore nuclear and has a well-defined trace.
\begin{lem}
The operator $\check{R}_\Gamma(h_q)$ has a smooth integral kernel
given by 
$$K_{\check{R}_\Gamma(h_q)}(  g,  k)= \sum_{\sigma\in\hat{M}}\frac{\dim(V_\sigma)}{4\pi\omega_X}\int_{   \imath\aaaa^*}
\Tr\:\pi^{\sigma,\lambda}(g) h_q(\pi^{\sigma,\lambda}) \pi^{\sigma,\lambda}(k^{-1}) ext\circ  \pi_*\: p_\sigma(\lambda) d\lambda\ .$$
\end{lem}
\proof
First of all note that the integral is well defined at $\lambda=0$.
If $ext\circ\pi_*$ is singular at this point, then it has
a pole of at most second order. But then the Plancherel density
vanishes at least of second order, too.

Let $\phi\in C_c^\infty(\Gamma\backslash G)$. In the Plancherel decomposition 
it is represented by the function $\pi\mapsto U_\Gamma(\phi)(\pi)\in M_\pi \hat{\otimes} V_\pi$.  
We fix $\lambda$ for a moment and choose orthonormal bases $\{v_i\}$ of
$L^2(B,V_B(\tilde{\sigma}_{-\lambda}))$ and $\{w_i\}$ of $H^{\sigma,\lambda}$
consisting of smooth sections. Furthermore, let
$\{v^i\}$ and $\{w^j\}$ be dual bases of 
$L^2(B,V_B(\sigma_{\lambda}))$ and $H^{\tilde{\sigma},-\lambda}$, respectively.
Then we have
$$U_\Gamma(\phi)(\pi^{\sigma,\lambda})=
\sum_{i,j}\int_{\Gamma\backslash G} <v^i\otimes \pi^{\tilde{\sigma},-\lambda}(k)w^j> \phi(k) \mu_G(dk)\: v_i\otimes w_j\ .$$
We have
\begin{equation}\label{nm1}
\check{R}_\Gamma(h_q)(\phi)(g)=\frac{\dim(V_\sigma)}{4\pi\omega_X} \sum_{\sigma\in\hat{M}}\int_{\imath\aaaa^*}
<(1\otimes \pi^{\sigma,\lambda}(g)h_q(\pi^{\sigma,\lambda}))U_\Gamma(\phi)(\pi^{\sigma,\lambda})>p_\sigma(\lambda) d\lambda\ .\end{equation}
We now compute
\begin{eqnarray*}
\lefteqn{
<(1\otimes \pi^{\sigma,\lambda}(g)h_q(\pi^{\sigma,\lambda}))U_\Gamma(\phi)(\pi^{\sigma,\lambda})>}&&\\
&=&\sum_{i,j}\int_{G} <v^i\otimes \pi^{\tilde{\sigma},-\lambda}(k)w^j> \phi(k) \mu_G(dk)<(1\otimes \pi^{\sigma,\lambda}(g)h_q(\pi^{\sigma,\lambda}))v_i\otimes w_j>\\
&=&\int_{G} \sum_{i,j}  <v^i\otimes  w^j> \phi(k) <(1\otimes \pi^{\sigma,\lambda}(g)h_q(\pi^{\sigma,\lambda})\pi^{\sigma,\lambda}(k^{-1}))v_i\otimes w_j>\mu_G(dk)\\
&=&\int_{G} \sum_{i,j} \phi(k) \langle v^i,\pi_* w^j\rangle
\langle ext(v_i), \pi^{\sigma,\lambda}(g)h_q(\pi^{\sigma,\lambda})\pi^{\sigma,\lambda}(k^{-1}) w_j\rangle\mu_G(dk)\\
&=&\int_{G} \sum_{j} \phi(k)  
\langle ext\circ \pi_*(w^j), \pi^{\sigma,\lambda}(g)h_q(\pi^{\sigma,\lambda})\pi^{\sigma,\lambda}(k^{-1}) w_j\rangle\mu_G(dk)\\
&=&\int_{G} \sum_{j,l} \phi(k)  \langle ext\circ \pi_*(w^j),w_l\rangle
\langle w^l, \pi^{\sigma,\lambda}(g)h_q(\pi^{\sigma,\lambda})\pi^{\sigma,\lambda}(k^{-1}) w_j\rangle\mu_G(dk)\\
&=&\int_{G} \sum_{l} \phi(k)  
\langle  w^l, \pi^{\sigma,\lambda}(g)h_q(\pi^{\sigma,\lambda})\pi^{\sigma,\lambda}(k^{-1})ext \circ \pi_* w_l\rangle\mu_G(dk)\\
&=&\int_{G}  \phi(k)  
 \Tr\: \pi^{\sigma,\lambda}(g)h_q(\pi^{\sigma,\lambda})\pi^{\sigma,\lambda}(k^{-1})ext \circ \pi_* \: \:\mu_G(dk)\ .
\end{eqnarray*}
Inserting this computation into (\ref{nm1}), the  we obtain the desired formula
for the integral kernel of $\check{R}_\Gamma(h_q)$.
\hB

\section{Poisson transforms and asymptotic computations}

\subsection{Motivation}

Let $q$ be symmetric and define and $h_q$ as in subsection \ref{zc1}. 
We want to show that
the function $g\mapsto [K_{\check{R}_\Gamma(h_q)}(  g,  g)-K_{\check{R}(h_q)}(g,g)]$ belongs to $L^1(\Gamma\backslash G)$.
It follows that
$$\Psi^\prime(\tilde{\check{h_q}}) =\int_{\Gamma\backslash G} 
[K_{\check{R}_\Gamma(h_q)}(\Gamma g,\Gamma g)-K_{\check{R}(h_q)}(g,g)] \mu_G(dg)$$
is well-defined, and we are asking for
an expression of $\Psi^\prime(\tilde{\check{h_q}})$ in terms
of $q$, respectively $h_q$. 

In the present section we show related results using the language of
Poisson transformations. In Subsection  \ref{trans} we will provide the translation
of these results and solve the problems above.  

\subsection{Poisson transformation, $c$-functions, and asymptotics}

We fix a $K$-type $\gamma$ and a $M$-type $\sigma$.
Let $T\in \Hom_{M}(V_\sigma,V_\gamma)$ and $\lambda\in\aca$.
If $w\in V_{\tilde{\gamma}}$, then by Frobenius reciprocity  we consider
$w\otimes T^*$ as an element of $H^{\tilde{\sigma},-\lambda}(\tilde{\gamma})$
which is given by the function $k\mapsto T^*(\tilde{\gamma}(k^{-1})w)$ under the canonical identification
$\phi_{-\lambda}:H^{\tilde{\sigma},-\lambda}_\infty\stackrel{\sim}{\rightarrow} C^\infty(K\times_M V_{\tilde{\sigma}})$. We further put $\Phi_{\lambda,\mu}:=\phi_{\lambda}^{-1}\circ \phi_{\mu}$.
We will also use the notation $\Phi_{0,\lambda}$ for $\phi_\lambda$.

The Poisson transformation
$$P^T_\lambda:H^{\sigma,\lambda}_{-\infty}\rightarrow C^\infty(G\times_K V_\gamma)$$ is a $G$-equivariant 
map which is defined by the relation
$$\langle w,P^T_\lambda(\psi)(g)\rangle = \langle  w\otimes T^*, \pi^{\sigma,\lambda}(g^{-1}) \psi\rangle\ ,$$
for all $\psi\in H^{\sigma,\lambda}_{-\infty}$, $w\in V_{\tilde{\gamma}}$. 

For the definition of the function
$c_\sigma$  we refer to \cite{bunkeolbrich982}, Sec. 5. 
We have the relation  $$c_\sigma(\lambda)c_{\tilde{\sigma}}(-\lambda)=p_\sigma(\lambda)^{-1}\ .$$
It turns out to be useful to introduce the normalized Poisson transformation  ${}^0P^T_\lambda:= c_\sigma(-\lambda)^{-1} P^T_\lambda$.

We introduce the family of operators $$\cP^T_{\lambda,a}:H^{\sigma,\lambda}_{-\infty}\rightarrow  
C^\infty(K\times_M V_\gamma)\ , \quad a\in A_+$$  by
$$\cP^T_{\lambda,a}(f)(k):={}^0P^T_\lambda(f)(ka)\ .$$

In order to discuss the asymptotic behaviour of $\cP^T_{\lambda,a}$
as $a\to\infty$ we need the normalized Knapp-Stein intertwining operators
$$J^w_{\sigma,\lambda}:H^{\sigma,\lambda}_{-\infty}\rightarrow
H^{\sigma^w,-\lambda}_{-\infty}\ .$$ 
Note that $J^w_{\sigma,\lambda}\circ J^w_{\sigma^w,-\lambda}=\id$.
We again refer to \cite{bunkeolbrich982}, Sec. 5, for more details.
The following is a consequence of \cite{bunkeolbrich982}, Lemma 6.2.
Let $\alpha\in\aaaa^*$ denote the short root of $(\aaaa,\naaa)$.
For $\lambda\in\imath\aaaa^*$ we have
 \begin{equation}\label{aass}\cP^{T}_{\lambda,a}=a^{\lambda-\rho} \frac{c_\gamma(\lambda)}{c_\sigma(-\lambda)}
T \circ \Phi_{0,\lambda}+ a^{-\lambda-\rho}\gamma(w) T \circ \Phi_{0,-\lambda}\circ J^w_{\sigma,\lambda}
+ a^{-\rho-\alpha} \cR_{-\infty}(\lambda,a)\ ,\end{equation}
where $\cR_{-\infty}(\lambda,a)\circ \Phi_{\lambda,0}$ remains bounded in
$C^\infty(\imath\aaaa^*,\Hom(H^{\sigma,0}_{-\infty},C^{-\infty}(K\times_M V_\gamma)))$ as $a\to\infty$. Multiplication by $T$
 (resp. $\gamma(w) T$)
is here considered as a map from $H^{\sigma,0}_{-\infty}$ (resp. $H^{\sigma^w,0}_{-\infty}$) to
$C^{-\infty}(K\times_M V_\gamma)$ in the natural way.
If $\chi,\tilde{\chi}$ are smooth functions on $K/M$ with disjoint support, then
$\chi\cR_{-\infty}(\lambda,a)\tilde{\chi}\circ \Phi_{\lambda,0}$ remains bounded in $C^\infty(\imath\aaaa^*,\Hom(H^{\sigma,0}_{-\infty},C^{\infty}(K\times_M V_\gamma)))$ as $a\to\infty$.

\subsection{An estimate}\label{esi}
 
In order to formulate the result appropriately we introduce
the following  space $C_{\Gamma}(G)$ of functions on $G$.
For each compact $V\subset \Omega$ and integer $N$ we consider the seminorm
$$|\phi|_{V,N}:=\sup_{kah\in VA_+M} (1+|\log(a)|)^{N} a^{2\rho} |\phi(kah)|\ , \phi\in C(G)\ .$$  
Here we consider $V$ as a subset of $K$ using the identification $G/P=K/M$.
We define the Fr\'echet space $C_{\Gamma}(G)$
as the space of all continuous functions $\phi$ on $G$ such that
$|\phi|_{V,N}<\infty$ for all compact $V\subset \Omega$ and $N\in\nat$.
If $\phi$ is $\Gamma$-invariant and belongs to
$C_{\Gamma}(G)$, then clearly $\phi\in L^1(\Gamma\backslash G)$.

Now let $T\in \Hom_{M}(V_\sigma,V_\gamma)$, $R\in \Hom_{M}(V_{\tilde{\sigma}},V_{\tilde{\gamma}})$, and 
$q\in C^\infty_c(\imath\aaaa^*)$. Then we can define
the operator
$$A_q= A_q(T,R):=\int_{\imath\aaaa^*} {}^0P^T_\lambda\circ \left(ext\circ \pi_* -1\right) \circ
({}^0P^R_{-\lambda})^*\:q(\lambda) d\lambda \in
\Hom(C_c^{-\infty}(G\times_KV_\gamma),C^\infty(G\times_KV_\gamma))\ .$$ This
operator has a smooth integral kernel $(g,h)\mapsto A_q(g,h)\in
\End(V_\gamma)$. The main result of the present subsection is the following
estimate. \begin{prop}\label{poo1} $$| A_q(g,g) |\in C_{\Gamma}(G)\ .$$
\end{prop}
\proof 
Note, that we only have to show finiteness of
the norms
$|.|_{V,N}$, where $V\subset \Omega$ is compact and has the additional
property that $V$ is contained in the interior of a compact subset $V_1\subset\Omega$ satisfying $\gamma V_1\cap V_1=\emptyset$ for all $1\not=\gamma\in\Gamma$. Indeed, any seminorm of $C_\Gamma(G)$
can be majorized by the maximum of finite number of these special ones.

We choose a smooth cut-off function
$\tilde\chi$ on  $\Omega$ such that $\supp(\tilde{\chi})\subset V_1$ and $\supp(1-\tilde\chi)\cap V =\emptyset$.
We further choose a compact $V_2$ containing $V$ in its interior and
being contained in the interior of $V_1$, and a cut-off function
$\chi$ on $\Omega$ such that $\supp(\chi)\subset V_2$ and $\supp(1-\chi)\cap V=\emptyset$.

Then we can write for $k\in V$ 
\begin{eqnarray}
\lefteqn{{}^0P^{T}_{\lambda}\circ \left(ext\circ \pi_* -1\right)\circ  ({}^0P^{R}_{-\lambda})^*(kah ,kah)}&&\nonumber\\&=&\gamma(h)^{-1}\chi(k)\circ 
[\cP^{T}_{\lambda,a}\circ \left(ext\circ \pi_* -1\right)\circ  (\cP^{R}_{-\lambda,a})^*](k,k)\circ \chi(k) \gamma(h)\ .\label{pw7}\end{eqnarray}

In order to employ the off-diagonal localization 
of the Poisson transformation we write
\begin{eqnarray}
\lefteqn{\chi\circ \cP^{T}_{\lambda,a}\circ \left(ext\circ \pi_* -1\right)\circ  (\cP^{R}_{-\lambda,a})^*\circ \chi }\hspace{1cm}&&\nonumber\\&=& 
\chi\circ \cP^{T}_{\lambda,a}\circ (1-\tilde{\chi})\circ  \left(ext\circ \pi_* -1\right)\circ  \tilde{\chi}\circ (\cP^{R}_{-\lambda,a})^*\circ \chi \label{pw1}\\
&& + 
\chi\circ \cP^{T}_{\lambda,a}\circ \left(ext\circ \pi_* -1\right)\circ  (1-\tilde{\chi})\circ (\cP^{R}_{-\lambda,a})^*\circ \chi
\ . \nonumber \end{eqnarray}
In (\ref{pw1}) we could insert the factor $(1-\tilde{\chi})$
since $\tilde{\chi}\circ \left(ext\circ \pi_* -1\right)\circ \tilde{\chi}=0$.
Using that $\supp(\chi)\cap \supp(1-\tilde{\chi})=\emptyset$
we have  
$$\chi\circ \cP^{T}_{\lambda,a}\circ (1-\tilde{\chi})=  a^{-\lambda-\rho}\chi\circ \gamma(w) T  \circ \Phi_{0,-\lambda}\circ J^w_{\sigma,\lambda}\circ (1-\tilde{\chi})
+ a^{-\rho-\alpha} \cR_\infty(\lambda,a)\ ,$$
where 
$\cR_\infty(\lambda,a)\circ \Phi_{\lambda,0}$ remains bounded in
$C^\infty(\imath\aaaa^*,\Hom(H^{\sigma,0}_{-\infty},C^{\infty}(K\times_M V_\gamma)))$ as $a\to\infty$. 

We obtain
\begin{eqnarray}
\lefteqn{\chi\circ \cP^{T}_{\lambda,a}\circ \left(ext\circ \pi_* -1\right)\circ  (\cP^{R}_{-\lambda,a})^*\circ \chi =}&&\nonumber\\&&
a^{-2\rho} a^{-2\lambda}
\chi\circ \gamma(w) T \circ \Phi_{0,-\lambda}\circ    J^w_{\sigma,\lambda}\circ  (1-\tilde{\chi})
\circ \left(ext\circ \pi_* -1\right)\circ\tilde{\chi}\circ\Phi_{\lambda,0}\circ   R^*\frac{c_\gamma(-\lambda)^*}{c_\sigma(\lambda)}
 \circ \chi\label{pw4}\\&&+
a^{-2\rho}  \chi\circ \gamma(w) T   \circ \Phi_{0,-\lambda}\circ  J^w_{\sigma,\lambda}\circ (1-\tilde{\chi}) \circ \left(ext\circ \pi_* -1\right)\circ\tilde{\chi}\circ (J^w_{\tilde{\sigma},-\lambda})^* \circ \Phi_{-\lambda,0}\circ  R^*\tilde{\gamma}(w)^*\circ \chi\label{pw2}\\
&&+
a^{-2\rho} a^{2\lambda} \chi \circ \frac{c_\gamma(\lambda)}{c_\sigma(-\lambda)}T\circ \Phi_{0,\lambda}\circ 
\left(ext\circ \pi_* -1\right)\circ 
(1-\tilde{\chi})\circ   (J^w_{\tilde{\sigma},-\lambda})^*\circ  \Phi_{-\lambda,0}\circ R^*\tilde{\gamma}(w)^*\circ \chi\label{pw5}\\
&&+a^{-2\rho} \chi\circ \gamma(w)T\circ \Phi_{0,-\lambda}\circ  J^w_{\sigma,\lambda}\circ 
\left(ext\circ \pi_* -1\right)\circ 
(1-\tilde{\chi})\circ (J^w_{\tilde{\sigma},-\lambda})^*\circ \Phi_{-\lambda,0}\circ R^*\tilde{\gamma}(w)^*\circ \chi\label{pw3}\\
&&+  Q(\lambda,a)\nonumber\ ,
\end{eqnarray}
where 
$a^{2\rho+\alpha}Q(\lambda,a)$ remains bounded in
$C^\infty(\imath\aaaa^*,\Hom(C^{-\infty}(K\times_M V_\gamma),C^{\infty}(K\times_M V_\gamma)))$ as $a\to\infty$. 

We further compute using that 
the intertwining operators commute with $ext\circ \pi_*$
(compare the proof of Lemma \ref{ll2} for a similar argument), the functional equation of the intertwining operators, and $\chi\circ \left(ext\circ \pi_* -1\right)\circ\chi=0$  
\begin{eqnarray*}
\lefteqn{(\ref{pw2}) + (\ref{pw3})} \hspace{1cm}&&\\
&=& a^{-2\rho}  \chi\circ \gamma(w) T   \circ \Phi_{0,-\lambda}\circ  J^w_{\sigma,\lambda}  \circ \left(ext\circ \pi_* -1\right)\circ\tilde{\chi}\circ (J^w_{\tilde{\sigma},-\lambda})^* \circ \Phi_{-\lambda,0}\circ  R^*\tilde{\gamma}(w)^*\circ \chi\\
&&+ a^{-2\rho} \chi\circ \gamma(w)T\circ \Phi_{0,-\lambda}\circ  J^w_{\sigma,\lambda}\circ 
\left(ext\circ \pi_* -1\right)\circ 
(1-\tilde{\chi})\circ (J^w_{\tilde{\sigma},-\lambda})^*\circ \Phi_{-\lambda,0}\circ R^*\tilde{\gamma}(w)^*\circ \chi\\
&=&a^{-2\rho} \chi\circ \gamma(w)T\circ \Phi_{0,-\lambda}\circ \left(ext\circ \pi_* -1\right)\circ J^w_{\sigma,\lambda}\circ 
(J^w_{\tilde{\sigma},-\lambda})^*\circ \Phi_{-\lambda,0}\circ R^*\tilde{\gamma}(w)^*\circ \chi\\
&=&a^{-2\rho} \chi\circ \gamma(w)T\circ \Phi_{0,-\lambda}\circ \left(ext\circ \pi_* -1\right)\circ  \Phi_{-\lambda,0}\circ R^*\tilde{\gamma}(w)^*\circ \chi\\
&=& a^{-2\rho}   \gamma(w)T\circ \Phi_{0,-\lambda}\circ \chi\circ \left(ext\circ \pi_* -1\right)\circ\chi \circ  \Phi_{-\lambda,0}\circ R^*\tilde{\gamma}(w)^* \\
&=&0\ .
\end{eqnarray*}
Note that in (\ref{pw4}) and (\ref{pw5}) one of the intertwining operators
is localized off-diagonally.
We conclude that the following families of operators
\begin{eqnarray*}
&&\chi\circ \gamma(w) T \circ \Phi_{0,-\lambda}\circ    J^w_{\sigma,\lambda}\circ   (1-\tilde{\chi})
\circ \left(ext\circ \pi_* -1\right)\circ\tilde{\chi}\circ\Phi_{\lambda,0}\circ  R^*\frac{c_\gamma(-\lambda)}{c_\sigma(\lambda)}
 \circ \chi\\
&&
 \chi \circ \frac{c_\gamma(\lambda)}{c_\sigma(-\lambda)}T\circ \Phi_{0,\lambda}
\circ \left(ext\circ \pi_* -1\right)\circ 
(1-\tilde{\chi})\circ   (J^w_{\tilde{\sigma},-\lambda})^*\circ  \Phi_{-\lambda,0}\circ R^*\tilde{\gamma}(w)^*\circ \chi
\end{eqnarray*}
belong to $C^\infty(\imath\aaaa^*,\Hom(C^{-\infty}(K\times_M V_\gamma),C^{\infty}(K\times_M V_\gamma)))$.
Restricting the smooth distribution kernel to the diagonal, multiplying by the
smooth compactly supported function $q$, and integrating over $\imath\aaaa^*$  
we obtain the following estimates using the standard theory of the Euclidean Fourier transform.
For any $N\in\nat$ 
\begin{eqnarray*}
&&\sup_{k\in K}\sup_{a\in A_+}|\int_{\imath\aaaa^*} Q(\lambda,a)(k,k) q(\lambda) d\lambda |
a^{2\rho}(1+|\log(a)|)^N <\infty\\&&
 \sup_{k\in K}\sup_{a\in A_+}\\&& |\int_{\imath\aaaa^*} \chi(k)^2 [\gamma(w) T \circ \Phi_{0,-\lambda}\circ    J^w_{\sigma,\lambda}\circ   (1-\tilde{\chi})
\circ \left(ext\circ \pi_* -1\right)\circ\Phi_{\lambda,0}\circ   R^*\frac{c_\gamma(-\lambda)}{c_\sigma(\lambda)}
](k,k)   q(\lambda) a^{-2\lambda} d\lambda |\\&&
 (1+|\log(a)|)^N <\infty\\
&&\sup_{k\in K} \sup_{a\in A_+}\\&&|\int_{\imath\aaaa^*}  \chi(k)^2 [ \frac{c_\gamma(\lambda)}{c_\sigma(-\lambda)}T\circ\Phi_{0,\lambda}\circ 
 \left(ext\circ \pi_* -1\right)\circ 
(1-\tilde{\chi})\circ  (J^w_{\tilde{\sigma},-\lambda})^*\circ  \Phi_{-\lambda,0}\circ R^*\tilde{\gamma}(w)^*] 
(k,k) q(\lambda) a^{2\lambda}d\lambda |\\&&
 (1+|\log(a)|)^N <\infty \ .
\end{eqnarray*}
This implies the proposition. \hB

\noindent
{\underline{Remark :}
We have shown in fact that $q\mapsto |(A_q)_{|diag}|$
is a continuous map from $C_c^\infty(\imath\aaaa^*)$ to $L^1(\Gamma\backslash
G)$. It would be desirable to extend
this map from $C_c^\infty(\imath\aaaa^*)$ to the Schwartz space
$\cS(\imath\aaaa^*)$. It is this technical problem that prevents us to prove
that the  Fouriertransform $\Phi$  of $\Psi$ restricted  to
the unitary principal series representations is a
tempered distribution. If this would true then it is in fact  a
measure and given by our computations below.

If we would like to show that the map $q\mapsto |(A_q)_{|diag}|$ extends
to a map from the Schwartz space to $L^1(\Gamma\backslash G)$ along the lines
above we need estimates on the growth of $ext$ as the parameter $\lambda$
tends to infinity along the imaginary axis. 
If the imaginary axis is in the domain of convergence of $ext$, i.e.
the critical exponent $\delta_\Gamma$ of $\Gamma$ is negative,
then such an estimate is easy to obtain. In the general case 
one has to estimate the meromorphic continuation of $ext$,  
and this is an open problem.

\subsection{A computation}

In this subsection we want to express 
$\int_{\Gamma \backslash G} \tr\: A_q(g,g)\mu_G(dg)$ in terms of $q$.

Recall that the symmetric space $X=G/K$ can be compactified
by adjoining the boundary $\partial X=G/P$. As a convex-cocompact
group $\Gamma$ acts freely and properly on $X\cup \Omega$ with compact quotient.
Therefore, we can choose a smooth function $\chi^\Gamma\in C^\infty_c(X\cup \Omega)$ such that $\sum_{\gamma\in \Gamma} \gamma^*\chi^\Gamma\equiv 1$
on $X\cup \Omega$. The restriction of $\chi^\Gamma$ to $X$
can be lifted to $G$ as a right-$K$-invariant function which we still denote
by $\chi^\Gamma$. We denote by $\chi_\infty^\Gamma$ the
right-$M$-invariant lift to $K$ of the restriction of $\chi^\Gamma$
to $\partial X=K/M$.  We write 
$$\int_{\Gamma \backslash G} \tr\: A_q(  g,  g)\mu_G(dg)= 
\int_G \chi^\Gamma(g) \tr\:A_q(g,g) \mu_G(dg)\ .$$
Let $\chi_U$ be the characteristic function of the ball $B_U$ in $X$ 
of radius $U$ centered at the origin $[K]$. Again, we denote its  right-$K$-invariant lift to $G$ by the same symbol.
Then we can write
$$\int_G \chi^\Gamma(g) \tr\:A_q(g,g) \mu_G(dg) = 
\lim_{U\to\infty} \int_G \chi^\Gamma(g)\chi_U(g) \tr\:A_q(g,g) \mu_G(dg)\ .$$
Given $U$ we fix a function $\chi_1\in C^\infty_c(G/K)$
such that 
\begin{equation}\label{sup1}\chi_1 \chi_U\chi^\Gamma =\chi_U\chi^\Gamma\ .\end{equation}
The operator
$\chi_U\chi^\Gamma {}^0P^{T}_{\lambda}\circ \left(ext\circ \pi_* -1\right)\circ  ({}^0P^{R}_{-\lambda})^*\chi_1$
has a an integral kernel of compact support. Since the kernel is smooth in the
interior of the support it is of trace class. We can write
$$\int_G \chi^\Gamma(g)\chi_U(g) \tr\:A_q(g,g) \mu_G(dg)=
  \int_{\imath\aaaa^*} \Tr\:[ \chi_U\chi^\Gamma {}^0P^{T}_{\lambda}\circ \left(ext\circ \pi_* -1\right)\circ  ({}^0P^{R}_{-\lambda})^*\chi_1] q(\lambda) d\lambda\ .$$
Note that 
$$\imath\aaaa^*\ni\lambda\mapsto\Tr\:[ \chi_U\chi^\Gamma {}^0P^{T}_{\lambda}\circ \left(ext\circ \pi_* -1\right)\circ  ({}^0P^{R}_{-\lambda})^*\chi_1]$$
is a smooth function. We want to compute its limit in the sense of distributions
as $U\to\infty$ using Green's formula.

Note that $V_B(1_{\rho+\alpha})$ is a complex bundle with a real structure
which is trivial together with this structure. Indeed,
$B$ is orientable, and $V_B(1_{\rho+\alpha})$ is a real power of
$\Lambda^{\max}T^*B$. We choose any non-vanishing positive section 
$\phi\in C^\infty(B,V_B(1_{\rho+\alpha}))$. For any $z\in \C$
we can form $\phi^z\in C^\infty(B,V_B(1_{\rho+z\alpha}))$.
In particular, if we choose $z$ such that $z\alpha=\lambda-\mu$, then
multiplication by $\phi^z$ gives an isomorphism $\bar{\Phi}_{\lambda,\mu}
:C^\infty(B,V_B(\sigma_\mu))\rightarrow C^\infty(B,V_B(\sigma_\lambda))$,
and similar isomorphisms of the spaces of $L^2$- and distribution sections.
If $\Ree(z)=0$, then $ext: C^\infty(B,V_B(1_{\rho+z\alpha}))
\rightarrow {}^\Gamma H^{1,\rho+z\alpha}_{-\infty}$ is regular (indeed
$\rho+z\alpha$ belongs to the domain of convergence).
Multiplication by $ext(\phi^z)$ gives a continuous map
$$\bar{\Phi}_{\lambda,\mu}:H^{\sigma,\mu}_{\infty}\rightarrow H^{\sigma,\lambda}_{-\infty}\ .$$
This map is $\Gamma$-equivariant and extends in fact to
larger subspaces of $H^{\sigma,\mu}_{-\infty}$ of distributions 
which are smooth on neighbourhoods of the limit set $\Lambda$.

The ususal trick to bring in Green's formula is to write
\begin{eqnarray*}\lefteqn{ \Tr\:[ \chi_U\chi^\Gamma {}^0P^{T}_{\lambda}\circ \left(ext\circ \pi_* -1\right)\circ  ({}^0P^{R}_{-\lambda})^*\chi_1]}\hspace{1cm}&&\\ 
&=& \lim_{\mu\to\lambda} \Tr\:[ \chi_U\chi^\Gamma {}^0P^{T}_{\lambda}\circ \left(ext\circ\bar{\Phi}_{\lambda,\mu}\circ  \pi_* - \bar{\Phi}_{\lambda,\mu}\right)\circ  ({}^0P^{R}_{-\mu})^*\chi_1]\ .
\end{eqnarray*}
Let $\nabla^\gamma$ denote the invariant connection of the bundle
$V(\gamma)=G\times_KV_\gamma$ over $X$
and $\Delta_\gamma:=-(\nabla^\gamma)^*\nabla^\gamma$ be the Laplace operator. 
Then there exists a constant
$c \in\R$ such that $(\Delta_\gamma+c +\lambda^2)\circ P^T_\lambda = 0$. Let $n$ denote the outer unit-normal vector field at $B_U$.

By Green's formula we have for $\lambda\not=\pm\mu$
\begin{eqnarray}
\lefteqn{ \Tr\:[ \chi_U\chi^\Gamma {}^0P^{T}_{\lambda}\circ \left(ext\circ\bar{\Phi}_{\lambda,\mu}\circ  \pi_* - \bar{\Phi}_{\lambda,\mu}\right)\circ  ({}^0P^{R}_{-\mu})^*\chi_1]}\hspace{1cm}&&\nonumber\\
 &=&\frac{1}{\lambda^2-\mu^2} \Tr\:  [\Delta_\gamma,\chi^\Gamma ]\chi_U\circ {}^0P^{T}_{\lambda}\circ \left(ext\circ\bar{\Phi}_{\lambda,\mu}\circ \pi_* -  \bar{\Phi}_{\lambda,\mu}\right)\circ ({}^0P^{R}_{-\mu})^*\chi_1 \label{ty1}\\
&&-\frac{1}{\lambda^2-\mu^2} \Tr\: \chi_{|\partial B_U}^\Gamma  \circ  ({}^0P^{T}_{\lambda})_{|\partial B_U}\circ \left(ext\circ\bar{\Phi}_{\lambda,\mu}\circ \pi_* -  \bar{\Phi}_{\lambda,\mu}\right) \circ (\nabla_n^\gamma  {}^0P^{R}_{-\mu})_{|\partial B_U}^*  \label{ty2}\\
&&+ \frac{1}{\lambda^2-\mu^2} \Tr   \: \chi_{|\partial B_U}^\Gamma \circ  (\nabla_n^\gamma {}^0P^{T}_{\lambda})_{|\partial B_U} \circ \left(ext\circ\bar{\Phi}_{\lambda,\mu}\circ \pi_* -  \bar{\Phi}_{\lambda,\mu}\right) \circ ({}^0P^{R}_{-\mu})^*_{|\partial B_U}  \label{ty3}\ .
\end{eqnarray}
Note that the derivatives of $\chi_1$ drop out because of (\ref{sup1}).
Moreover, 
$({}^0P^{T}_{\lambda})_{|\partial B_U}:H^{\sigma,\lambda}_{-\infty}\rightarrow C^\infty(\partial B_U, V(\gamma)_{|\partial B_U})$ denotes the composition of the Poisson transform and restriction to the boundary of $B_U$,
and this operator can be expressed in terms of $\cP^T_{\lambda,a}$.

 We introduce the following  notation.
Let $a_U\in A_+$ be such that $\dist_X(a_UK,K)=U$.
We define $\omega(U):=a_U^{-2\rho}\vol(\partial B_U)$.
Note that $\omega_X:=\lim_{U\to\infty} \omega(U)$ exists.
Let $\chi^\Gamma_a\in C^\infty(K)$ denote the function $k\mapsto \chi^\Gamma(ka)$. Note that $\lim_{a\to\infty} \chi_a=\chi_\infty$. 
Then we can write
\begin{eqnarray*}
 \lefteqn{(\ref{ty2})+(\ref{ty3})}\\ &=& 
-\frac{\omega_Ua_U^{2\rho}}{\lambda^2-\mu^2} \Tr \:     \chi_{a_U}^\Gamma   \circ  \cP^{T}_{\lambda,a_U} \circ  \left(ext\circ\bar{\Phi}_{\lambda,\mu}\circ  \pi_* - \bar{\Phi}_{\lambda,\mu}\right) \circ  (\partial  \cP^{R}_{-\mu,a_U})^* \\
&&+ \frac{\omega_Ua_U^{2\rho}}{\lambda^2-\mu^2} \Tr\:  \chi_{a_U}^\Gamma   \circ \partial \cP^{T}_{\lambda,a_U} \circ  \left(ext\circ\bar{\Phi}_{\lambda,\mu}\circ  \pi_* - \bar{\Phi}_{\lambda,\mu}\right) \circ  (  \cP^{R}_{-\mu,a_U})^*    \ ,
\end{eqnarray*}
Here $\partial \cP^{T}_{\lambda,a}$ stands for the derivative of the function
$a\mapsto  \cP^{T}_{\lambda,a}$ with respect to the positive fundamental
unit vector field on $A$.

Let $\chi$ be a smooth cut-off function on $X\cup\Omega$
of compact support such that 
\begin{equation}\label{chiass}
\gamma\supp(\chi)\cap\supp(\chi)=\emptyset\ ,\quad \forall \gamma\in\Gamma\ , \gamma\not=1\ . \end{equation} Note that
$\chi^\Gamma$ can be decomposed into a finite sum $\chi^\Gamma=\sum_i \chi^i$
such that each $\chi^i$ satisfies (\ref{chiass}). We fix a cut-off function
$\tilde{\chi}$ on $\Omega$ satisfying (\ref{chiass})
and $\tilde{\chi}\equiv 1$ on a neighbourhood of $\supp(\chi)\cap\Omega$.
We further define $\chi_a(k)=\chi(kaK)$, $\chi_\infty(k)=\chi_{|\partial X}(kM)$ and observe that $|\chi_a-\chi_\infty|=O(a^{-\alpha})$ for   any seminorm $|.|$ of $C^\infty(K)$.
Using that $\tilde{\chi}\left(ext\circ\bar{\Phi}_{\lambda,\mu}\circ  \pi_* - \bar{\Phi}_{\lambda,\mu}\right)\tilde{\chi}=0$
 we can write
\begin{eqnarray}
 \lefteqn{(\ref{ty2})+(\ref{ty3})\mbox{\hspace{1cm}($\chi^\Gamma$ replaced by $\chi$)}}&&\nonumber\\
&=&-\frac{\omega_Ua_U^{2\rho}}{\lambda^2-\mu^2} \Tr\:      \chi_{a_U} \circ  \cP^{T}_{\lambda,a_U} \circ (1-\tilde{\chi})\circ \left(ext\circ\bar{\Phi}_{\lambda,\mu}\circ  \pi_* - \bar{\Phi}_{\lambda,\mu}\right) \circ (\partial  \cP^{R}_{-\mu,a_U})^*  \nonumber\\
&&-\frac{\omega_Ua_U^{2\rho}}{\lambda^2-\mu^2} \Tr \:   \chi_{a_U} \circ    \cP^{T}_{\lambda,a_U} \circ  \tilde{\chi}\circ \left(ext\circ\bar{\Phi}_{\lambda,\mu}\circ  \pi_* - \bar{\Phi}_{\lambda,\mu}\right)\circ (1-\tilde{\chi})\circ (\partial  \cP^{R}_{-\mu,a_U})^*  \nonumber\\
&&+ \frac{\omega_Ua_U^{2\rho}}{\lambda^2-\mu^2} \Tr\:     \chi_{a_U}\circ  \partial \cP^{T}_{\lambda,a_U} \circ (1-\tilde{\chi})\circ \left(ext\circ\bar{\Phi}_{\lambda,\mu}\circ  \pi_* - \bar{\Phi}_{\lambda,\mu}\right) \circ(  \cP^{R}_{-\mu,a_U})^*    \nonumber\\
&&+ \frac{\omega_Ua_U^{2\rho}}{\lambda^2-\mu^2} \Tr\:     \chi_{a_U}  \circ \partial \cP^{T}_{\lambda,a_U}  \circ \tilde{\chi}\circ \left(ext\circ\bar{\Phi}_{\lambda,\mu}\circ  \pi_* - \bar{\Phi}_{\lambda,\mu}\right)\circ (1-\tilde{\chi})\circ (  \cP^{R}_{-\mu,a_U})^*  \ .\label{ti1}
\end{eqnarray}

We now insert the asymptotic decomposition (\ref{aass}) of the operators 
$\cP^T_{\lambda,a}$ as $a\to\infty$ noting that
in each line one of these operators is localized  off-diagonally.
In order to stay in trace class operators we choose a function
$\chi_1\in C^\infty(K)$ such that  $\supp(1-\chi_1)\cap \supp(\chi)=\emptyset$ and $\supp(1-\tilde{\chi})\cap\supp(\chi_1)=\emptyset$. We obtain 
\begin{eqnarray}
\lefteqn{(\ref{ty2})+(\ref{ty3}) \mbox{\hspace{1cm}($\chi^\Gamma$ replaced by $\chi$)}}&&\nonumber\\
&=&
-\frac{\omega_Ua_U^{-\mu-\lambda}(-\mu-\rho)}{\lambda^2-\mu^2}    \Tr \:  R^* \frac{c_{\tilde{\gamma}}(-\mu)^*}{c_\sigma(\mu)} \gamma(w)T\chi_\infty   \Phi_{0,-\lambda}J^w_{\sigma,\lambda} (1-\tilde{\chi})  \left(ext\circ\bar{\Phi}_{\lambda,\mu}\circ  \pi_* - \bar{\Phi}_{\lambda,\mu}\right)\Phi_{\mu, 0}\chi_1
\nonumber\\
&&-\frac{\omega_Ua_U^{\mu-\lambda}(\mu-\rho)}{\lambda^2-\mu^2}    \Tr\:  
    R^*  T   \chi_\infty      \Phi_{0,-\lambda}J^w_{\sigma,\lambda} (1-\tilde{\chi}) \left(ext\circ\bar{\Phi}_{\lambda,\mu}\circ  \pi_* - \bar{\Phi}_{\lambda,\mu}\right)(J^w_{\tilde{\sigma},-\mu})^* \Phi_{-\mu,0}\chi_1\nonumber\\
&&-\frac{\omega_Ua_U^{\mu+\lambda}(\mu-\rho)}{\lambda^2-\mu^2} \Tr\:  R^* \gamma(w)^{-1}    \frac{c_\gamma(\lambda)}{c_\sigma(-\lambda)} T \chi_\infty    \Phi_{0,\lambda}   \tilde{\chi}   \left(ext\circ\bar{\Phi}_{\lambda,\mu}\circ  \pi_* - \bar{\Phi}_{\lambda,\mu}\right)(1-\tilde{\chi})    (J^w_{\tilde{\sigma},-\mu})^* \Phi_{-\mu,0}\chi_1 \nonumber\\
&&-\frac{\omega_Ua_U^{\mu-\lambda}(\mu-\rho)}{\lambda^2-\mu^2} \Tr\:R^* T    \chi_\infty       \Phi_{0,\lambda}  J^w_{\sigma,\lambda}  \tilde{\chi}   \left(ext\circ\bar{\Phi}_{\lambda,\mu}\circ  \pi_* - \bar{\Phi}_{\lambda,\mu}\right)(1-\tilde{\chi})     (J^w_{\tilde{\sigma},-\mu})^* \Phi_{-\mu,0}\chi_1\nonumber\\
&&+ \frac{\omega_Ua_U^{-\mu-\lambda}(-\lambda-\rho)}{\lambda^2-\mu^2} \Tr \:   R^* \frac{c_{\tilde{\gamma}}(-\mu)^*}{c_\sigma(\mu)}\gamma(w) T\chi_\infty         \Phi_{0,-\lambda} J^w_{\sigma,\lambda} (1-\tilde{\chi})   \left(ext\circ\bar{\Phi}_{\lambda,\mu}\circ  \pi_* - \bar{\Phi}_{\lambda,\mu}\right)\Phi_{\mu, 0}  \chi_1 \nonumber\\
&&+ \frac{\omega_Ua_U^{\mu-\lambda}(-\lambda-\rho)}{\lambda^2-\mu^2} \Tr\:     R^* T \chi_\infty         \Phi_{0,-\lambda}  J^w_{\sigma,\lambda} (1-\tilde{\chi})   \left(ext\circ\bar{\Phi}_{\lambda,\mu}\circ  \pi_* - \bar{\Phi}_{\lambda,\mu}\right)(J^w_{\tilde{\sigma},-\mu})^* \Phi_{-\mu,0}  \chi_1\nonumber\\
&&+ \frac{\omega_Ua_U^{\mu+\lambda}(\lambda-\rho)}{\lambda^2-\mu^2} \Tr\:    R^* \gamma(w)^{-1}    \frac{c_\gamma(\lambda)}{c_\sigma(-\lambda)} T   \chi_\infty  \Phi_{0,\lambda} \left(ext\circ\bar{\Phi}_{\lambda,\mu}\circ  \pi_* - \bar{\Phi}_{\lambda,\mu}\right)(1-\tilde{\chi})         (J^w_{\tilde{\sigma},-\mu})^* \Phi_{-\mu,0} \chi_1\nonumber\\
&&+ \frac{\omega_Ua_U^{\mu-\lambda}(-\lambda-\rho)}{\lambda^2-\mu^2} \Tr\:     R^* T\chi_\infty     \Phi_{0,-\lambda}J^w_\lambda  \tilde{\chi}   \left(ext\circ\bar{\Phi}_{\lambda,\mu}\circ  \pi_* - \bar{\Phi}_{\lambda,\mu}\right)(1-\tilde{\chi})    (J^w_{\tilde{\sigma},-\mu})^* \Phi_{-\mu,0}\chi_1\nonumber\\
&&+a_U^{-\alpha}\frac{1}{\lambda^2-\mu^2}R_\chi(\lambda,\mu,a_U)\nonumber\\
&=&
-\frac{\omega_Ua_U^{-\mu-\lambda} }{\lambda+\mu} \Tr\:
 R^*  \frac{c_{\tilde{\gamma}}(-\mu)^*}{c_\sigma(-\mu)} \gamma(w)T \chi_\infty \Phi_{0,-\lambda}J^w_{\sigma,\lambda} (1-\tilde{\chi})  \left(ext\circ\bar{\Phi}_{\lambda,\mu}\circ  \pi_* - \bar{\Phi}_{\lambda,\mu}\right) \Phi_{\mu, 0}\chi_1\nonumber\\
&&-\frac{\omega_Ua_U^{\mu-\lambda}}{\lambda-\mu} \Tr\:  R^*    T\chi_\infty  \Phi_{0,-\lambda}J^w_{\sigma,\lambda} (1-\tilde{\chi})\left(ext\circ\bar{\Phi}_{\lambda,\mu}\circ  \pi_* - \bar{\Phi}_{\lambda,\mu}\right)J^w_{\sigma,-\mu} \Phi_{-\mu,0}\chi_1\nonumber\\
&&+\frac{\omega_Ua_U^{\mu+\lambda} }{\lambda+\mu} \Tr\:  R^*\gamma(w)^{-1}    \frac{c_\gamma(\lambda)}{c_\sigma(-\lambda)} T   \chi_\infty   \Phi_{0,\lambda}      \left(ext\circ\bar{\Phi}_{\lambda,\mu}\circ  \pi_* - \bar{\Phi}_{\lambda,\mu}\right)(1-\tilde{\chi})   (J^w_{\tilde{\sigma},-\mu})^*\Phi_{-\mu,0} \chi_1\nonumber\\
&&-\frac{\omega_Ua_U^{\mu-\lambda}}{\lambda-\mu} \Tr\:R^* T    \chi_\infty       \Phi_{0,-\lambda}  J^w_{\sigma,\lambda}  \tilde{\chi}   \left(ext\circ\bar{\Phi}_{\lambda,\mu}\circ  \pi_* - \bar{\Phi}_{\lambda,\mu}\right)(1-\tilde{\chi})    (J^w_{\tilde{\sigma},-\mu})^* \Phi_{-\mu,0}\chi_1\nonumber\\
&&+a_U^{-\alpha}\frac{1}{\lambda^2-\mu^2}R_\chi(\lambda,\mu,a_U)\nonumber\ .
\end{eqnarray}
The remainder $R_\chi(\lambda,\mu,a)$  is holomorphic and  
  remains bounded in $C^\infty(\imath\aaaa^*\times\imath\aaaa^*)$  as $a\to\infty$.

We define  $\langle R,T\rangle\in\C$ such that
$R^*\circ T = \langle R,T\rangle \id_{V_\sigma}$.
If $\sigma$ is not Weyl-invariant, i.e. $\sigma^w\not\cong \sigma$,
then the compositions $R^* c_{\tilde{\gamma}}(-\mu)^* \gamma(w)T$, $R^*\gamma(w)^{-1}    c_\gamma(\lambda) T$ vanish.
If the representation $\sigma$ is Weyl-invariant, then it can be extended
to the normalizer $N_K(M)$ of $M$. In particular, we can define
$\sigma(w)$. In this case we define
$ \langle R,T\rangle(\lambda)\in\C$ such that
$$\sigma(w)R^*\gamma(w)^{-1}\frac{c_\gamma(\lambda)}{c_\sigma(-\lambda)} T  =\langle R,T\rangle(\lambda) \id_{V_\sigma}\ . $$
Note that $R^* \frac{c_{\tilde{\gamma}}(-\mu)^*}{c_\sigma(\mu)} \gamma(w)T\sigma(w)^{-1}=\langle R,T\rangle(-\mu) \id_{V_\sigma}$.
Further we put
$J_{\sigma,\lambda}=\sigma(w)J^w_{\sigma,\lambda}:H^{\sigma,\lambda}_{-\infty}
\rightarrow H^{\sigma,-\lambda}_{-\infty}$.
Then we can write
 \begin{eqnarray}
\lefteqn{ \Tr [    \chi  \chi_U\circ {}^0 P^{T}_{\lambda}\circ \left(ext\circ\bar{\Phi}_{\lambda,\mu}\circ  \pi_* - \bar{\Phi}_{\lambda,\mu}\right)\circ ({}^0P^{R}_{-\mu})^* \chi_1 ]} &&\nonumber\\
&=&
-\frac{\omega_Ua_U^{-\mu-\lambda} \langle R,T\rangle(-\mu)}{\lambda+\mu} \Tr  \: \Phi_{\mu, -\lambda}\circ 
 \chi_\infty     \circ J_{\sigma,\lambda}\circ  (1-\tilde{\chi})  \circ \left(ext\circ\bar{\Phi}_{\lambda,\mu}\circ  \pi_* - \bar{\Phi}_{\lambda,\mu}\right)\circ \chi_1\label{term3}\\
&&-\frac{\omega_Ua_U^{\mu-\lambda}\langle R,T\rangle}{\lambda-\mu} \Tr \:     \chi_\infty       \circ J^w_{\sigma,\lambda} \circ   \left(ext\circ\bar{\Phi}_{\lambda,\mu}\circ  \pi_* - \bar{\Phi}_{\lambda,\mu}\right)\circ(J^w_{\tilde{\sigma},-\mu})^*\circ  \Phi_{-\mu,-\lambda}\circ \chi_1\nonumber\\
&&+\frac{\omega_Ua_U^{\mu+\lambda}\langle R,T\rangle(\lambda) }{\lambda+\mu} \Tr  \:   \chi_\infty      \circ     \left(ext\circ\bar{\Phi}_{\lambda,\mu}\circ  \pi_* - \bar{\Phi}_{\lambda,\mu}\right)\circ(1-\tilde{\chi}) \circ    (J_{\tilde{\sigma},-\mu})^* \circ \Phi_{-\mu,\lambda} \circ \chi_1  \label{term4}\\
&& + \frac{1}{\lambda^2-\mu^2}\Tr^\prime \:     [\Delta_\gamma,\chi ] \circ {}^0P^{T}_{\lambda}\circ  \left(ext\circ\bar{\Phi}_{\lambda,\mu}\circ  \pi_* - \bar{\Phi}_{\lambda,\mu}\right)\circ ({}^0P^{R}_{-\mu})^*  \label{term1}\\&&+\frac{a_U^{-\alpha}}{\lambda^2-\mu^2} Q_\chi(\lambda,\mu,a_U)\nonumber\\
 \lefteqn{Q_\chi(\lambda,\mu,a_U)}\\
&:=&
 R_{\chi}(\lambda,\mu,a_U)\nonumber\\
&&- a_U^\alpha \Tr^\prime\:   [\Delta_\gamma,\chi ](1-\chi_U)\circ {}^0P^{T}_{\lambda}\circ \left(ext\circ\bar{\Phi}_{\lambda,\mu}\circ  \pi_* - \bar{\Phi}_{\lambda,\mu}\right)\circ   ({}^0P^{R^*}_{-\mu})^* \label{term2}\end{eqnarray}
We defer the justification of the terms (\ref{term1}), (\ref{term2})
to Lemma \ref{ll3} below. The functional $\Tr^\prime$ here
is applied to operators with distribution kernels which are continuous on the diagonal, and it takes the integral of its local trace over the diagonal.
Note that the remainder $Q_\chi$ is independent of the choice of $\chi_1$.

The left-hand side of this formula is holomorphic on $\aca\times \aca$.
The terms on the right-hand side may have poles. To aim of the following
discussion is to understand these singularities properly.

  \begin{lem}\label{ll2}
$$\frac{1}{\lambda-\mu}   \Tr \:     \chi_\infty       \circ J^w_{\sigma,\lambda} \circ   \left(ext\circ\bar{\Phi}_{\lambda,\mu}\circ  \pi_* - \bar{\Phi}_{\lambda,\mu}\right)\circ (J^w_{\tilde\sigma,-\mu})^*\circ  \Phi_{-\mu,-\lambda}\circ \chi_1$$ is regular for $\mu=\lambda$.
\end{lem}
\proof
We must show that 
$$\Tr \:      \chi_\infty \circ      J^w_{\sigma,\lambda} \circ  \left(ext \circ  \pi_* - 1\right)\circ(J^w_{\tilde{\sigma},-\lambda})^*\circ\chi_1 =0\ .$$
Recall the definition of the scattering matrix $S^w_{\sigma,\lambda}$
from \cite{bunkeolbrich982}, Def. 5.6. We are going to employ the relations
$$ext\circ S^w_{\sigma,\lambda}= J^w_{\sigma,\lambda}\circ ext\ , \: \pi_*\circ (J^w_{\tilde{\sigma},-\lambda})^* = 
 (S^w_{\tilde{\sigma},-\lambda})^*\circ \pi_*\ , \: S^w_{\sigma,\lambda}\circ (S^w_{\tilde{\sigma},-\lambda})^* = \id\ .$$
We now compute
\begin{eqnarray*}
\lefteqn{  \Tr\: \chi_\infty\circ   J^w_{\sigma,\lambda} \circ \left(ext\circ   \pi_* - 1\right) \circ (J^w_{\tilde{\sigma},-\lambda})^*  \circ\chi_1}\\&=&
 \Tr\:(\chi_\infty\circ   ext  \circ   \hat{S}^w_{\sigma,\lambda}  \circ  (S^w_{\tilde{\sigma},-\lambda})^*\circ \pi_*\circ \chi_1    -         \chi_\infty\circ J^w_{\sigma,\lambda}\circ (J^w_{\tilde{\sigma},-\lambda})^*\circ\chi_1\\
&=&   \Tr\: \chi_\infty \circ ( ext  \circ \pi_*-1)\circ\chi_1  \\
&=&0\ .\end{eqnarray*}
\hB  
In particular we have
\begin{eqnarray}
 &&\lim_{U\to\infty}\lim_{\mu\to\lambda} 
 -\frac{\omega(U)a_U^{\mu-\lambda}\langle R,T\rangle}{\lambda-\mu} \Tr \:     \chi_\infty       \circ J^w_{\sigma,\lambda} \circ   \left(ext\circ\bar{\Phi}_{\lambda,\mu}\circ  \pi_* - \bar{\Phi}_{\lambda,\mu}\right)\circ (J^w_{\tilde{\sigma},-\mu})^*\circ  \Phi_{-\mu,-\lambda}\circ\chi_1\nonumber\\
&=&  \omega_X \langle R,T\rangle  \frac{d}{d\mu}_{|\mu=\lambda} \Tr \:     \chi       \circ J^w_{\sigma,\lambda} \circ   \left(ext\circ\bar{\Phi}_{\lambda,\mu}\circ  \pi_* - \bar{\Phi}_{\lambda,\mu}\right)\circ (J^w_{\tilde{\sigma},-\mu})^*\circ  \Phi_{-\mu,-\lambda}\circ \chi_1\label{pl1}\ .
\end{eqnarray}

If the distribution kernel of an operator 
 admits a continuous restriction to the diagonal, then let $\cD A$ denote this restriction. 
\begin{lem}\label{ll3}
\begin{enumerate}
\item
For any compact subset $Q\subset \imath\aaaa^*$ there
is a constant $C$ such that for all $k\in K$ and $\mu,\lambda\in Q$
$$
 |\cD [\Delta_\gamma,\chi ]  \circ {}^0P^{T}_{\lambda}\circ \left(ext\circ\bar{\Phi}_{\lambda,\mu}\circ  \pi_* - \bar{\Phi}_{\lambda,\mu}\right)\circ ({}^0P^{R}_{-\mu})^*(ka)| 
< C a^{-2\rho-\alpha}\ .$$
 \item We have  $$\Tr^\prime \:  [\Delta_\gamma,\chi^\Gamma ] \circ {}^0 P^{T}_{\lambda}\circ \left(ext\circ\bar{\Phi}_{\lambda,\mu}\circ  \pi_* - \bar{\Phi}_{\lambda,\mu}\right)\circ ({}^0P^{R^*}_{-\mu})^*  =0$$
(note that we consider the cut-off function $\chi^\Gamma$ here).
\end{enumerate}
\end{lem}
\proof
The reason that 1. holds true is that $|d\chi(ka)|\le C a^{-\alpha}$ and 
$|\Delta \chi (ka)|\le C a^{-\alpha}$
uniformly in $k\in K$ and $a\in A$. 
We use the decomposition 
\begin{eqnarray}
\lefteqn{ \cD  [\Delta_\gamma,\chi ]  \circ {}^0P^{T}_{\lambda}\circ \left(ext\circ\bar{\Phi}_{\lambda,\mu}\circ  \pi_* - \bar{\Phi}_{\lambda,\mu}\right)\circ ({}^0P^{R}_{-\mu})^*(ka) }&&\nonumber\\
&=& 
  \cD  [\Delta_\gamma,\chi ] \circ {}^0P^{T}_{\lambda}\circ \tilde{\chi}\circ \left(ext\circ\bar{\Phi}_{\lambda,\mu}\circ  \pi_* - \bar{\Phi}_{\lambda,\mu}\right)\circ (1-\tilde{\chi}) \circ  ({}^0P^{R}_{-\mu})^* (ka)  \label{to1}\\
&& + 
 \cD    [\Delta_\gamma,\chi ] \circ {}^0P^{T}_{\lambda}\circ (1-\tilde{\chi})\circ \left(ext\circ\bar{\Phi}_{\lambda,\mu}\circ  \pi_* - \bar{\Phi}_{\lambda,\mu}\right)\circ ({}^0P^{R}_{-\mu})^* (ka)  \label{to2}
\ . \end{eqnarray}
The asymptotic behaviour (\ref{aass}) of the operators
$\cP^T_{\lambda,a}$ is uniform for $\lambda$ in compact subsets of $\imath \aaaa^*$ and can be differentiated with respect to $a$.
We conclude that for any compact subset $Q\subset \imath\aaaa^*$ there is a constant
$C\in\R$ such for all $\lambda,\mu\in Q$ we have
\begin{eqnarray*}
 \sup_k |\cD   [\Delta_\gamma,\chi ]_a  \circ {}^0P^{T}_{\lambda}\circ (1-\tilde{\chi})\circ \left(ext\circ\bar{\Phi}_{\lambda,\mu}\circ  \pi_* - \bar{\Phi}_{\lambda,\mu}\right)\circ ({}^0P^{R}_{-\mu})^* (ka)|  &\le& C a^{-2\rho-\alpha}\\
 \sup_k|\cD [\Delta_\gamma,\chi ]_a  \circ {}^0P^{T}_{\lambda}\circ (1-\tilde{\chi})\circ \left(ext\circ\bar{\Phi}_{\lambda,\mu}\circ  \pi_* - \bar{\Phi}_{\lambda,\mu}\right)({}^0P^{R}_{-\mu})^* (ka) | &\le& C a^{-2\rho-\alpha}\ .
\end{eqnarray*}

We can write $\chi^\Gamma$ as a finite sum $\chi^\Gamma=\sum_i \chi^i$, where
the cut-off functions $\chi^i$ obeying (\ref{chiass}). 
For each index $i$ we choose
an appropriate cut-off function $\chi_1^i$ as above.
It follows from 1. that
$\tr\:\cD    [\Delta_\gamma,\chi^\Gamma ] \circ {}^0P^{T}_{\lambda}\circ   \left(ext\circ\bar{\Phi}_{\lambda,\mu}\circ  \pi_* - \bar{\Phi}_{\lambda,\mu}\right)\circ ({}^0P^{R}_{-\mu})^*$
is integrable over $G$.
We compute
\begin{eqnarray*}
\lefteqn{\Tr^\prime\:   [\Delta_\gamma,\chi^\Gamma ] \circ {}^0P^{T}_{\lambda}\circ \left(ext\circ\bar{\Phi}_{\lambda,\mu}\circ  \pi_* - \bar{\Phi}_{\lambda,\mu}\right)\circ ({}^0P^{R}_{-\mu})^*}&&\\
&=&\sum_{\gamma\in\Gamma }\Tr^\prime\:    \pi(\gamma)^{-1} \chi^\Gamma \pi(\gamma)  \circ[\Delta_\gamma,\chi^\Gamma ] \circ{}^0 P^{T}_{\lambda}\circ \left(ext\circ\bar{\Phi}_{\lambda,\mu}\circ  \pi_* - \bar{\Phi}_{\lambda,\mu}\right)\circ ({}^0P^{R}_{-\mu})^* \\
&=&\sum_{\gamma\in\Gamma }\Tr^\prime\:    \chi^\Gamma\circ [\Delta_\gamma,(\gamma^{-1})^*\chi^\Gamma ] \circ \pi(\gamma) \circ {}^0P^{T}_{\lambda}\circ \left(ext\circ\bar{\Phi}_{\lambda,\mu}\circ  \pi_* - \bar{\Phi}_{\lambda,\mu}\right)\circ \pi^{\sigma,\mu}(\gamma)^{-1}\circ ({}^0P^{R}_{-\mu})^* \\
&=&\sum_{\gamma\in\Gamma }\Tr^\prime\:   \chi^\Gamma\circ [\Delta_\gamma,(\gamma^{-1})^*\chi^\Gamma ]  \circ {}^0P^{T}_{\lambda}\circ \left(\pi^{\sigma,\lambda}(\gamma) \circ ext\circ\bar{\Phi}_{\lambda,\mu}\circ  \pi_* - \pi^{\sigma,\lambda}(\gamma) \circ \bar{\Phi}_{\lambda,\mu}\right)\circ \pi^{\sigma,\mu}(\gamma)^{-1}\circ ({}^0P^{R^*}_{-\mu})^* \\
&=&\sum_{\gamma\in\Gamma }\Tr^\prime\:  \chi^\Gamma\circ [\Delta_\gamma,(\gamma^{-1})^*\chi^\Gamma ]  \circ {}^0P^{T}_{\lambda}\circ  \left(ext\circ\bar{\Phi}_{\lambda,\mu}\circ  \pi_* - \bar{\Phi}_{\lambda,\mu}\circ \pi^{\sigma,\mu}(\gamma) \right)\circ  \pi^{\sigma,\mu}(\gamma)^{-1}\circ ({}^0P^{R}_{-\mu})^*  \\
&=&\sum_{\gamma\in\Gamma }\Tr^\prime\:    \chi^\Gamma\circ [\Delta_\gamma,(\gamma^{-1})^*\chi^\Gamma ]  \circ {}^0P^{T}_{\lambda}\circ \left(ext\circ\bar{\Phi}_{\lambda,\mu}\circ  \pi_* - \bar{\Phi}_{\lambda,\mu}\right)\circ  \pi^{\sigma,\mu}(\gamma)\circ \pi^{\sigma,\mu}(\gamma)^{-1}\circ ({}^0P^{R}_{-\mu})^* \\
&=&\sum_{\gamma\in\Gamma }\Tr^\prime\:   \chi^\Gamma\circ [\Delta_\gamma,(\gamma^{-1})^*\chi^\Gamma ]  \circ {}^0P^{T}_{\lambda}\circ \left(ext\circ\bar{\Phi}_{\lambda,\mu}\circ  \pi_* - \bar{\Phi}_{\lambda,\mu}\right)\circ ({}^0P^{R}_{-\mu})^*  \\
&=&0
\end{eqnarray*}
since $ \sum_{\gamma\in\Gamma } (\gamma^{-1})^*\chi^\Gamma \equiv 1$.
\hB
Note that the second assertion of the lemma implies
$$\sum_i \Tr  \:  [\Delta_\gamma,\chi^i ] \circ {}^0 P^{T}_{\lambda}\circ \left(ext\circ\bar{\Phi}_{\lambda,\mu}\circ  \pi_* - \bar{\Phi}_{\lambda,\mu}\right)\circ ({}^0P^{R^*}_{-\mu})^*\circ \chi_1^i  =0$$

We now combine (\ref{term3}) and (\ref{term4}) and write
\begin{eqnarray}
\lefteqn{(\ref{term3}) + (\ref{term4})}&&\nonumber\\&=&\frac{\omega(U)(a_U^{\mu+\lambda}-a_U^{-\mu-\lambda})}{\lambda+\mu} \langle R,T\rangle(-\mu) \nonumber\\&&\hspace{1cm} \Tr  \: \Phi_{\mu, -\lambda}\circ 
 \chi_\infty     \circ J_{\sigma,\lambda}\circ  (1-\tilde{\chi})\circ   \left(ext\circ\bar{\Phi}_{\lambda,\mu}\circ  \pi_* - \bar{\Phi}_{\lambda,\mu}\right)\circ \chi_1\label{tu1} \\
&&+\frac{\omega(U)a_U^{\mu+\lambda} }{\lambda+\mu} \Tr  \: \left[ \langle R,T\rangle(\lambda)  \chi_\infty      \circ     \left(ext\circ\bar{\Phi}_{\lambda,\mu}\circ  \pi_* - \bar{\Phi}_{\lambda,\mu}\right)\circ(1-\tilde{\chi}) \circ    (J_{\tilde{\sigma},-\mu})^* \circ \Phi_{-\mu,\lambda} \circ \chi_1 \right. \nonumber\\&&\hspace{2cm}   \left.   -      
  \langle R,T\rangle(-\mu)   \Phi_{\mu, -\lambda}\circ 
 \chi_\infty     \circ J_{\sigma,\lambda}\circ  (1-\tilde{\chi}) \circ  \left(ext\circ\bar{\Phi}_{\lambda,\mu}\circ  \pi_* - \bar{\Phi}_{\lambda,\mu}\right)\circ \chi_1\right]\label{tu2}
\end{eqnarray}
Note that (\ref{tu1}) is regular at $\lambda=\mu=0$.
In fact, if $ext$ has a pole at $\lambda=0$, then it is of first order and 
$J_{\sigma,0}=\id$ (see \cite{bunkeolbrich982}, Prop. 7.4). Hence, the composition  $\chi_\infty \circ J_{\sigma,0}\circ  (1-\tilde{\chi})$ vanishes.

 In order to see
that (\ref{tu2}) is regular at $\lambda=\mu=0$, too, observe in addition that
\begin{equation}\label{re2}
\Tr   \left[ \chi      \circ     \left(ext\circ  \pi_* -  1\right)\circ(1-\tilde{\chi}) \circ    (J_{\tilde{\sigma},0})^* \circ  \chi_1       -\chi     \circ J_{\sigma,0}\circ  (1-\tilde{\chi})  \left(ext \circ  \pi_* - 1\right)\circ \chi_1\right]=0\ .
\end{equation}

Combining Lemma \ref{ll2} and \ref{ll3}, 2., and (\ref{re2}) we 
conclude that $\frac{1}{\lambda^2-\mu^2}Q_{\chi^\Gamma}(\lambda,\mu,a_U)$ is 
regular at $\mu=\lambda$, where we set
$Q_{\chi^\Gamma} :=\sum_i Q_{\chi^i}$.
 Furthermore, locally uniformly in $\lambda$ 
$$\lim_{U\to\infty} \frac{a_U^{-\alpha}}{\lambda^2-\mu^2} Q_{\chi^\Gamma}(\lambda,\mu,a_U)=0\ .
$$ 
By the Lemma of Riemann-Lebesgue we have 
\begin{eqnarray*} &&\lim_{U\to\infty} \frac{\omega(U)a_U^{2\lambda} }{2\lambda} \Tr  \: \left[ \langle R,T\rangle(\lambda)  \chi_\infty      \circ     \left(ext\circ  \pi_* - 1\right)\circ(1-\tilde{\chi}) \circ    (J_{\tilde{\sigma},-\lambda})^* \circ \Phi_{-\lambda,\lambda} \circ \chi_1 \right. \nonumber\\&&\hspace{2cm}   \left.   -      
  \langle R,T\rangle(-\lambda)   \Phi_{\lambda, -\lambda}\circ 
 \chi_\infty     \circ J_{\sigma,\lambda}\circ  (1-\tilde{\chi}) \circ  \left(ext\circ   \pi_* - 1\right)\circ \chi_1\right] \\
&=&0
\end{eqnarray*}
as distributions on $\imath\aaaa^*$.
Moreover,
\begin{eqnarray} \lefteqn{\lim_{U\to\infty} 
\frac{\omega(U)(a_U^{2\lambda}-a_U^{-2\lambda})}{2\lambda} \langle R,T\rangle(-\lambda) }\nonumber\hspace{1cm}&&\\&& \Tr  \: \Phi_{\lambda, -\lambda}\circ 
 \chi_\infty     \circ J_{\sigma,\lambda}\circ  (1-\tilde{\chi}) \circ \left(ext\circ   \pi_* - 1\right)\circ \chi_1  \\
 &=&   \pi \omega_X \langle R,T\rangle(0) \delta_0(\lambda)
\lim_{\lambda\to 0}\left[\Tr  \:   
  \Phi_{\lambda, -\lambda}\circ \chi_\infty     \circ J_{\sigma,\lambda}\circ  (1-\tilde{\chi}) \circ  \left(ext\circ   \pi_* - 1\right)\circ \chi_1\right] \label{pl2}
\end{eqnarray}
in the sense of distributions on $\imath\aaaa^*$.
Combining (\ref{pl1}) and (\ref{pl2})
we now have shown the following proposition.

\begin{prop}\label{poo2}
\begin{eqnarray}
\lefteqn{\int_{\Gamma\backslash G} \tr \: A_q(g,g) \mu_G(dg)}\label{mm2}\\
&=& \omega_X \langle R,T\rangle 
\int_{\imath\aaaa^*}
\sum_i \frac{d}{d\mu}_{|\mu=\lambda} \Tr \: \left[    \chi^i_\infty       \circ J^w_{\sigma,\lambda} \circ   \left(ext\circ\bar{\Phi}_{\lambda,\mu}\circ  \pi_* - \bar{\Phi}_{\lambda,\mu}\right)\circ (J^w_{\tilde{\sigma},-\mu})^*\circ  \Phi_{-\mu,-\lambda}\circ \chi_1^i\right] q(\lambda) d\lambda
 \nonumber\\&+&   \pi \omega_X \langle R,T\rangle(0)  
\sum_i \lim_{\lambda\to 0}\left(\Tr  \: \left[ \Phi_{\lambda, -\lambda}\circ 
 \chi^i_\infty     \circ J_{\sigma,\lambda}\circ   \left(ext\circ  \pi_* - 1\right)\circ\chi_1^i\right]\right) q(0)\nonumber\ .\end{eqnarray}
\end{prop}

Observe that we can rewrite this in the more invariant form
\begin{eqnarray*}
&&\omega_X \langle R,T\rangle 
\int_{\imath\aaaa^*}
 \frac{d}{d\mu}_{|\mu=\lambda} \Tr^\prime \: \left[    \chi^\Gamma_\infty       \circ J^w_{\sigma,\lambda} \circ   \left(ext\circ\bar{\Phi}_{\lambda,\mu}\circ  \pi_* - \bar{\Phi}_{\lambda,\mu}\right)\circ (J^w_{\tilde{\sigma},-\mu})^*\circ  \Phi_{-\mu,-\lambda}   \right] q(\lambda) d\lambda
 \nonumber\\&+& \pi\omega_X \langle R,T\rangle(0)  
 \lim_{\lambda\to 0}\left( \Tr^\prime\:\left[ \Phi_{\lambda, -\lambda}\circ  
 \chi^\Gamma_\infty    \circ J_{\sigma,\lambda}\circ   \left(ext\circ  \pi_* - 1\right) \right]\right)q(0)\ .
\end{eqnarray*}

\section{The Fourier transform of $\Psi$}

\subsection{The contribution of the scattering
matrix}\label{trans}

We consider a symmetric function $q\in C^\infty_c(\imath\aaaa^*, \End(V_\gamma\otimes
\Hom_M(V_\gamma,V_\sigma)))$. For $\pi\in \hat{G}$ we form
$h_q(\pi)\in\End(V_\pi)$ as in Subsection \ref{zc1}. There we have seen
that $\check{R}(h_q)$ and $\check{R}_\Gamma(h_q)$ have smooth integral kernels. 

We choose a basis $\{v_i\}_{i=1,\dots,\dim(\gamma)}$
of $V_\gamma$ and a basis $\{T_j\}_{j=1,\dots, \dim(\Hom_M(V_\gamma,V_\sigma))}$ of $\Hom_M(V_\gamma,V_\sigma)$. 
Let $\{v^i\}$ be a dual basis of $V_{\tilde{\gamma}}$ and $\{T^j\}$
be a dual basis of $\Hom_M(V_\sigma,V_\gamma)=\Hom_M(V_\gamma,V_\sigma)^*$
(with respect to the pairing $\langle T^\prime, T\rangle =\tr_{V_\gamma} T^\prime\circ T$, $T\in  \Hom_M(V_\gamma,V_\sigma)$, $T^\prime\in \Hom_M(V_\sigma,V_\gamma)$).
Then $\{\phi_{ij}:=v_i\otimes T_j\}$ can be considered as a basis of
$H^{\sigma,\lambda}(\gamma)$. Furthermore,
$\{\phi^{ij}:=v^i\otimes (T^j)^*\}$ can be considered as a basis
of $H^{\tilde{\sigma},-\lambda}(\tilde\gamma)=H^{\sigma,\lambda}(\gamma)^*$, and we have
\begin{eqnarray}
\langle \phi_{ij},\phi^{hl}\rangle&=&\int_K \langle T_j \circ \gamma(k)^{-1}
(v_i),(T^l)^* \circ \tilde{\gamma}(k)^{-1} (v^h)\rangle dk\nonumber\\ &=&\int_K
\langle \gamma(k)\circ T^l\circ T_j \circ \gamma(k)^{-1} (v_i),   v^h\rangle
dk\nonumber\\ &=&\frac{\langle T^l,T_j\rangle}{\dim(V_\gamma)} \langle v_i,v^h\rangle\nonumber\\
&=& \frac{1}{\dim(V_\gamma) }\delta_j^l\delta_i^h\ . 
\label{diesehier}\end{eqnarray}
For  $g\in G$ and $w\in V_{\tilde{\gamma}}$ we
define $p^T_{\lambda,w}(g)\in H^{\tilde{\sigma},-\lambda}_\infty$ such that
$\langle P^T_\lambda(\phi)(g),w\rangle = \langle p^T_{\lambda,w}(g),\phi\rangle$
for all $\phi \in H^{\sigma,\lambda}_{-\infty}$. Using (\ref{fr1})
we can write $p^T_{\lambda,w}(g)=\pi^{\tilde{\sigma},-\lambda}(g)(w\otimes T^*)$.
We can write
$$q=\sum_{i,j,k,l} q_{ijkl} v_i\otimes T_j\otimes v^k\otimes T^l\ ,$$
where the functions 
$q_{ijjk}:=\langle v^i\otimes T^j ,q(v_k\otimes T_l)\rangle$ belong to $C^\infty_c(\imath\aaaa^*)$.
Now we can compute
\begin{eqnarray*}
\Tr \:\pi^{\sigma,\lambda}(g)h_q(\pi^{\sigma,\lambda}) &=&
  \dim(V_\gamma)\sum_{i,j} \langle \phi^{ij}, 
\pi^{\sigma,\lambda}(g)h_q(\pi^{\sigma,\lambda}) \phi_{ij}\rangle\\
 &=& \dim(V_\gamma)\sum_{i,j}\langle  P^{T^{j}}_{\lambda }(h_q(\pi^{\sigma,\lambda}) \phi_{ij})(g^{-1}),v^i\rangle\\
&=& \dim(V_\gamma)\sum_{i,j} \langle  h_q(\pi^{\sigma,\lambda}) \phi_{ij}, p^{T^j}_{\lambda,v^i}(g^{-1})\rangle\\
&=& \dim(V_\gamma)\sum_{i,j,k,l} q_{klij}(\lambda)\langle   v_k\otimes T_l, p^{T^j}_{\lambda,v^i}(g^{-1})\rangle\\
&=& \dim(V_\gamma)\sum_{i,j,k,l}q_{klij}(\lambda)\langle P^{T_l^*}_{-\lambda}(p^{T^j}_{\lambda,v^i}(g^{-1}))(1),v_k\rangle\\
&=& \dim(V_\gamma)\sum_{i,j,k,l}q_{klij}(\lambda)    \langle v^i, P^{T^j}_{\lambda}\circ (P^{ T_l^*}_{-\lambda})^*(g^{-1},1)(v_k)\rangle \ .
\end{eqnarray*}
In the last line of this computation 
$P^{T^j}_{\lambda}\circ (P^{ T_l^*}_{-\lambda})^*(g,g^\prime)$
is the integral kernel of the  $G$-equivariant operator $$P^{T^j}_{\lambda}\circ (P^{ T_j^*}_{-\lambda})^*:C^{-\infty}_c(G\times_K V_\gamma)\rightarrow C^{\infty}(G\times_K V_\gamma)\ .$$
We can express the integral kernel of $\check{R}(h_q)$  via Poisson transforms as follows:
\begin{eqnarray}
K_{\check{R}(h_q)}(g,g_1)&=& \int_{\hat{G}}\Tr  \: \pi(g) h_q(\pi) \pi(g_1^{-1})  p(d\pi)\nonumber\\
&=& \frac{\dim(V_\sigma)}{4\pi\omega_X}\int_{\imath\aaaa^*}
\Tr   \: \pi(g_1^{-1}g) h_q(\pi^{\sigma,\lambda})    p_{\sigma} (\lambda) d\lambda\nonumber\\
&=&
\sum_{i,j,k,l}\frac{\dim(V_\sigma)\dim(V_\gamma)}{4\pi\omega_X}\int_{\imath\aaaa^*}q_{klij}(\lambda)    \langle v^i, P^{T^j}_{\lambda}\circ (P^{ T_l^*}_{-\lambda})^*(g^{-1}g_1,1)(v_k)\rangle p_{\sigma} (\lambda) d\lambda\nonumber\\ &=& 
\sum_{i,j,k,l}\frac{\dim(V_\sigma)\dim(V_\gamma)}{4\pi\omega_X}\int_{\imath\aaaa^*}q_{klij}(\lambda)    \langle v^i, P^{T^j}_{\lambda}\circ (P^{ T_l^*}_{-\lambda})^*( g_1,g)(v_k)\rangle p_{\sigma} (\lambda) d\lambda\nonumber\\ &=& 
 \sum_{i,j,k,l}\frac{\dim(V_\sigma)\dim(V_\gamma)}{4\pi\omega_X} \int_{\imath\aaaa^*}q_{klij}(\lambda)    \langle v^i, {}^0 P^{T^j}_{\lambda}\circ ({}^0 P^{ T_l^*}_{-\lambda})^*( g_1,g)(v_k)\rangle   d\lambda \label{bm1} 
\end{eqnarray}
In a similar manner we obtain
\begin{eqnarray*}
\lefteqn{\Tr \:\pi^{\sigma,\lambda}(g)h_q(\pi^{\sigma,\lambda}) \pi^{\sigma,\lambda}(g_1^{-1})ext\:\pi_*}\\&=&\Tr \:h_q(\pi^{\sigma,\lambda}) \pi^{\sigma,\lambda}(g_1^{-1})ext\:\pi_*\:\pi^{\sigma,\lambda}(g)\\
&=& \dim(V_\gamma)\sum_{k,l} \langle \phi^{kl}, 
h_q(\pi^{\sigma,\lambda}) \pi^{\sigma,\lambda}(g_1^{-1})ext\:\pi_*\:\pi^{\sigma,\lambda}(g) \phi_{kl}\rangle\\
&=& \dim(V_\gamma)\sum_{k,l} \langle h_q(\pi^{\sigma,\lambda})^* \phi^{kl},  \pi^{\sigma,\lambda}(g_1^{-1})ext\:\pi_*\:\pi^{\sigma,\lambda}(g) \phi_{kl}\rangle\\
&=& \dim(V_\gamma)\sum_{i,j,k,l} q_{klij}(\lambda)
\langle v^i\otimes T^j, \pi^{\sigma,\lambda}(g_1^{-1})ext\:\pi_*\:\pi^{\sigma,\lambda}(g) \phi_{kl}\rangle\\
&=& \dim(V_\gamma)\sum_{i,j,k,l} q_{klij}(\lambda)
\langle   v^i, P_{\lambda}^{T^j}(ext\:\pi_*\:\pi^{\sigma,\lambda}(g) \phi_{kl})(g_1)\rangle\\
&=& \dim(V_\gamma)\sum_{i,j,k,l} q_{klij}(\lambda)
\langle ext\:\pi_*\:\pi^{\sigma,\lambda}(g) \phi_{kl},p_{\lambda,v^i}^{T^j}(g_1)\rangle\\
&=& \dim(V_\gamma)\sum_{i,j,k,l} q_{klij}(\lambda)
\langle  \phi_{kl},\pi^{\tilde{\sigma},-\lambda}(g^{-1})ext\:\pi_*(p_{\lambda,v^i}^{T^j}(g_1))\rangle\\
&=&\dim(V_\gamma)\sum_{i,j,k,l} q_{klij}(\lambda)
\langle v_k, P^{T_l^*}_{-\lambda}(ext\:\pi_*(p_{\lambda,v^i}^{T^j}(g_1)))(g)\rangle\\
&=&\dim(V_\gamma)\sum_{i,j,k,l} q_{klij}(\lambda)
\langle v^i, P^{T^j}_{\lambda}\circ ext\circ \pi_*\circ (P_{-\lambda}^{T_l^*})^*(g_1,g)(v_k)\rangle
\end{eqnarray*}
and thus
\begin{eqnarray}
K_{\check{R}_\Gamma(h_q)}(g,g_1)&=&
 \frac{\dim(V_\sigma) }{4\pi\omega_X}
\int_{\imath\aaaa^*}
\Tr   \: \pi^{\sigma,\lambda}(g) h_q(\pi^{\sigma,\lambda})\pi^{\sigma,\lambda}(g_1^{-1})    p_{\sigma} (\lambda) d\lambda\nonumber\\
&=&
\sum_{i,j,k,l}\frac{\dim(V_\sigma)\dim(V_\gamma)}{4\pi\omega_X}\int_{\imath\aaaa^*} q_{klij}(\lambda)
\langle v^i, P^{T^j}_{\lambda}\circ ext\circ \pi_*\circ
(P_{-\lambda}^{T_l^*})^*(g_1,g)(v_k)\rangle p_{\sigma} (\lambda)
d\lambda\nonumber\\ &=&  
\sum_{i,j,k,l}\frac{\dim(V_\sigma)\dim(V_\gamma)}{4\pi\omega_X}\int_{\imath\aaaa^*} q_{klij}(\lambda)
\langle v^i, {}^0P^{T^j}_{\lambda}\circ ext\circ \pi_*\circ
({}^0P_{-\lambda}^{T_l^*})^*(g_1,g)(v_k)\rangle   d\lambda \label{bm2}
\end{eqnarray} We conclude that
\begin{eqnarray*}
\lefteqn{K_{\check{R}_\Gamma(h_q)}(g,g_1)-K_{\check{R}(h_q)}(g,g_1)}\\&=&
 \sum_{i,j,k,l}  \frac{\dim(V_\sigma)\dim(V_\gamma)}{4\pi\omega_X}\int_{\imath\aaaa^*}
q_{klij}(\lambda)
\langle v^i, {}^0P^{T^j}_{\lambda}\circ \left(ext\circ \pi_*-1\right)\circ ({}^0P_{-\lambda}^{T_l^*})^*(g_1,g)(v_k)\rangle   d\lambda\\
&=& \frac{\dim(V_\sigma)\dim(V_\gamma)}{4\pi\omega_X}\sum_{i,j,k,l} \langle v^i, A_{q_{klij}}(T^j,T_l^*)(g_1,g)v_k \rangle \ .
\end{eqnarray*}
By Proposition \ref{poo1} the difference $$g\mapsto (K_{\check{R}_\Gamma(h_q)}(g,g)-K_{\check{R}(h_q)}(g,g))$$
is integrable over $\Gamma\backslash G$.
Using the fact that $\gamma$ is irreducible we compute
\begin{eqnarray*}
\lefteqn{\sum_{i,k}\int_{\Gamma\backslash G}
\langle v^i, A_{q_{klij}}(T^j,T_l^*)(g,g)v_k \rangle \mu_{G}(dg)}\hspace{2cm}\\
&=&\sum_{i,k}\int_{\Gamma\backslash G}\int_K
\langle v^i, A_{q_{klij}}(T^j,T_l^*)(gh,gh)v_k\rangle \mu_K(dh)\mu_{G}(dg)\\
&=&\sum_{i,k}\int_{\Gamma\backslash G}\int_K
\langle  v^i, \gamma(h)^{-1}A_{q_{klij}}(T^j,T_l^*)(g,g)\gamma(h)v_k \rangle \mu_K(dh)\mu_{G}(dg)\\
&=&\sum_k\int_{\Gamma\backslash G} 
  \tr\: A_{q_{klkj}}(T^j,T_l^*)(g,g)    \mu_{G}(dg)
\end{eqnarray*}
The following formula is now   an immediate consequence of
Proposition \ref{poo2}.
\begin{eqnarray}
\int_{\Gamma\backslash G}
\lefteqn{\left(K_{\check{R}_\Gamma(h_q)}(g,g)-K_{\check{R} (h_q)}(g,g)\right)\mu_G(dg)}\label{fast}\\
&=&    
 \frac{\dim(V_\sigma)\dim(V_\gamma)}{4\pi } \int_{\imath\aaaa^*}
\sum_i \sum_{k,j,l}\langle T_l^*,T^j\rangle \nonumber\\&& \frac{d}{d\mu}_{|\mu=\lambda} \Tr \: \left[    \chi^i_\infty       \circ J^w_{\sigma,\lambda} \circ   \left(ext\circ\bar{\Phi}_{\lambda,\mu}\circ  \pi_* - \bar{\Phi}_{\lambda,\mu}\right)\circ (J^w_{\tilde{\sigma},-\mu})^*\circ  \Phi_{-\mu,-\lambda}\circ \chi_1^i\right] q_{klkj}(\lambda) d\lambda
 \nonumber\\&+& \frac{\dim(V_\sigma)\dim(V_\gamma)}{4 }   
\sum_i \sum_{k,j,l}\langle T_l^*,T^j\rangle(0) 
\lim_{\lambda\to 0}\left( \Tr  \: \left[  \Phi_{\lambda, -\lambda}\circ
 \chi^i_\infty     \circ J_{\sigma,\lambda}\circ   \left(ext\circ  \pi_* - 1\right)\circ\chi_1^i\right]\right) q_{klkj}(0)\nonumber\ .\nonumber\end{eqnarray}
Now we will rewrite this formula in a more invariant fashion.
Using (\ref{diesehier}) we first compute 
\begin{eqnarray}
&&\sum_{k,j,l} \langle T_{l}^*, T^{j} \rangle q_{klkj}(\lambda)\\ &=&\frac{1}{\dim(V_\sigma)} 
\sum_{j,l,k} T_{l}(T^{j})  
\langle v^k\otimes T^l,q(\lambda) v_k\otimes T_j \rangle\nonumber\\
&=&\frac{1}{\dim(V_\sigma)} 
\sum_{l,k}    
\langle v^k\otimes T^l,q(\lambda) v_k\otimes T_l \rangle\nonumber\\
&=&\frac{1}{\dim(V_\gamma)\dim(V_\sigma)} \Tr\:
h_q(\pi^{\sigma,\lambda})\label{fast1}
\end{eqnarray}
We assume for a moment that $\sigma$ is Weyl-invariant.
For $T\in \Hom_M(V_\sigma,V_\gamma)$ let
$T^\sharp\in \Hom(H^{\sigma,\lambda}_{-\infty},V_\gamma)$ be given by
$T^\sharp(f)=P^T_\lambda(f)(1)$. Recall the relation (\cite{bunkeolbrich982}, Lemma 5.5, 1.)
$$T^\sharp\circ \sigma(w)\circ  J^w_{\sigma,\lambda}= c_\sigma(-\lambda)
[\gamma(w)\circ c_\gamma(\lambda) \circ T \circ \sigma(w)^{-1}]^\sharp\ .$$
We compute
\begin{eqnarray*}
\langle R,T\rangle (0) &=& \frac{1}{\dim(V_\sigma) c_\sigma(0)}\Tr \: \sigma(w)\circ R^*\circ\gamma(w)^{-1}\circ c_\gamma(0) \circ T\\
&=&\frac{1}{\dim(V_\sigma) c_\sigma(0)}\Tr \: \gamma(w)\circ c_\gamma(0) \circ T\circ \sigma(w)^{-1}\circ R^*\\
&=&\frac{1}{\dim(V_\sigma) c_\sigma(0)}\int_K \Tr\:  \gamma(k)^{-1} \circ \gamma(w)\circ c_\gamma(0) \circ T\circ \sigma(w)^{-1}\circ R^*\circ \gamma(k) \mu_K(dk)\\
&=&\frac{1}{\dim(V_\sigma) c_\sigma(0)} \sum_{i} \int_K \langle \gamma(k)^{-1}\circ \gamma(w)\circ c_\gamma(0) \circ T\circ \sigma(w)^{-1}\circ R^*\circ \gamma(k)(v_i),v^i\rangle \mu_K(dk)\\
&=&\frac{1}{\dim(V_\sigma) c_\sigma(0)}\sum_{i} \langle  P^{\gamma(w)\circ c_\gamma(0) \circ T\circ \sigma(w)^{-1}}_0(v_i\otimes R^*)(1),v^i\rangle\\
&=&\frac{1}{\dim(V_\sigma) c_\sigma(0)} \sum_{i} \langle [\gamma(w)\circ c_\gamma(0) \circ T\circ \sigma(w)^{-1}]^\sharp(v_i\otimes R^*),v^i\rangle\\
&=&\frac{1}{\dim(V_\sigma)}\sum_{i} \langle    T^\sharp \circ   \sigma(w)\circ  J^w_{\sigma,0}(v_i\otimes R^*),v^i\rangle\\
&=&\frac{1}{\dim(V_\sigma)}\sum_{i} 
P^T_0\circ    J_{\sigma,0}(v_i\otimes R)(1),v^i\rangle\\
&=&\frac{1}{\dim(V_\sigma)}\sum_{i}\langle     J_{\sigma,0}(v_i\otimes R^*), v^i\otimes T^* \rangle
\end{eqnarray*}
Since $J_{\sigma,0}$ is $K$-equivariant, we can write
$$J_{\sigma,0}(v_i\otimes R)=v_i\otimes j_{\sigma,0}(R)$$
for some $j_{\sigma,0}\in\End(\Hom_M(V_\sigma,V_\gamma))$.
We now have
\begin{eqnarray}
\sum_{k,j,l} \langle T_{l}^*, T^{j} \rangle(0) q_{klkj}(\lambda)
&=&\frac{1}{\dim(V_\sigma)}\sum_{k,j,l}\sum_{i}\langle     v_i\otimes j_{\sigma,0}(T_{l}^*), v^i\otimes T^{j} \rangle
 \langle  v^k\otimes T^l,q(\lambda) v_k\otimes T_j \rangle\nonumber\\
&=&\frac{1}{\dim(V_\sigma)}\sum_{k,l}
 \langle  v^k\otimes T^l,q(\lambda) (v_k\otimes j_{\sigma,0}(T_l^*))  \rangle\\
 \nonumber\\
&=&\frac{1}{\dim(V_\sigma)}\sum_{k,l}
 \langle  v^k\otimes  T^l, q(\lambda) \circ J_{\sigma,0}(v_k\otimes T_l^*) \rangle\nonumber\\
&=&\frac{1}{\dim(V_\sigma)\dim(V_\gamma)} \Tr \: h_q(\pi^{\sigma,0})\circ J_{\sigma,0} \label{fast2}
\end{eqnarray}
Inserting (\ref{fast1}) and (\ref{fast2}) into (\ref{fast})
we obtain the following theorem.

\begin{theorem}\label{main}
If $q$ is a smooth compactly supported symmetric function on $\hat M\times\imath\aaaa^*$
such that $q(\sigma,\lambda)\in \End(V_\gamma\otimes
\Hom_M(V_\gamma,V_\sigma))$, then the difference
$$g\mapsto K_{\check{R}_\Gamma(h_q)}(g,g)-K_{\check{R}(h_q)}(g,g)$$
is integrable over $\Gamma\backslash G$, and we have
\begin{eqnarray*}
\int_{\Gamma\backslash G}
\lefteqn{\left(K_{\check{R}_\Gamma(h_q)}(g,g)-K_{\check{R} (h_q)}(g,g)\right)\mu_G(dg)}\\
&=& \sum_{\sigma\in\hat{M}}\frac{1}{4\pi}  \sum_i 
\int_{\imath\aaaa^*}\\&&
  \frac{d}{d\mu}_{|\mu=\lambda} \Tr \: \left[    \chi^i_\infty       \circ J^w_{\sigma,\lambda} \circ   \left(ext\circ\bar{\Phi}_{\lambda,\mu}\circ  \pi_* - \bar{\Phi}_{\lambda,\mu}\right)\circ (J^w_{\tilde{\sigma},-\mu})^*\circ  \Phi_{-\mu,-\lambda}\circ \chi_1^i\right]  \Tr\:
 h_q(\pi^{\sigma,\lambda})  d\lambda
 \nonumber\\&+&\frac{1}{4}   
\sum_i  \lim_{\lambda\to 0}\left( \Tr  \: \left[\Phi_{\lambda,-\lambda}\circ   
 \chi^i_\infty     \circ J_{\sigma,\lambda}\circ   \left(ext\circ  \pi_* - 1\right)\circ\chi_1^i\right] \right)\Tr \:  h_q(\pi^{\sigma,0})\circ J_{\sigma,0} \ .\end{eqnarray*}
\end{theorem}
We can again rewrite this formula as follows
\begin{eqnarray*}
\int_{\Gamma\backslash G}
\lefteqn{\left(K_{\check{R}_\Gamma(h_q)}(g,g)-K_{\check{R} (h_q)}(g,g)\right)\mu_G(dg)}\\
&=& \sum_{\sigma\in\hat{M}}\frac{1}{4\pi}    
\int_{\imath\aaaa^*}\\&&
  \frac{d}{d\mu}_{|\mu=\lambda} \Tr^\prime \: \left[     \chi^\Gamma_\infty      \circ J^w_{\sigma,\lambda} \circ   \left(ext\circ\bar{\Phi}_{\lambda,\mu}\circ  \pi_* - \bar{\Phi}_{\lambda,\mu}\right)\circ (J^w_{\tilde{\sigma},-\mu})^*\circ  \Phi_{-\mu,-\lambda}  \right]  \Tr\:
 h_q(\pi^{\sigma,\lambda})  d\lambda
 \nonumber\\&+&\frac{1}{4}   
    \lim_{\lambda\to 0}\left( \Tr^\prime  \: \left[  \Phi_{\lambda,-\lambda}\circ
      \chi^\Gamma_\infty\circ J_{\sigma,0}\circ   \left(ext\circ  \pi_* - 1\right) \right]\right) \Tr \:  h_q(\pi^{\sigma,0})\circ J_{\sigma,0} \ .\end{eqnarray*}

\subsection{The Fourier transform of $\Psi$. A conjecture.}\label{conj}

Observe that Theorem \ref{main} does not solve our initial problem
of computing the Fourier transform of the distribution $\Psi$.
The point is that there is no function $f\in C^\infty_c(G)$ such that
its Fourier transform $\hat{f}$ has compact support.
In order to extend Theorem \ref{main} to $\hat{f}$ we must
extend Proposition \ref{poo1} to Schwartz functions.
As explained in the remark at the end of Subsection \ref{esi}
the main obstacle to do this is an estimate
of the growth of the extension map $ext=ext_\lambda:C^{-\infty}(B,V_B(\sigma,\lambda))\rightarrow H^{\sigma,\lambda}_{-\infty}$ as 
$\lambda$ tends to infinity along the imaginary axis. 
  
The goal of the present subsection is to rewrite the result of the computation
of $\Psi^\prime$ in terms of characters thus obtaining the candidate of the
measure $\Phi$. We will also take the discrete spectrum of
$L^2(\Gamma\backslash G)$ into account.

Recall that $$\Tr\:\hat{f}(\pi^{\sigma,\lambda})=\theta_{\pi^{\sigma,\lambda}}(f)\ .$$

If $\pi^{\sigma,0}$ is reducible, then it decomposes
into a sum of $\pi^{\sigma,+}\oplus \pi^{\sigma,-}$ of limits of discrete
series representations which are just the $\pm 1$ eigenspaces of $J_{\sigma,0}$.
In this case $ext$ is regular at $\lambda=0$ (\cite{bunkeolbrich982}, Prop. 7.4.). We can write $$\Tr \:  \hat{f}(\pi^{\sigma,0})\circ J_{\sigma,0} = \theta_{\pi^{\sigma,+}}(f)-\theta_{\pi^{\sigma,-}}(f)\ .$$

If we replace $h_q$ by $\hat{f}$ then formulas
(\ref{bm1}) and (\ref{bm2}) just give the contributions of
the continuous spectrum $K_{R^{ac}(f)}$ and $K_{R^{ac}_\Gamma(f)}$ to the integral kernels $K_{R(f)}$ and $K_{R_\Gamma(f)}$.
We have $K_{R(f)}=K_{R^{ac}(f)}+K_{R^{d}(f)}$, where
$$K_{R^{d}(f)}(g,g_1)=\sum_{\pi\in\hat{G}_d} \Tr\: \pi(g)\hat{f}(\pi)\pi(g_1^{-1})\ .$$
Since we assume that $f$ is $K$-finite, and there are only finitely
many discrete series representations containing a given $K$-type,
this sum is finite.

Furthermore, $K_{R_\Gamma(f)}=K_{R^{ac}_\Gamma(f)}+K_{R^{d}_\Gamma(f)}+
K_{R^{c}_\Gamma(f)}$. Here
$K_{R^{d}_\Gamma(f)}=\sum_{\pi\in \hat{G}_d}K_{R^{d,\pi}_\Gamma(f)}$
is the contribution of discrete series representations
(again a finite sum), 
$$ 
K_{R^{d,\pi}_\Gamma(f)}(g,g_1)=\sum_{i,j} \overline{\langle \psi_i,\phi_j\rangle}
\langle \psi_i,\pi(g)\hat{f}(\pi)\pi(g_1^{-1})\phi_j\rangle\ ,
$$
where $\{\phi_j\}$ and $\{\psi_i\}$ are orthonormal bases
of the infinite-dimensional Hilbert spaces $V_\pi$ and $M_\pi$, respectively.
The finite sum $K_{R^{c}_\Gamma(f)}=\sum_{\pi\in \hat{G}_c}K_{R^{c,\pi}_\Gamma(f)}$ is the discrete contribution 
of representations belonging to $\hat{G}_c$. If we choose orthogonal bases
$\{\phi_j\}$ and $\{\psi_i\}$ of the Hilbert spaces $V_\pi$ and $M_\pi$
(note that $\dim(M_\pi)<\infty$), then we can write
$$K_{R^{c,\pi}_\Gamma(f)}(g,g_1)=\sum_{i,j} \overline{\langle \psi_i,\phi_j\rangle}
\langle \psi_i,\pi(g)\hat{f}(\pi)\pi(g_1^{-1})\phi_j\rangle\ .$$
We define the multiplicity of $\pi$ by 
$N_\Gamma(\pi):=\dim(M_\pi)$.
It is clear that 
$$\int_{\Gamma\backslash G} K_{R^{c}_\Gamma(f)}(g,g) \mu_G(dg)=\sum_{\pi\in\hat{G}_c}
N_\Gamma(\pi) \theta_{\pi}(f)\ .$$

In  Lemma \ref{disfin} we will show that if $f$ is $K$-finite and invariant
under conjugation by $K$, then for each $\pi\in \hat{G}_d$ we have
$$K_{R^{d,\pi}_\Gamma(f)}(g,g)-K_{R^{d,\pi}(f)}(g,g)\in L^1(\Gamma\backslash
G)\ . $$

Given any $A\in V_{\tilde\pi,K}\otimes V_{\pi,K}$ we define
the function $\hat{G}\ni\pi^\prime\mapsto h_A(\pi^\prime)$ to be zero for all
$\pi^\prime\not=\pi$ and $h_A(\pi):=\bar A$, where $\bar{A}:=\int_{K}
\tilde\pi(k)\otimes \pi(k) A dk$. The map 
$$V_{\tilde\pi,K}\otimes V_{\pi,K}\ni A\mapsto T(A):=\int_{\Gamma\backslash G}
[K_{\check R^{d}_\Gamma(h_{\bar A})}(g,g)-K_{\check R^{d}(h_{\bar A})}(g,g)]
\mu_G(dg)$$ is well-defined and a $(\gaaa,K)$-invariant functional on
$V_{\tilde\pi,K}\otimes V_{\pi,K}$. Since $V_{\pi,K}$ is irreducible it
follows that 
$$T(A)= N_\Gamma(\pi) \theta_\pi(A)$$ for some
number $N_\Gamma(\pi)\in\C$ which plays the role of the multiplicity of $\pi$.
Here  we consider $A$ as a finite-dimensional operator
on  $V_\pi$.

It follows from Lemma \ref{disfin} that
$K_{R^{ac}_\Gamma(f)}(g,g)-K_{R^{ac}(f)}(g,g)$ belongs to
$L^1(\Gamma\backslash G)$ not only if $\hat{f}$ is
smooth of compact support, $K$-finite and  $K$-invariant
(Proposition \ref{poo1}), but also
in the case that $f\in C^\infty_c(G)$ is $K$-finite and $K$-conjugation
invariant.  

 The following conjecture provides the candidate for the
measure $\Phi$. Its discrete part is expressed in terms of multiplicities
$N_\Gamma(\pi)$.
If  $\pi\in \hat{G}_c$, then $N_\Gamma(\pi)$
is just the dimension of the space of multiplicities $M_\pi$ and thus a
non-negative integer. If $\pi\in\hat{G}_d$,
then $N_\Gamma(\pi)$ is a sort of regularized dimension of $M_\pi$. We will
show in Proposition \ref{zpro} that $N_\Gamma(\pi)$
is an integer in this case, too. 
The continuous part of the spectrum will contribute a point measure
supported on the irreducible constitutents of the representations $\pi^{\sigma,0}$, and the
corresponding weight will be denoted by $\tilde{N}_\Gamma(\pi)$. The remaining
contribution of the continuous spectrum is absolute continuous
to the Lebesgue measure $d\lambda$ on $\hat G_{ac}=\hat{M}\times
\imath\aaaa^*_+$ and will be described by the density
$L_\Gamma(\pi^{\sigma,\lambda})$.
By Theorem \ref{zetfun} this density appears in the functional equation of
the Selberg zeta function. 
In particular, as it can be already seen from the definition below, $L_\Gamma(\pi^{\sigma,\lambda})$
admits a meromorphic continuation to all of $\aca$ as a function of $\lambda$. Its residues are closely related to the multiplicities of resonances.

\begin{con}\label{zzz}
If $f\in C^\infty_c(G)$ is bi-$K$-finite, then
we have
\begin{eqnarray*}
\Psi(f)&=&\sum_{\sigma\in\hat{M}}
 \frac{1}{4\pi} \int_{\imath\aaaa^*} L_\Gamma(\pi^{\sigma,\lambda})
    \theta_{\pi^{\sigma,\lambda}}(f)
    d\lambda
  \\&+&\sum_{\sigma\in\hat{M}, \sigma\cong \sigma^w, \pi^{\sigma,0}\mbox{\scriptsize red. }} \sum_{\epsilon\in \{+,-\}}
       \tilde{N}_\Gamma(\pi^{\sigma,\epsilon})  \theta_{\pi^{\sigma,\epsilon}}(f) \\
 \\&+&\sum_{\sigma\in\hat{M}, \sigma\cong \sigma^w, \pi^{\sigma,0} \:\mbox{\scriptsize  irred. } } 
      \tilde{N}_\Gamma(\pi^{\sigma,0})  \theta_{\pi^{\sigma,0}}(f) \\
&+&\sum_{\pi\in \hat{G}_{c}\cup \hat{G}_d} N_\Gamma( \pi) \theta_\pi(f) \ ,\end{eqnarray*}
where
$$ L_\Gamma(\pi^{\sigma,\lambda}):=  \frac{d}{d\mu}_{|\mu=\lambda} \Tr^\prime \: \left[      \chi^\Gamma_\infty\circ      J^w_{\sigma,\lambda} \circ   \left(ext\circ\bar{\Phi}_{\lambda,\mu}\circ  \pi_* - \bar{\Phi}_{\lambda,\mu}\right)\circ (J^w_{\tilde{\sigma},-\mu})^*\circ  \Phi_{-\mu,-\lambda} \right]$$
and 
\begin{eqnarray*}\tilde{N}_\Gamma(\pi^{\sigma,\pm})&:=&\pm \frac{ 1}{4}   
\Tr^\prime  \: \lim_{\lambda\to 0}\left(\left[  
   \Phi_{\lambda,-\lambda}\circ \chi^\Gamma_\infty \circ   J_{\sigma,\lambda}\circ   \left(ext\circ  \pi_* - 1\right)  \right] \right)\\
\tilde{N}_\Gamma(\pi^{\sigma,0})&:=&\frac{ 1}{4}   
\Tr^\prime  \: \lim_{\lambda\to 0}\left(\left[  
   \Phi_{\lambda,-\lambda}\circ \chi^\Gamma_\infty \circ   J_{\sigma,\lambda}\circ   \left(ext\circ  \pi_* - 1\right)  \right] \right) \ .
\end{eqnarray*}
\end{con}

Note that the discussion above does not prove this conjecture. What it does
prove is the following theorem.

\begin{theorem}\label{dadd}
If $\hat{G}\ni \pi\mapsto h(\pi)\in \End_K(V_\pi)$ is smooth, of  compact
support, and factorizes over finitely many $K$-types, then
we have
\begin{eqnarray*}
\Psi^\prime(\tilde{\check h})&=&\sum_{\sigma\in\hat{M}}
  \frac{1}{4\pi}\int_{\imath\aaaa^*} L_\Gamma(\pi^{\sigma,\lambda})
    \Tr\: h(\pi^{\sigma,\lambda})
    d\lambda
  \\&+&\sum_{\sigma\in\hat{M}, \sigma\cong \sigma^w, \pi^{\sigma,0}\mbox{\scriptsize red. }} \sum_{\epsilon\in \{+,-\}}
       \tilde{N}_\Gamma(\pi^{\sigma,\epsilon})   \Tr\:
h(\pi^{\sigma,\epsilon}) \\
&+&\sum_{\sigma\in\hat{M}, \sigma\cong \sigma^w, \pi^{\sigma,0} \:\mbox{\scriptsize  irred. } } 
      \tilde{N}_\Gamma(\pi^{\sigma,0})  \Tr\:
h(\pi^{\sigma,0}) \\
 &+&\sum_{\pi\in \hat{G}_{c}\cup \hat{G}_d}
N_\Gamma( \pi) \Tr \: h(\pi) \ .\end{eqnarray*} 
\end{theorem}

        \section{Resolvent kernels and Selberg zeta functions}

\subsection{Meromorphic continuation of resovent kernels}

We fix any $K$-type $\gamma$. The Casimir operator \(
\Omega  \) of \( G \) gives rise to an unbounded selfadjoint operator on the
Hilbert space of sections \( L^{2}(X,V(\gamma )) \). To be precise it is the
unique selfadjoint extension of the restriction of $\Omega$ to the space of
smooth sections with compact support, where we normalize $\Omega$ such that it is bounded from below.
For each complex number $z$ which is
not contained in the spectrum $\sigma (\Omega)$  we let \( R(z) \)
be the operator \( (z-\Omega )^{-1} \). 
By \( r(z) \) we denote its
distribution kernel. 

The Casimir operator descends to an operator acting on sections of
$V_Y(\gamma)$ and induces an unique unbounded selfadjoint operator $\Omega_Y$
on  $L^{2}(Y,V_{Y}(\gamma ))$.
For \( z\not \in \sigma (\Omega _{Y})
\) we define the resolvent \( R_{Y}(z) \) of the operator \( (z-\Omega
_{Y})^{-1} \) and denote by \( r_{Y}(z) \) its
distribution kernel. 

We consider both, \( r(z) \) and \( r_{Y}(z) \), as
distribution sections of the bundle \( V(\gamma )\otimes V(\tilde\gamma )
\) over \( X\times X \). We will in particular  be
interested in the difference \( d(z):=r_{Y}(z)-r(z) \). 

Let \( \hat{M}(\gamma )\subset \hat{M} \) denote the set of irreducible representations
of \( M \) which appear in the restriction of \(
\gamma  \) to \( M \). For each \( \sigma \in \hat{M} \) and \( \lambda \in
\aaaa_{\C}^{*} \) let \( z_{\sigma }(\lambda ) \) denote the
value of \( \Omega  \) on the principal series representation \(
H^{\sigma,\lambda } \). Note that \( z_{\sigma }(\lambda ) \) is of the form
\( z_{\sigma }(\lambda )=c_{\sigma }-\lambda ^{2} \) for some \( c_{\sigma
}\in \R \). We define \( \C_{\gamma } \) to be the branched cover of
\( \C \) to which the inverse functions $$ \lambda _{\sigma
}(z)=\sqrt{c_{\sigma }-z} $$ extend holomorphically for all \( \sigma \in
\hat{M}(\gamma ) \). We fix one sheet \( \C^{phys} \) of \( \C_{\gamma } \)
over the set \( \C\setminus [b,\infty ) \), \( b:=\min _{\sigma \in
\hat{M}}c_{\sigma } \), which we call physical.  We will often consider
$\C^{phys}$ as a subset of $\C$. It follows from the Plancherel theorem for \(
L^{2}(X,V(\gamma )) \) and \( L^{2}(Y,V_{Y}(\gamma )) \) that \( [b,\infty )
\) is the continuous spectrum of \( \Omega  \) and \( \Omega _{Y} \). Thus \(
d(z) \) is defined on the complement of finitely many points of \( \C^{phys}
\) which belong to the discrete spectrum of \( \Omega  \) and \( \Omega _{Y}
\). 

Let \( \Delta  \) denote the diagonal in \( X\times X \) and define \( S:= \)\( \bigcup _{1\not =\gamma \in \Gamma }(1\times \gamma )\Delta \subset X\times X \).
Let \( \Omega _{i} \), \( i=1,2 \), denote the Casimir operators of \( G \)
acting on the first and the second variable of the product \( X\times X \).
The distribution \( d(z) \) satisfies the elliptic differential equation

\begin{equation}
\label{eq3}
(2z-\Omega _{1}-\Omega _{2})d(z)=0
\end{equation}
 on \( X\times X\setminus S \) and is therefore  smooth on this set.

\begin{lem}\label{lem1} 
\( d(z) \)  extends to \( \C_{\gamma } \) as a meromorphic family
of smooth sections of \( V(\gamma )\otimes V(\tilde\gamma ) \)  on \( X\times
X\setminus S \). \end{lem}
\proof
We first show
\begin{lem}
 \( r(z) \) and \( r_{Y}(z) \) extend to \emph{\( \C_{\gamma } \)}
as meromorphic families of distributions.
\end{lem}
\proof 
 We give the argument for \( r_{Y}(z) \)
since \( r(z) \) can be considered as a special case where \( \Gamma  \) is
trivial. Let \( V_{i}\subset X \), \( i=1,2 \) be open subsets such that the
restriction to \( V_{i} \) of the projection \( X\rightarrow Y \) is a diffeomorphism. We consider \( \phi \in C_{c}^{\infty }(V_{1},V(\tilde\gamma
)),\psi \in C_{c}^{\infty }(V_{2},V(\gamma )) \) as compactly supported
sections over $Y$.  

We now employ the Plancherel theorem for \(
L^{2}(Y,V_{Y}(\gamma )) \) in order to show that $$r_{\phi ,\psi
}(z):=r_{Y}(z)(\phi \otimes \psi ) = \langle \phi ,(z-\Omega_Y)^{-1}\psi
\rangle $$ extends meromorphically to \emph{\( \C_{\gamma } \).} We
decompose \( \phi =\sum _{s\in \sigma _{p}(\Omega _{Y})}\phi _{s}+\phi _{ac}
\), \( \psi =\sum _{s\in \sigma _{p}(\Omega _{Y})}\psi _{s}+\psi _{ac} \)
according to the discrete and continuous spectrum of \( \Omega _{Y} \). We
have \( r_{\phi ,\psi }(z)=\sum _{s\in \sigma _{p}(\Omega _{Y})}(z-s)^{-1}
\langle \phi_{s},\psi _{s}\rangle  + \langle  \phi _{ac},(z-\Omega )^{-1}\psi
_{ac}\rangle  \). 
We now employ the Eisenstein Fourier transformation in order to rewrite
the last term of this equation.

For each \( \sigma \in \hat{M}(\gamma ) \), we consider the
normalized Eisenstein series as a meromorphic family of maps
$$C^{-\infty}(B,V_B(\sigma_\lambda))\otimes
\Hom_M(V_\sigma,V_\gamma)\ni f\otimes T\mapsto {}^0 E^T_\lambda(f):={}^0
P^T_\lambda\circ ext(f)\in C^\infty(Y,V_Y(\gamma))\ .$$
For each $\sigma$ let $\{T_i(\sigma)\}_i$, $T_i(\sigma)\in
\Hom(V_{\sigma},V_{\gamma})$, be a base, and let $T^j(\sigma)\in
\Hom(V_{\tilde\sigma},V_{\tilde\gamma})$ be the dual base such that
$T_i(\sigma)^* T^j(\sigma)=\delta_i^j$. If $\phi\in
C^\infty_c(Y,V_Y(\tilde\gamma))$, then its Eisenstein Fourier transform
$EFT_{\tilde\sigma}(\phi)(\lambda)\in C^\infty(B,V_B(\tilde\sigma_{\lambda}))\otimes
\Hom_M(V_{\tilde\sigma},V_{\tilde\gamma})$  is given by  $$\langle
EFT_{\tilde\sigma}(\phi)(\lambda),f\otimes T_i(\sigma)\rangle:=  \langle \phi,
{}^0 E_{-\lambda}^{T_i(\sigma)}(f)\rangle= \langle
({}^0 E^{T_i(\sigma)}_{-\lambda})^*(\phi),f\rangle \ .$$ As a consequence of the
Plancherel theorem we obtain for $ \phi \in C_{c}^{\infty
}(V_{1},V(\tilde\gamma )),\psi \in C_{c}^{\infty }(V_{2},V(\gamma )) $ that
\begin{eqnarray*}
\langle \phi_{ac},\psi_{ac} \rangle &= &\sum _{\sigma \in \hat{M}(\gamma
)}\frac{1}{4\pi  \omega_X} \int _{i\aaaa^{*}}  \langle 
EFT_{\tilde\sigma}(\phi)(-\lambda), EFT_\sigma(\psi)(\lambda)\rangle
d\lambda\\ &=& \sum _{\sigma \in \hat{M}(\gamma )}\frac{1}{4\pi  \omega_X}
\sum_j \int _{i\aaaa^{*}} \langle ({}^0 E^{T_j(\sigma)}_{\lambda})^*(\phi)
,  ({}^0  E^{T^j(\sigma)}_{-\lambda})^*
(\psi)  \rangle d\lambda\\
&=& \sum _{\sigma \in \hat{M}(\gamma )}\frac{1}{4\pi  \omega_X}
\sum_j \int _{i\aaaa^{*}} \langle  \phi
,  {}^0 E^{T_j(\sigma)}_{\lambda} \circ ({}^0  E^{T^j(\sigma)}_{-\lambda})^*
(\psi)  \rangle d\lambda\ .
\end{eqnarray*}
In a similar manner we obtain 
\begin{equation}
\label{eq1}
\langle  \phi_{ac} ,(z-\Omega )^{-1}\psi_{ac}  \rangle =
\sum _{\sigma \in \hat{M}(\gamma )}\frac{1}{4\pi  \omega_X}
\sum_j \int _{i\aaaa^{*}}(z-z_{\sigma }(\lambda ))^{-1} \langle 
 \phi ,  {}^0 E^{T_j(\sigma)}_{\lambda} \circ {}^0  (E^{T^j(\sigma)}_{-\lambda})^*
(\psi )  \rangle d\lambda . \end{equation}
We further investigate the summands in (\ref{eq1}) for each $\sigma$
seperately. So for \( z\in \C^{phys} \) we put 
\begin{equation} \label{eq2} u(z):=\int_{i\aaaa^{*}}(z-z_{\sigma }(\lambda ))^{-1}  \sum_j 
  \langle 
 \phi,  {}^0 E^{T_j(\sigma)}_{\lambda} \circ {}^0  (E^{T^j(\sigma)}_{-\lambda})^*
(\psi )  \rangle
d\lambda  \ . \end{equation}
If we define \( U(\lambda ):=u(z_{\sigma }(\lambda )) \), then it is defined
for \( \Ree (\lambda )>0 \). We claim that 
\[
U(-\lambda ):=U(\lambda )-\frac{2\pi}{ \lambda } \sum_j
 \langle 
 \phi,  {}^0 E^{T_j(\sigma)}_{\lambda} \circ {}^0  (E^{T^j(\sigma)}_{-\lambda})^*
(\psi)  \rangle
\]  provides a meromorphic continuation of \( U \) to all of \(
\aca \). Indeed, \( U(\lambda ) \) is meromorphic for
\( \Ree (\lambda )<0 \), too. We let \( \lambda =\epsilon +i\mu  \) and show
that the jump \( U(\epsilon +i\mu )-U(-\epsilon +i\mu ) \) vanishes as \(
\epsilon \to 0 \). Using the functional equation of the Eisenstein series and
unitarity of the scattering matrix we see that the integrand in (\ref{eq2}) is
a symmetric function. We obtain  
\begin{eqnarray*}\hspace{-0.5cm}\lefteqn{
 \lim _{\epsilon \to 0}[u(z_{\sigma
}(\epsilon +i\mu ))-u(z_{\sigma }(-\epsilon +i\mu ))]}&&\\
&=&\lim _{\epsilon \to 0}
\int_{i\aaaa^{*}} [\frac{1}{z_{\sigma
}(\epsilon +i\mu)-z_{\sigma }(\lambda )} - \frac{1}{z_{\sigma
}(-\epsilon +i\mu)-z_{\sigma }(\lambda )}]
 \sum_j 
  \langle 
 \phi,  {}^0 E^{T_j(\sigma)}_{\lambda} \circ {}^0  (E^{T^j(\sigma)}_{-\lambda})^*
(\psi )  \rangle
d\lambda\\
&=&
 \lim _{\epsilon \to 0}
\int_{\aaaa^{*}} [\frac{1}{-\epsilon^2-2\imath\epsilon\mu+\mu^2-\lambda^2 } - \frac{1}{-\epsilon^2+2\imath\epsilon\mu+\mu^2-\lambda^2 }]
 \sum_j 
  \langle 
 \phi,  {}^0 E^{T_j(\sigma)}_{\imath\lambda} \circ {}^0  (E^{T^j(\sigma)}_{-\imath\lambda})^*
(\psi )  \rangle
d\lambda\\
&=&
 \lim _{\epsilon \to 0}
\int_0^\infty [\frac{1}{-\epsilon^2-2\imath\epsilon\mu+\mu^2-\nu} - \frac{1}{-\epsilon^2+2\imath\epsilon\mu+\mu^2-\nu}]
 \sum_j 
  \langle 
 \phi,  {}^0 E^{T_j(\sigma)}_{\imath\sqrt{\nu}} \circ {}^0  (E^{T^j(\sigma)}_{-\imath\sqrt{\nu}})^*
(\psi )  \rangle
\frac{1}{\sqrt{\nu}}d\nu\\
&=&
 \lim _{\epsilon \to 0}
\int_0^\infty [\frac{1}{-\epsilon-2\imath\mu+\frac{\mu^2-\nu}{\epsilon}} - \frac{1}{-\epsilon+2\imath\mu+\frac{\mu^2-\nu}{\epsilon}}]
 \sum_j 
  \langle 
 \phi,  {}^0 E^{T_j(\sigma)}_{\imath\sqrt{\nu}} \circ {}^0  (E^{T^j(\sigma)}_{-\imath\sqrt{\nu}})^*
(\psi )  \rangle
\frac{1}{\epsilon\sqrt{\nu}}d\nu\\
&=&
 \lim _{\epsilon \to 0}
\int_{-\mu^2/\epsilon}^\infty [\frac{1}{-\epsilon-2\imath\mu-t} - \frac{1}{-\epsilon+2\imath\mu-t}]
 \sum_j 
  \langle 
 \phi,  {}^0 E^{T_j(\sigma)}_{\imath\sqrt{t\epsilon+\mu^2}} \circ {}^0  (E^{T^j(\sigma)}_{-\imath\sqrt{t\epsilon+\mu^2}})^*
(\psi)  \rangle
\frac{1}{\sqrt{t\epsilon+\mu^2}}dt\\
&=&
 \int_{-\infty}^\infty [\frac{1}{-2\imath\mu-t} -\frac{1}{2\imath\mu-t }]
 \sum_j 
  \langle 
 \phi,  {}^0 E^{T_j(\sigma)}_{\imath\sqrt{\mu^2}} \circ {}^0  (E^{T^j(\sigma)}_{-\imath\sqrt{\mu^2}})^*
(\psi )  \rangle
\frac{1}{ \sqrt{ \mu^2}}dt\\
&=&
- \frac{2\pi }{\imath\mu}
 \sum_j 
  \langle 
 \phi,  {}^0 E^{T_j(\sigma)}_{\imath \mu } \circ {}^0  (E^{T^j(\sigma)}_{-\imath \mu })^*
(\psi )  \rangle\ .
  \end{eqnarray*}
It follows from this formula that $U$ is continuous at the imaginary axis.
This finishes the proof
of the claim. We conclude that (\ref{eq1}) extends meromorphically to \(
\C_{\gamma } \). 

Thus we have shown that \( r_{Y}(z) \) is a meromorphic family of distributions
if we put the topology induced by the evaluations against sections of the form
\( \phi \otimes \psi  \). We now argue that it is indeed meromorphic with
respect to the strong topology. It is clear that \( r_{Y}(z) \) is so on the
sheet \( \C^{phys} \). If \( z_{1},z_{2}\in \C_{\gamma } \)
project to \( z\in \C^{phys} \), then by the construction above 
the difference
\( r_{Y}(z_{1})-r_{Y}(z_{2}) \) is a meromorphic family of of smooth sections.
We conclude that \( r_{Y}(z) \) is a meromorphic family of distributions
on \( \C_{\gamma } \). \hB

We  see that \( d(z) \) extends to a meromorphic family of distributions
on \( \C_{\gamma } \). Since it fulfills the differential equation
(\ref{eq3}) on \( \Delta \setminus S \) we conclude that its restriction to
this set is in fact a meromorphic family of smooth sections. This finishes the
proof of Lemma \ref{lem1}. \hB

\subsection{Finite propagation speed estimates}

On \( X\times X \) we define the function \( d_{o}(x,y):=\inf _{1\not =\gamma
\in \Gamma }d(x,\gamma y) \). Note that $S=\{d_0=0\}$. 
Given $\epsilon>0$ we
define the neighbourhood $S_\epsilon:=\{d_0\le \epsilon\}$ of $S$.
Let \( b:=\inf \sigma (\Omega _{Y}) \).

\begin{lem} \label{lem2}
Given $\epsilon > 0$ and a compact subset $W\subset \{\Ree(z)<b\}$
there is a  constant \( C>0 \) 
such that \(
|d(z)(x,y)|<C e^{-(d_{0}(x,y)-\epsilon )\sqrt{b-\Ree(z)}} \) for
all \( (x,y)\not \in S_{2\epsilon} \) and $z\in W$. \end{lem}
\proof
The proof is based on the finite propagation speed of the wave operators \( \cos (tA) \),
 \( \cos (tA_{Y}) \) where \( A:=\sqrt{\Omega -b} \), \( A_{Y}:=\sqrt{\Omega
_{Y}-b} \). We write 
\[
R(z)=\int _{0}^{\infty }e^{-t\sqrt{b-z}}\cos (tA)dt,R_{Y}(z)=\int _{0}^{\infty }e^{-t\sqrt{b-z}}\cos (tA_{Y})dt.\]

 Finite propagation speed gives 
$$ d(z)(x,y)=\int _{d_{0}(x,y)-\epsilon }^{\infty }e^{-t\sqrt{b-z}}[\cos
(tA_{Y})-\cos (tA)]dt\:(x,y) $$ on the level of distribution kernels. By partial
integration \[
(\Omega _{1}-b)^{N}(\Omega _{2}-b)^{N}d(z)(x,y)=(b-z)^{2N}\int _{d_{0}(x,y)-\epsilon }^{\infty }e^{-t\sqrt{b-z}}[\cos (tA_{Y})-\cos (tA)]dt(x,y).\]
 We now employ the fact that \( \int _{d_{0}(x,y)-\epsilon }^{\infty }e^{-t\sqrt{b-z}}\cos (tA_{Y})dt \)
is a bounded operator on \( L^{2}(Y,V_{Y}(\gamma )) \) with norm bounded by
\( C_1 e^{-(d_{0}(x,y)-\epsilon )\Ree (\sqrt{b-z})} \). A similar
estimate holds for the other term. If we choose \( N>\dim(X)/4 \), then we can
conclude that
\[
|d(z)(x,y)|<C e^{-(d_{0}(x,y)-\epsilon )\Ree (\sqrt{b-z})},\]
where \( C \) depends on $W$, \( C_1 \), and a uniform estimate of
norms of delta distributions as functionals on the Sobolev spaces \(
W^{2N,2}(X,V(\gamma )) \) and \( W^{2N,2}(Y,V_{Y}(\gamma )) \) which hold
because \( X,Y \) as well as the bundles \( V(\gamma ),V_{Y}(\gamma ) \) have
bounded geometry. 
We further have employed the fact (which is again a consequence of bounded
geometry) that we can use powers of the  operator \( \Omega ,\Omega _{Y} \) in
order to define the norm of the Sobolev spaces. The assertion of the lemma now
follows from
\( \Ree (\sqrt{b-z})\ge \sqrt{b-\Ree (z)} \).
\hB

The distribution \( r(z) \) is smooth outside the diagonal \( \Delta  \) because
it satisfies a differential equation similar to (\ref{eq3}). For $\epsilon>0$
we define the neighbourhood  $\Delta_\epsilon:=\{d\le \epsilon\}$ of $\Delta$.

\begin{lem} \label{lem3}
For $\epsilon>0$ and a compact subset $W\subset \{\Ree(z)<b\}$
there is a constant and \( C>0 \)
such that \( |r(z)(x,y)|<C e^{-(d(x,y)-\epsilon
)\sqrt{b-\Ree (z)}} \) for all \( (x,y)\not \in \Delta_{2\epsilon}  \) and \(
z\in W\). \end{lem}
\proof
The proof is similar to that of Lemma \ref{lem2}. Using finite propagation
speed we can write
\[
r(z)(x,y)=\int _{d(x,y)-\epsilon }^{\infty }e^{-t\sqrt{b-z}}\cos (tA)dt(x,y),\]

\[
(\Omega _{1}-b)^{N}(\Omega _{2}-b)^{N}r(z)(x,y)=(b-z)^{2N}\int _{d(x,y)-\epsilon }^{\infty }e^{-t\sqrt{b-z}}\cos (tA)dt(x,y).\]
We now argue as in the proof of Lemma \ref{lem2} in order to conclude the estimate.
\hB

Let $L_\gamma$ denote the action of $\gamma\in\Gamma$ on sections of
$V(\gamma)$.
 \begin{lem}\label{lem4}
If \( \sqrt{b-\Ree (z)}>\delta _{\Gamma }+\rho  \), then on \( X\times
X\setminus S \) we have \( d(z)=\sum _{1 \not =\gamma \in \Gamma }(1\otimes
L_{\gamma })r(z) \).  \end{lem}
\proof
It follows from Lemma \ref{lem3} that \( |(1\otimes L_{\gamma
})r(z)(x,y)|<Ce^{-d(x,\gamma ^{-1}y)\sqrt{b-\Ree (z)}} \). Therefore the sum
converges locally uniformly on \( X\times X\setminus S \). The distribution \(
u(z):=r_{Y}(z)-\sum _{\gamma \in \Gamma }(1\otimes L_{\gamma })r(z) \)
satisfies the differential equations  \[
(z-\Omega _{1})u(z)=0,(z-\Omega _{2})u(z)=0,\]
and is therefore a smooth section depending meromorphically on \( z\in
\C^{phys} \). We further have the estimate \( \sum _{1 \not =\gamma \in
\Gamma }|(1\otimes L_{\gamma })r(z)(x,y)|<Ce^{-rd_{0}(x,y)} \), where \(
r<\sqrt{b-\Ree (z)}-\delta _{\Gamma }-\rho  \). For \( \Ree (z)\ll 0 \) we see
that \( u(z) \) defines a bounded operator on \( L^{2}(Y,V_{Y}(\gamma )) \)
and therefore vanishes. Since $u$ is meromorphic in $z$ it vanishes for all
$z$ with \( \sqrt{b-\Ree (z)}>\delta _{\Gamma }+\rho  \). This proves the
lemma. \hB

\subsection{Boundary values}

The meromorphic family of eigenfunctions $d(z)$ on $X\times X\setminus S$ has
meromorphic families of hyperfunction boundary values. Since we consider a
product of rank one spaces it is easy to determine the leading exponents
of a joint eigenfunction of $\Omega_1,\Omega_2$ with eigenvalue $z$.
These exponents are pairs of elements of $\aca$.

\begin{lem}
The set of leading exponents of a joint eigenfunction of $\Omega_1,\Omega_2$
in the bundle $V(\gamma)\otimes V(\tilde\gamma)$ with generic eigenvalue $z$ is 
$$\{\mu_{(\sigma,\epsilon),(\sigma^\prime,\epsilon^\prime)}(z)\:|\:\sigma\in\hat M(\gamma),\sigma^\prime\in\hat M(\tilde\gamma),\:\epsilon,\epsilon^\prime\in\{+,-\}\},$$ where
$\mu_{(\sigma,\epsilon),(\sigma^\prime,\epsilon^\prime)}(z)=(-\rho+\epsilon\lambda_\sigma(z),-\rho+\epsilon^\prime\lambda_{\sigma^\prime}(z))$.
The corresponding boundary value is a section of the bundle
$V(\gamma(\sigma)_{\epsilon\lambda_\sigma(z)})\otimes
V(\tilde\gamma(\sigma^\prime)_{\epsilon^\prime\lambda_{\sigma^\prime}(z)})\rightarrow \partial X\times \partial X$,
where $\gamma(\sigma),\tilde\gamma(\sigma^\prime)$ denote the isotypic
components. \end{lem} 
\proof
An eigensection of $\Omega$ in $V(\gamma)\rightarrow X$  has leading exponents
$-\rho+\epsilon\lambda_\sigma(z)$,  $\epsilon\in\{+,-\}$, and the corresponding
boundary value is a section of $V(\gamma(\sigma)_{\epsilon\lambda_\sigma(z)})$.
This implies the lemma. \hB

Let $\partial\Delta\subset \partial X\times\partial X$ be the diagonal in the
boundary and define $\partial
S:=\bigcup_{1\not=\gamma\Gamma} (\gamma\times 1) (\partial \Delta
\cap\Omega\times\Omega)$. Note that $d(z)$ is a joint eigenfunction in a
neighbourhood in $X\times X$ of $\Omega\times\Omega\setminus \partial S$.
Therefore, for generic $z$ it has boundary values along this set \cite{kashiwaraoshima77}.  
We denote the boundary value associated to the leading exponent
$\tau:=\mu_{(\sigma,\epsilon),(\sigma^\prime,\epsilon^\prime)}(z)$ by
$\beta_\tau(f)$.

 \begin{lem}\label{u8u}
We have  $\beta_\tau(d(z))=0$ (the meromorphic family of hyperfunctions
vanishes) except for $\tau=\mu_{(\sigma,-),(\tilde\sigma^w,-)}(z)$, $\sigma\in
\hat M(\gamma)$, in which case $\beta_\tau(d(z))$ is  a meromorphic family of
real analytic sections.
 \end{lem}
\proof
We employ the fact that $\beta_\tau(d(z))$ depends meromorphically on $z$.
Let $U\subset \Omega$ be such that the restriction of the projection
$\Omega\rightarrow B$ is a diffeomorphism. There is a constant $c>0$ such
that for all $(k_1,k_2)\in U\times U$ we have
$\ee^{d_0(k_1a_1,k_2a_2)}>c|\max(a_1,a_2)|$ (see \cite{bunkeolbrich982}, Cor. 2.4). 
Using Lemma \ref{lem2} we see that for
$\Ree(z)<b$ we have
$|d(z)(k_1a_1,k_2a_2)|< C |\max(a_1,a_2)|^{-\sqrt{b-\Ree(z)}}$, where $C$
depends on $z$. If one of the signs $\epsilon,\epsilon^\prime$ is positive,
for $z\ll 0$ we have $\sqrt{b-z}-2\rho+\epsilon\lambda_\sigma(z)+\epsilon^\prime
\lambda_{\sigma^\prime}(z)>0$. For those $z$ we have
$\lim_{\min(a_1,a_2)\to\infty}  d(z)(k_1a_1,k_2a_2)
(a_1,a_2)^{-\mu_{(\sigma,\epsilon),(\sigma^\prime,\epsilon^\prime)}(z)}=0$
uniformly in $(k_1,k_2)$, where $(a,b)^{(\mu,\nu)}:=a^\mu b^\nu$. 
This  shows that $\beta_\tau(d(z))=0$ 
if one of $\epsilon,\epsilon^\prime$ is positive.

We now consider the kernel $r(z)$ on $X\times X\setminus \Delta$.
 It is a joint eigenfunction of $\Omega_1,\Omega_2$ to the
eigenvalue $z$ on a neighbourhood of $\partial X\times\partial X\setminus
\partial\Delta$ and therefore has hyperfunction boundary values along this set.
A similar argument as above but using Lemma \ref{lem3} instead of
\ref{lem2}  shows that $\beta_\tau(r(z))=0$ except for
$\epsilon=\epsilon^\prime=-$.

Note that $r(z)$ is $G$-invariant in the sense that for $g\in G$ we have
$L_g\otimes L_{g} r(z)=r(z)$. If
$\tau=\mu_{(\sigma,-),(\sigma^\prime,-)}(z)$, then
$\beta_\tau(r(z))$ is a $G$-invariant hyperfunction section of
$V(\gamma(\sigma)_{-\lambda_\sigma(z)})\otimes
V(\tilde\gamma(\sigma^\prime)_{-\lambda_{\sigma^\prime(z)}})$ over 
$\partial X\times\partial X\setminus \partial \Delta$.
Since this set is an orbit of $G$, an invariant hyperfunction
on this set is smooth, and the evaluation at the point $(w,1)\in\partial
X\times \partial X$ provides an injection of the space of invariant
sections into $V_\gamma(\sigma)\otimes V_{\tilde\gamma}(\sigma^\prime)$. If
$b$ is such a $G$-invariant section, then we have for $ma\in MA$
\begin{eqnarray*}
b(w,1)&=&b(maw,ma)\\
&=&b(wm^{w}a^{-1},ma)\\
&=&\gamma(wm^{-1}w)\otimes \tilde\gamma(m^{-1})
a^{\lambda_{\sigma^\prime}(z)-\lambda_\sigma(z)} b(w,1)\ .
\end{eqnarray*}  
Thus $b(w,1)\in [V_{\gamma^w}(\sigma^w)\otimes
V_{\tilde\gamma}(\sigma^\prime)]^M$. We conclude that
$\sigma^\prime\cong\tilde\sigma^w$, and in this case
$\lambda_{\sigma^\prime}(z)=\lambda_\sigma(z)$ holds automatically.
Thus $\beta_\tau(r(z))=0$ if $\sigma^\prime\not\cong\tilde\sigma^w$.

We write 
$$V(\gamma(\sigma)_{-\lambda})\otimes
V(\tilde\gamma(\tilde\sigma^w)_{-\lambda})=V(\sigma_{-\lambda})\otimes
V(\tilde\sigma^w_{-\lambda})\otimes \Hom_M(V_\sigma,V_\gamma)\otimes
\Hom_M(V_{\tilde\sigma^w},V_{\tilde\gamma})\ .$$
The space of invariant sections of $V_{\sigma_{-\lambda}}\otimes
V_{\tilde\sigma^w_{-\lambda}}$ over $\partial X\times\partial X\setminus
\partial\Delta$ is spanned by the distribution kernel $\hat
j^{w}_{\sigma^w_\lambda}$ of the Knapp-Stein intertwining operator
$\hat J^{w}_{\sigma^w_\lambda}$.
We conclude that for each $\sigma\in\hat M(\gamma)$ there is a meromorphic
family $A_\sigma(z)\in \Hom_M(V_\sigma,V_\gamma)\otimes
\Hom_M(V_{\tilde\sigma^w},V_{\tilde\gamma}) $ such that for
$\tau=\mu_{(\sigma,-),(\tilde\sigma^w,-)}(z)$  we have
$\beta_\tau(r(z))=
 \hat j^{w}_{\sigma^w_{\lambda_\sigma(z)}}\otimes A_\sigma(z)$
under the identifications above.

Let $\tau=\mu_{(\sigma,-),(\tilde\sigma^w,-)}(z)$.
We now employ Lemma \ref{lem4} which states that for $\Ree(z)\ll 0$ we have
$d(z)=\sum _{1 \not =\gamma \in \Gamma }(L_\gamma\otimes
1)r(z)$. The sum converges locally uniformly and thus in the space of smooth
section over $X\times X\setminus S$. We further see that
convergence holds locally uniformly in a neighbourhood of
$\Omega\times\Omega\setminus \partial S$. Thus by \cite{oshimasekiguchi80} we can consider distribution boundary values, and by continuity of the boundary value map we
have on $\Omega\times\Omega\setminus \partial S$
 \begin{eqnarray*}
 \beta_{\tau}(d(z))&=&\sum_{1\not=\gamma\in\Gamma} 
(\pi^{\sigma_{-\lambda_\sigma(z)}}(\gamma)\otimes 1) \beta_{\tau}(d(z))\\
&=&\sum_{1\not=\gamma\in\Gamma} 
(\pi^{\sigma_{-\lambda_\sigma(z)}}(\gamma)\otimes 1)
\hat j^{w}_{\sigma^{w},\lambda_\sigma(z)}\otimes A_\sigma(z)\\
&=&(\hat s^{w}_{\sigma^{w},\lambda_\sigma(z)}-\hat
j^{w}_{\sigma^{w},\lambda_\sigma(z)})\otimes A_\sigma(z)\ ,
\end{eqnarray*}
where $\hat s^{w}_{\sigma^{w},\lambda_\sigma(z)}$ is the
distribution kernel of the scattering matrix $\hat
S^{w}_{\sigma^{w},\lambda_\sigma(z)}$.
Here we use the identity $\pi_*\circ \hat
J^{w}_{\sigma^{w},\lambda_\sigma(z)}= \hat
S^{w}_{\sigma^{w},\lambda_\sigma(z)} \circ \pi_*$ which implies that
the distribution kernel of the scattering matrix $\hat
S^w_{\sigma^w,\lambda}$ can be obtained by averaging the distribution kernel of
the Knapp-Stein intertwining operator $\hat
J^w_{\sigma^w,\lambda}$  for $\Ree(\lambda)\gg  0$. 

It follows from the results of  \cite{bunkeolbrich982} that
$\hat s^{w}_{\sigma^{w}_{\lambda_\sigma(z)}}-\hat
j^{w}_{\sigma^{w}_{\lambda_\sigma(z)}}$ extends to a meromorphic family of
smooth sections on all of $\aca$. It follows from \cite{bunkeolbrich991},
Lemma 2.19, 2.20,  that it is indeed a meromorphic family of real analytic sections.
Strictly speaking,  in  \cite{bunkeolbrich991} we only considered the
spherical $M$-type   for $G=SO(1,n)$, but the same arguments can be applied in
the general case. We conclude 
that $\beta_{\tau}(d(z))$ is real analytic as required.

A similar reasoning shows that $\beta_\tau(d(z))=0$ for all $\tau$ which are
not of the form $\mu_{(\sigma,-),(\tilde\sigma^w,-)}(z)$ for some
$\sigma\in\hat M(\gamma)$.
 \hB

\begin{lem}\label{asa}
We have an asymptotic expansion (for generic $z$)
\begin{equation}\label{aq1}
\:d(z)(k_1a,k_2a)\stackrel{a\to\infty}{\sim}\sum_{\sigma\in\hat M(\gamma)}
\sum_{n=0}^\infty a^{-2\rho-2\lambda_\sigma(z)-n\alpha}
p_{z,\sigma,n}(k_1,k_2)\end{equation} which holds locally uniformly for
$k\in\Omega\times\Omega\setminus \partial S$, and where the real analytic
sections $p_{z,\sigma,n}(k_1,k_2)$ of  
$V(\gamma(\sigma)_{-\lambda})\otimes
V(\tilde\gamma(\tilde\sigma^w)_{-\lambda})$   depend meromorphically on $z$.
\end{lem} \proof
Since the boundary value of $d(z)$ along $\Omega\times\Omega\setminus \partial
S$ is real analytic we can employ \cite{oshimasekiguchi80}, Prop. 2.16, in
order to  conclude that $d(z)$ has an asymptotic expansion
with coefficients which depend meromorphically  on $z$.
The formula follows from an inspection of the list of leading exponents Lemma \ref{u8u}.
\hB

 Lemma \ref{asa} has the  following consequence.
For generic $z$ we have
$$
\tr\:d(z)(ka,ka)\stackrel{a\to\infty}{\sim}\sum_{\sigma\in\hat M(\gamma)}
\sum_{n=0}^\infty a^{-2\rho-2\lambda_\sigma(z)-n\alpha}
p_{z,\sigma,n}(k)$$ which holds locally uniformly for $k\in\Omega$,
and where the real analytic functions $p_{z,\sigma,n}$ depend meromorphically
on $z$.

\subsection{The regularized trace of the resolvent}

 \begin{lem}\label{lem6}
The integral $Q_\gamma(z):=\int_X \chi^\Gamma(x) \tr\: d(z) dx$ converges for
$\Ree(z)\ll 0$ and admits a meromorphic continuation to all of $\C_\gamma$. 
\end{lem} 
\proof
Convergence for $\Ree(z)\ll 0$ follows from Lemma \ref{lem2}.
Fix $R\in A$. We write
$Q_\gamma(z)=Q_1(z,R)+Q_2(z,R)$, where
$$Q_1(z,R):=\int_1^R \int_{K} \chi^\Gamma(ka) \tr\: d(z)(ka) dk v(a) da\ ,$$
where $v$ is such that $dk \: v(a) da$ is the volume measure on $X$. Note that
$v(a)\sim a^{2\rho} (\omega_X+a^{-\alpha} c_1+a^{-2\alpha} c_2+...)$ as
$a\to\infty$.  

It is clear that $Q_1(z,R)$ admits a meromorphic continuation.
We have an asymptotic expansion as $a\to\infty$.
$$ u(z,a):=\int_{K} \chi^\Gamma(ka) \tr\: d(z)(ka) dk v(a) \sim
\sum_{\sigma\in\hat M(\gamma)} \sum_{n=0}^\infty
a^{-2\lambda_\sigma(z)-n\alpha} q_{z,\sigma,n}\ ,$$ where $q$ depends
meromorphically on $z$. For $m\in\nat$ let
$$u_m(z,a):=u(z,a)-\sum_{\sigma\in\hat M(\gamma)} \sum_{n=0}^m
a^{-2\lambda_\sigma(z)-n\alpha} q_{z,\sigma,n}\ .$$
Given a compact subset $W$ of $\C_\gamma$ we can choose $m\in \nat_0$ such
that $\int_R^\infty   u_m(z,a) da$ converges (for generic $z$)  and depends
meromorphically on $z$ for all $z\in W$. We further have
$$\sum_{\sigma\in\hat M(\gamma)} \sum_{n=0}^m  \int_R^\infty 
a^{-2\lambda_\sigma(z)-n\alpha} q_{z,\sigma,n}da = \sum_{\sigma\in\hat M(\gamma)} \sum_{n=0}^m
\frac{R^{-2\lambda_\sigma(z)-n\alpha}q_{z,\sigma,n}}{2\lambda_
\sigma(z)+n\alpha}\ ,$$
and this function extends meromorphically to $C_\gamma$.
Since we can choose $W$ arbitrary large we conclude that $Q_\gamma(z)$ admits
a meromorphic continuation to all of $\C_\gamma$. \hB

\subsection{A functional equation}

Let $L(\gamma):=\{c_\sigma \:|\: \sigma\in\hat{M}(\gamma)\}$ be the set of
ramification points of $\C_\gamma$, define $\C^\sharp:=\C\setminus L(\gamma)$, and let
$\C_\gamma^\sharp\subset \C_\gamma$ be the preimage of $\C^\sharp$ under the projection
$\C_\gamma\rightarrow \C$. Then $\C_\gamma^\sharp\rightarrow \C^\sharp$ is a Galois
covering with group of deck transformations $\Pi:=\oplus_{L(\gamma)} \Z_2$.
The action of
$\Pi$ extends to $\C_\gamma$ such $\C_\gamma\setminus \C_\gamma^\sharp$ consists of
fixed points. For $l\in L(\gamma)$ let  $q_l\in \Pi$ be the corresponding
generator. Then we have $\lambda_\sigma(q_lz)=-\lambda_\sigma(z)$ for all
$\sigma\in \hat M(\gamma)$ with $c_\sigma=l$ and 
$\lambda_{\sigma^\prime}(q_l z)=\lambda_{\sigma^\prime}(z)$ else.

\begin{lem}\label{lem7}
For $l\in L(\gamma)$  we have
$$
Q_\gamma(q_lz)-Q_\gamma( z)=\sum_{\sigma\in\hat M(\gamma),
c_\sigma=l}    \frac{-[\gamma:\sigma]}{2
\lambda_\sigma(z)}L_\Gamma(\pi^{\sigma,\lambda_\sigma(z)})
\ .$$\end{lem} 
\proof 
In the
proof of Lemma \ref{lem1} we have seen that (using the notation introduced
there) $$\langle r_Y(q_l z)-r_Y(z), \phi\otimes \psi \rangle  = \sum _{\sigma
\in \hat{M}(\gamma ), c_\sigma=l}\frac{-1}{2 \omega_X  \lambda_\sigma(z)
} \sum_j   \langle   \phi,  {}^0 E^{T_j(\sigma)}_{\lambda_\sigma(z)}
\circ {}^0  (E^{T^j(\sigma)}_{-\lambda_\sigma(z)})^* (\psi )  \rangle  \
. $$  
We conclude that 
$$r_Y(q_l z)-r_Y(z) = \sum_{\sigma\in\hat M(\gamma),
c_\sigma=l}     \frac{-1}{2\omega_X  \lambda_\sigma(z) } \sum_j  
  {}^0 E^{T_j(\sigma)}_{\lambda_\sigma(z)} \circ {}^0 
(E^{T^j(\sigma)}_{-\lambda_\sigma(z)})^*    \ .$$
The same reasoning applies to the trivial group $\Gamma$, where the Eisenstein
series get replaced by the Poisson transformations.
Thus we can write 
$$d(q_l z)-d(z) =
\sum_{\sigma\in\hat M(\gamma), c_\sigma=l}  \frac{-1}{2 \omega_X
 \lambda_\sigma(z) } \sum_j                 {}^0
P^{T_j(\sigma)}_{\lambda_\sigma(z)}             \circ (ext\circ \pi_*-1) \circ
({}^0P^{T^j(\sigma)}_{-\lambda_\sigma(z)})^*     \ .$$ 

The proof of Lemma \ref{lem6} shows that $Q_1(z,R)$ has an asymptotic expansion
$$Q_1(z,R)\sim Q_\gamma(z)+\sum_{\sigma\in\hat M(\gamma)} \sum_{n=0}^\infty
\frac{R^{-2\lambda_\sigma(z)-n\alpha}q_{z,\sigma,n}}{2\lambda_
\sigma(z)+n\alpha}\ .$$
In particular, if  $2\lambda_\sigma(z)\not\in -\nat_0 \alpha$ for all $\sigma$,
then $Q_\gamma(z)$ is the constant term in the asymptotic expansion of $Q_1(z,R)$.

We can now apply Proposition \ref{poo2} which  can be interpreted
as the determination of the constant term (as $R\to\infty$)  of    
$$\int_1^R\int_K\chi^\Gamma(ka) \tr\: \left[{}^0
P^{T_j(\sigma)}_{\lambda }             \circ (ext\circ \pi_*-1) \circ
({}^0P^{T^j(\sigma)}_{-\lambda })^*\right](ka,ka) dk v(a) da$$
as a distribution on $\imath\aaaa^*\setminus\{0\}$.
This shows the desired equation first on $z_\sigma(\imath\aaaa^*)$
and then everywhere by meromorphic continuation. \hB 

\subsection{Selberg zeta functions}

For a detailed investigation of Selberg zeta functions associated to bundles
(for cocompact $\Gamma$) we refer to \cite{bunkeolbrichbuch}.
We assume that $\sigma$ is irreducible and Weyl invariant, or that it is of
the form $\sigma^\prime\oplus(\sigma^\prime)^w$ for some non-Weyl invariant
irreducible $M$-type $\sigma^\prime$. In the latter case we define
$L_\Gamma(\pi^{\sigma,\lambda}):=L_\Gamma(\pi^{\sigma^\prime,\lambda})+L_\Gamma(\pi^{(\sigma^\prime)^w,\lambda})$.

Let $P=MAN$ be a parabolic subgroup of
$M$. If $\gamma\in\Gamma$, then it can be conjugated in $G$ to an element
$m_ga_g\in MA$  with $a_g>1$. Let $\bar{\naaa}$ be the negative root space of
$(\gaaa,\aaaa)$. For $\Ree(\lambda)>\rho$ we can define the Selberg zeta
function $Z_S(\lambda,\sigma)$ by the converging infinite product
$$Z_S(\lambda,\sigma):=\prod_{1\not=[g]\in C\Gamma} \prod_{k=0}^\infty
\det\left( 1-\sigma(m_g)\otimes
S^k(\Ad(m_ga_g)_{|\bar\naaa})a_g^{-\lambda-\rho}\right) \ .$$
In the case of cocompact $\Gamma$ it was shown by \cite{fried86} that
$Z_S(\lambda,\sigma)$ has a meromorphic continuation to all of $\aca$.
In \cite{pattersonperry95} it was explained that the argument of \cite{fried86}
extends to the case of convex cocompact subgroups since it is the compactness
of the non-wandering set of the geodesic flow of $Y$ that matters and not the
compactness of $Y$. Strictly speaking,  \cite{pattersonperry95} deals with 
the spherical case of $SO(1,2n)$, but the argument extends to the general
case. 

There  is a virtual representation $\gamma$ of $K$   (i.e. an element of the
integral representation ring of $K$)  such that $\gamma_{|M}=\sigma$
in the integral representation ring of $M$ (see \cite{miatellovargas83}, \cite{bunkeolbrichbuch}). We call $\gamma$ a lift of
$\sigma$. Note that $\gamma$ is not unique. We can extend the material
developed above to virtual $K$-types by taking the traces with corresponding
signs. Because of the factor $[\gamma:\sigma]$ Lemma \ref{lem7} has the
following corollary.

\begin{kor}\label{kk1}
If  $\gamma$ is a lift of $\sigma$, then 
$Q_\gamma(z)$ extends to a twofold branched cover of $\C$
associated with $\lambda_\sigma(z)$.
\end{kor}

\begin{theorem}\label{zetfun}
The Selberg zeta function satisfies 
$$\frac{Z_S(\lambda,\sigma)}{Z_S(-\lambda,\sigma)}
=\exp\:  \int_0^\lambda L_\Gamma(\pi^{\sigma,\lambda}) d\lambda\ .$$
In particular, the residues of $L_\Gamma(\pi^{\sigma,\lambda})$ are integral.
\end{theorem}
\proof
By \cite{bunkeolbrichbuch} Prop. 3.8. we have
\begin{equation}\label{t5}Q_\gamma(z_\sigma(\lambda))=\frac{1}{2\lambda}
Z^\prime_S(\lambda,\sigma)/Z_S(\lambda,\sigma)\end{equation}
for $\Ree(\lambda)\gg 0$. Indeed, $Q_\gamma(z_\sigma(\lambda))$ is just what
is called in \cite{bunkeolbrichbuch} the hyperbolic contribution associated to
the resolvent. So Corollary   \ref{kk1} and Lemma \ref{lem7} yields the
functional equation of the logarithmic derivative of the Selberg zeta function
$$\frac{Z^\prime_S(-\lambda,\sigma)}{Z_S(-\lambda,\sigma)}+\frac{Z_S^\prime(\lambda,\sigma)}{Z_S(\lambda,\sigma)}=
 L_\Gamma(\pi^{\sigma,\lambda})\ .$$ 
Integrating and employing the apriori information that $Z_S(\lambda,\sigma)$
is meromorphic we obtain the desired functional equation. \hB

\underline{Remarks :}
\begin{enumerate}
\item 
As explained in the introduction it is known 
(from the approach to $Z_S$ using symbolic dynamics and Ruelles thermodynamic formalism)
that $Z_S(\lambda)$ is a meromorphic function of finite order.
It follows that $ L_\Gamma(\pi^{\sigma,\lambda})$, as a function of $\lambda$,
grows at most polynomially.
\item 
In order to describe the singularities of 
$Z_S(\lambda,\sigma)$ we assume that $\gamma$ is an admissible
lift of $\sigma$ (see \cite{bunkeolbrichbuch}). Let $n_{\lambda,\sigma}$ denote the (virtual) dimension of the subspace of the $L^2$-kernel of $\Omega_Y-z_\sigma(\lambda)$ on $V_Y(\gamma)$ which is generated by
non-discrete series representations of $G$.
It follows from Theorem \ref{zetfun} that
$$\ord_{\lambda=\mu} Z_S(\lambda,\sigma)=\left(\begin{array}{cc}
  \res_{\lambda=\mu}L_\Gamma(\pi^{\sigma,\lambda})+n_{-\mu,\sigma}&\quad  \Ree(\mu)<0\\
n_{\mu,\sigma}&\quad 0 <\Ree(\mu)\end{array}\right.$$ 
\item If $\mu$ is non-integral and if $ext$ has a pole at $\mu$ 
of at most first order, then
one has $$\res_{\lambda=\mu}L_\Gamma(\pi^{\sigma,\lambda})=\dim\: {}^\Gamma C^{-\infty}(\Lambda,V(\sigma,\mu))\ ,$$ where ${}^\Gamma C^{-\infty}(\Lambda,V(\sigma,\mu))$ is the space of 
invariant distributions with support on the limit set $\Lambda$.
If $ext$ has  higher order singularity, then the residue has a similar interpretation (see \cite{bunkeolbrich991}, Prop. 5.6)
This provides an independent argument for the integrality of the residues of
$L_\Gamma(\pi^{\sigma,\lambda})$. 
\end{enumerate}

\subsection{Integrality of $N_\Gamma(\pi)$ for discrete series representations}

Let $\pi$ be a discrete series representation of $G$ containing the $K$-type
$\tilde\gamma$. There are embeddings $M_\pi\otimes
V_\pi(\tilde\gamma)\hookrightarrow L^2(\Gamma\backslash G)(\tilde\gamma)$
and $V_{\tilde\pi}\otimes V_\pi(\tilde\gamma)\hookrightarrow 
L^2(G)$. Let $A\in \End_K(V_\pi(\tilde\gamma))$ be given. We extend
$A$ by zero to the orthogonal complement of $V_\pi(\tilde\gamma)$ thus
obtaining an operator in  $\End_K(V_\pi)$ which we will still denote by $A$.
The operator $A$ induces operators $\check R_\Gamma(h_A)$ and $\check
R_\Gamma(h_A)$ on $L^2(\Gamma\backslash G)$ and $L^2(G)$,
where $h_A$ is supported on $\{\pi\}\subset \hat G$ and $h_A(\pi):=A$.

\begin{lem}\label{disfin}
We have
$$ K_{\check R_\Gamma(h_A)}(g,g)- K_{\check R(h_A)}(g,g) \in
L^1(\Gamma\backslash G)\ .$$ \end{lem}
\proof
Let $D(G,\gamma)$ be the algebra of invariant differential operators
on $V(\gamma)$. It is isomorphic to $(\cU(\gaaa)\otimes_{\cU(\kaaa)}
\End(V_{\gamma}))^K$.
If $\pi^\prime$ is an admissible representation of $G$, then
$D(G, \gamma)$ acts in a natural way on
$ (V_{\pi^\prime}\otimes V_{\gamma})^K$.
If $\pi^\prime$ is irreducible, then $ (V_{\pi^\prime}\otimes V_{ \gamma})^K$ is an
irreducible representation of $D(G, \gamma)$.
The correspondence $\pi^\prime\mapsto  (V_{\pi^\prime}\otimes V_{ \gamma})$
provides a bijection between the sets of equivalence classes of irreducible
representations of $G$  containing the $K$-type $\tilde\gamma$
and irreducible representations of $D(G, \gamma)$.

Note that $\End_K(V_\pi( \tilde\gamma))   \cong \End( (V_\pi\otimes V_{ \gamma})^K)$. We conclude that
there is $D_A\in D(G, \gamma)$ that induces the endomorphism $A$ on
$V_\pi(\tilde \gamma)$. 
Let $z_0$ be the eigenvalue of the Casimir operator on $\pi$
and $Z$ be the finite set of irreducible representations
of $G$ containing the $K$-type $\tilde \gamma$ such that $\Omega$ acts with
eigenvalue $z_0$. Then we can choose $D_A$ such that it vanishes on all
$(V_{\pi^\prime}\otimes V_\gamma)^K$ for $\pi^\prime\in Z$,
$\pi^\prime\not=\pi$.

For simplicity we assume that $z_0$ is not a branching point of $R(z)$.
In the latter case the following argument can easily be modified.
The operators $D_AR_Y(z)$ and $D_A R(z)$ have  poles at $z_0$
with residues $K_{\check R_\Gamma(h_A)}$ and $K_{\check R(h_A)}$.
The difference $(D_A)_1 d(z) := (D_A)_1 r_Y(z)-(D_A)_1 r(z)$
of distribution kernels is still a meromorphic family of joint eigenfunctions
with real analytic boundary values along $\Omega\times\Omega\setminus \partial
S$. We have the asymptotic expansion
\begin{equation}\label{aq2}
\: (D_A)_1 d(z)(k_1a,k_2a)\stackrel{a\to\infty}{\sim}\sum_{\sigma\in\hat M(\gamma)}
\sum_{n=0}^\infty a^{-2\rho-2\lambda_\sigma(z)-n\alpha}
p_{z,\sigma,n,A}(k_1,k_2)\ .\end{equation} 
The residue of $(D_A)_1 d(z)$ at $z_0$ can be computed by integrating
$(D_A)_1 d(z)$  along a small circle counter-clockwise surrounding $z_0$.
 If we insert the asymptotic expansion  (\ref{aq2}) into this integral, then we
 obtain an asymptotic expansion 
\begin{eqnarray*}
\res_{z=z_0}(D_A)_1
d(z)(k_1a,k_2a)&\stackrel{a\to\infty}{\sim}&\sum_{\sigma\in\hat M(\gamma)}
\sum_{n=0}^\infty     \res_{z=z_0} a^{-2\rho-2\lambda_\sigma(z)-n\alpha}
p_{z,\sigma,n,A}(k_1,k_2)\\ &\sim&\sum_{\sigma\in\hat M(\gamma)}
\sum_{n=0}^\infty   \sum_{m=0}^{finite} \log(a)^m a^{-2\rho-2\lambda_\sigma(z_0)-n\alpha}
  p_{z_0,\sigma,n,m,A}(k_1,k_2)\ ,
\end{eqnarray*}
where $p_{z_0,\sigma,n,m,A}$ is a real analytic section on
$\Omega\times\Omega\setminus  \partial S$.
 Since $K_{\check R_\Gamma(h_A)}$ and $K_{\check R(h_A)}$
project onto eigenspaces of square integrable sections
we have for $k_1\not= k_2$
\begin{eqnarray*}
\lefteqn{K_{\check R(h_A)}(k_1a_1,k_2a_2)\stackrel{a_i\to\infty}{\sim}}&&\\
&&\sum_{n_1,n_2=0}^\infty\sum_{m_1,m_2}^{finite}
 \log(a_1)^{m_1} \log(a_2)^{m_2} a_1^{-\rho-\lambda_\sigma(z_0)-n_1\alpha}
a_2^{-\rho-\lambda_\sigma(z_0)-n_2\alpha}
  p_{z_0,\sigma,n_1,n_2,m_1,m_2,A}(k_1,k_2)\end{eqnarray*}
with $p_{z_0,\sigma,n_1,n_2,m_1,m_2,A}(k_1,k_2)=0$ as long as
$-\lambda_\sigma(z_0)-n_1\alpha\ge 0$, $-\lambda_\sigma(z_0)-n_2\alpha\ge 0$,
and similar for $K_{\check R_\Gamma(h_A)}$.
We conclude that $\Ree(-2\lambda_\sigma(z)-n\alpha)<0$ if $p_{z_0,\sigma,n,m,A}\not=0$.
The assertion of the lemma now follows. \hB

\begin{prop}\label{zpro}
If $\pi$ be a representation of the discrete series of $G$, then
$N_\Gamma(\pi)\in\Z$.
\end{prop}
\proof
There exists an invariant generalized Dirac operator
$D$ acting on a graded vector bundle $E\rightarrow X$, $E=E^+\oplus E^-$,
such that $V_{\tilde\pi}\oplus \{0\}$ is the kernel of $D$. If $\gamma$ is the
graded $K$-type associated to $E$, then 
\begin{equation}\label{klopt}\gamma_{|M}=0\ .\end{equation} 
Let $D_Y$ be the
induced operator on $Y$. The distribution kernels of $r(z):=(z-D^2)^{-1}$ and
$r_Y(z):=(z-D_Y^2)^{-1}$  have  meromorphic continuations to a branched
covering of $\C$. Their difference goes into the functional $\Psi^\prime$. The
function $$Q(z):=\int_{\Gamma\backslash G} \tr\:(r_Y(z)-r(z))(g,g) dg$$ has a
meromorphic continuation to all of $\C$ by (\ref{klopt}) and Lemma \ref{lem7}. Its residue at $z=0$ is given by
$$\res_{z=0} Q(z)=N_\Gamma(\pi)+\sum_{\pi^\prime\in \hat G_c}  n(\gamma,\pi^\prime)
N_\Gamma(\pi^\prime)\ ,$$
where the  sum reflects the fact that a finite number of representations
belonging to $\hat G_c$  
may contribute to the kernel of $D_Y$. Note that $n(\gamma,\pi^\prime)\in \Z$
and $N_\Gamma(\pi^\prime)\in\nat_0$. We show that $\res_{z=0} Q(z)=0$ in order
to conclude that $N_\Gamma(\pi)=-\sum_{\pi^\prime\in G_c}  n(\gamma,\pi^\prime)
N_\Gamma(\pi^\prime)\in \Z$.

It suffices to show that $Q(z)\equiv 0$ for $\Ree(z)\ll 0$.
This follows from (\ref{t5}), but we will give an independent argument.
We can write
$$r_Y(z)-r(z)=\frac{1}{z}((D_Y^2)_1 r_Y(z)-(D^2)_1 r(z))\ .$$
Integrating the restriction of this difference to the diagonal over
$\Gamma\backslash G$ we obtain
\begin{eqnarray*}
Q(z)&=& \frac{1}{z}\int_{\Gamma\backslash G} \tr\: ((D_Y^2)_1 r_Y(z)-(D^2)_1
 r(z))(g,g) dg\\
&=& \frac{1}{z}\int_{G}  \chi^\Gamma(g) \tr\: ((D_Y^2)_1 r_Y(z)-(D^2)_1
 r(z))(g,g) dg\\
&=& -  \frac{1}{z}\int_{G}  \chi^\Gamma(g) \tr\: ((D_Y)_1 (D_Y)_2 r_Y(z)-D_1
D_2 r(z))(g,g) dg\\&& -   \frac{1}{z}\int_{G}  \tr\:
c(d\chi^\Gamma)_1 ((D_Y)_1 r_Y(z)-D_1 r(z))(g,g) dg \\ 
&=&-  \frac{1}{z}\int_{G}    \chi^\Gamma(g) \tr\: ((D_Y^2)_1 
r_Y(z)-D_1^2  r(z))(g,g) dg\\&& -   \frac{1}{z}\int_{G}  \tr\:
c(d\chi^\Gamma)_1 ((D_Y)_1 r_Y(z)-D_1 r(z))(g,g) dg\ ,
\end{eqnarray*}
where $c(d\chi^\Gamma)$ denotes Clifford multiplication.
We conclude that
$$Q(z)=-   \frac{1}{2z}\int_{G}  \tr\:
c(d\chi^\Gamma)_1 ((D_Y)_1 r_Y(z)-D_1 r(z))(g,g) dg\ .$$
The right hand side of this equation vanishes as a consequence of 
$\sum_{\gamma\in\Gamma} \gamma^*\chi^\Gamma\equiv 1$ and the
$\Gamma$-invariance of $((D_Y)_1 r_Y(z)-D_1 r(z))(g,g)$.
This finishes the proof of the proposition. \hB

\bibliographystyle{plain}

\end{document}